\documentclass[12pt]{article}
\usepackage{amsfonts,amsmath,amsthm,amscd,amssymb,latexsym,mathrsfs}
\usepackage[mathscr]{eucal}
\usepackage[ruled,vlined]{algorithm2e}
\usepackage[all]{xy}
\usepackage{rotating}
\usepackage{lscape}
\usepackage{setspace}
\usepackage{tkz-graph}
\usepackage{tikz-cd}
\usetikzlibrary{calc}
\usepackage{tikz}
\usetikzlibrary{arrows}
\usepackage[T1]{fontenc}
\usepackage{float}

\usepackage{textcomp}
\usepackage{multicol}
\usepackage{enumerate}
\usepackage{harmony}
\usepackage{abc}
\usepackage{wasysym}
\usepackage{ textcomp}

\usepackage{picinpar}
\usepackage{longtable}
\usepackage{calc}
\usepackage{ifthen}
\usepackage{dcpic,pictexwd}
\usepackage[all]{xy}
\usepackage{color}
\usepackage{tabularx}
\usepackage{makeidx}

\usepackage{tikz-cd}
\usetikzlibrary{calc}
\usepackage{tikz}
\usetikzlibrary{arrows}
\usepackage{stackrel}
\usepackage{float}
\usepackage{picinpar}
\usepackage{longtable}
\usepackage{calc}
\usepackage{ifthen}
\usepackage{dcpic,pictexwd}
\usepackage[all]{xy}
\usepackage{color}
\usepackage{tabularx}
\usepackage{makeidx}

\usepackage{tikz-cd}
\usetikzlibrary{calc}
\usepackage{tikz}
\usetikzlibrary{arrows}
\usepackage{stackrel}
\usepackage{float}
\usepackage{textcomp}
\usepackage{tikz}
\usepackage{amssymb,latexsym}

    \usepackage[breakable]{tcolorbox}
    \usepackage{fancyvrb} 
    \usepackage[Export]{adjustbox} 
    \adjustboxset{max size={0.9\linewidth}{0.9\paperheight}}

\usepackage{mathtools}

\DeclareFontFamily{U}{musix}{}
\DeclareFontShape{U}{musix}{m}{n}{<-> s*[1.01] musix11}{}

\usepackage[minimal]{leadsheets}
\useleadsheetslibraries{musicsymbols}
\useleadsheetslibraries{chords}
\usepackage{harmony}

\date{}

\theoremstyle{plain}
\newtheorem{Theorem}{Theorem}

\newtheorem{Proposition}{Proposition}

\theoremstyle{Definition}

\newtheorem{Remark}{Remark}

\begin{document}

\textbf{Interactions Between Brauer Configuration Algebras and Classical Cryptanalysis to Analyze Bach's Canons}

\par\bigskip
Agust\'{\i}n Moreno Ca\~{n}adas\footnote{Universidad Nacional de Colombia (Bogot\'a)}, Pedro Fernando Fern\'andez Espinosa \footnote{Universidad de Caldas (Manizales)},  Jos\'e Gregorio Rodr\'{\i}guez Nieto \footnote{Universidad Nacional de Colombia (Medell\'{\i}n)}, Odette M. Mendez \footnote{Universidad Nacional de Colombia (Manizales)}, Ricardo Hugo Arteaga-Bastidas \footnote{Universidad Nacional de Colombia (Bogot\'a)}.

\par\bigskip

\abstract{Since their introduction, Brauer configuration algebras (BCAs) and their specialized messages have helped research in several fields of mathematics and sciences. This paper deals with a new perspective on using such algebras as a theoretical framework in classical cryptography and music theory. It is proved that some block cyphers define labeled Brauer configuration algebras. Particularly, the dimension of the  BCA associated with a ciphertext-only attack of the Vigenere cryptosystem is given by the corresponding key's length and the captured ciphertext's coincidence index. \par\bigskip On the other hand, historically, Bach's canons have been considered solved music puzzles. However, due to how Bach posed such canons, the question remains whether their solutions are only limited to musical issues. This paper gives alternative solutions based on the theory of Brauer configuration algebras to some of the puzzle canons proposed by Bach in his Musical Offering (BWV 1079) and the canon \^a 4 Voc: Perpetuus (BWV 1073). Specifically to the canon \^a 6 Voc  (BWV 1076), canon 1 \^a2 (also known as the crab canon), and canon \^a4 Quaerendo Invenietis. These solutions are obtained by interpreting such canons as ciphertexts (via route and transposition cyphers) of some specialized Brauer messages. In particular, it is noted that the structure or form of the notes used in such canons can be described via the shape of the most used symbols in Bach's works.}

\section{Introduction}\label{Introduction}


Brauer configuration algebras were introduced in 2017 to research algebras of wild representation type \cite{Green, Schroll}. Soon afterwards, it was discovered that such algebras are helpful tools for investigating different fields of mathematics and sciences due to their combinatorial nature. These algebras are induced by special systems of multisets, which are suitable tools to deal with the theory of integer partitions and enumerative combinatorics \cite{Stanley, Andrews}. Brauer messages associated with Brauer configurations have been used in cryptography to study the behavior of some provably secure hash functions or to give alternative descriptions to the schedule of an AES (Advanced Encryption Standard) key \cite{Canadas7, Canadas5}. This paper shows that behind the definition of some block cyphers is a Brauer configuration algebra whose invariants allow giving properties of the corresponding cryptographic system.

\par\bigskip

We will see that cryptanalysis of a Vigenere cryptosystem provides formulas for the dimension of its associated Brauer configuration algebra and its center. Similar results can be obtained for the Brauer configuration algebra associated with the permutation (transposition) cryptosystem. In fact, the dimension of the Brauer configuration algebra associated with a ciphertext obtained by permutation (or via a permutation network) equals the dimension of the Brauer configuration algebra associated with the corresponding plaintext.\par\bigskip

Musical compositions are another example of the presence of Brauer's messages. These musical contents can be seen as ciphertexts of specialized Brauer messages via an appropriated transposition.    \par\bigskip

In the Baroque period, musical canons were musical pieces proposed by composers as puzzles to be solved with the help of some predefined hints. Bach proposed some of the most celebrated canons in music history in his Musical Offering \cite{Rasmussen}. We remember that in this work, Bach proposed ten canons based on a single theme that King Frederick the Great proposed. Up-to-date, it is accepted that Bach's canons were solved by two of Bach's former students, Johann Friedrich Agricola and Johann Philip Kirnberger \cite{Hofstadter,  Milka, Harvard, Tatlow}. 

\par\bigskip
According to some of Bach's biographers \cite{Milka, Harvard, Shafer, Wissner, Sylvestre, Collins, Collins1, deCoul, Butler, Greer}, it is easy to think that many of the canons need to be solved taking into account another point of view because it is generally accepted that Bach used a particular type of symbolism in his music. We recall that Bach was used to including musical cryptograms in his works. For instance, his well-known motif is included in musical pieces such as the Art of Fugue or the Brandenburg Concerto No. 2 \cite{Barber}. In this line, this paper proposes to consider some of Bach's canons as ciphertexts obtained by the transposition of some Brauer messages in such a way that such messages allow to recover part of Bach's data and symbolism via a route deciphering.

\subsection{Motivations}

One of the main characteristics of Brauer configuration algebras is their capacity to be adapted in different contexts, provided its combinatorial definition \cite{Green}. The interaction of such algebras with different sciences and mathematics fields has allowed new perspectives to investigate the underlying theoretical frameworks \cite{Canadas7, Canadas5, Canadas1, Canadas10, Canadas2, Angarita, Canadas3, Canadas4, Canadas6}. \par\bigskip

This paper is focused on the interaction between classical cryptanalysis, music theory, and Brauer configuration algebras. We give formulas for the dimension of the Brauer configuration algebras induced by ciphertexts obtained via some block cyphers and music pieces interpreted as ciphertexts of some route and transposition cyphers. These procedures permit us to give alternative solutions to some of the canons proposed by Bach.

\subsection{Contributions}

The main results of this paper are Theorems  \ref{Permutation}, \ref{Vigenere1}, \ref{Vigenere2}, \ref{Musical} with Proposition \ref{Vigenere3}. Furthermore, the Brauer analysis of Bach canons \^a  6 Voc (BWV 1076), 1 \^a2 (crab canon), and \^a4 Quaerendo Invenietis realized in Section \ref{Bach}.  \par\bigskip

Theorem \ref{Permutation} proves that the Brauer configuration algebras associated with the plaintext and ciphertext of a permutation cryptosystem coincide. Theorems \ref{Vigenere1}, \ref{Vigenere2} and Proposition \ref{Vigenere3} give formulas for the dimensions of the Brauer configuration algebra and the corresponding center associated with a Vigenere ciphertext. \par\bigskip

Theorem \ref{Musical} gives conditions for a musical piece written on the usual Western staff to be a specialization of a suitable Brauer configuration. Section 3.5 is devoted to the Brauer analysis of some of Bach's canons. In particular, some commonly used symbols in Bach's works are built by representing musical notes as points in the plane. Such points are vertices of graphs induced by the canons. The graphs' edges are obtained by connecting consecutive classes of musical notes (i.e., consecutive points) appropriately.

\par\bigskip

The organization of this paper goes as follows: Background, main definitions and notation are given in section \ref{Background}; we remind definitions of multisets (Section \ref{OSSP}), Brauer configuration algebras (Section \ref{BCAs}) and some of their properties. Section \ref{Cryptography} is devoted to reminding basic facts about cryptography. We present the main results in section \ref{Main}. Concluding remarks are given in section \ref{conclusions}.

\section{Preliminaries}\label{Background}
This section provides basic definitions and notations regarding Brauer
configuration algebras, and cryptography \cite{Green, Schroll, Canadas7, Canadas5, Stinson, Fernandez, Rios, Canadas8, Sierra}.\par\bigskip

\subsection{Background}
 Brauer configuration algebras (BCAs) were introduced by Green and Schroll in 2017 \cite{Green, Schroll}. Soon afterwards, several works were written to apply these algebras in different fields of mathematics and sciences. On one hand, Espinosa and Rios wrote their doctoral dissertations based on these algebras \cite{Fernandez, Rios}. Espinosa et al. introduced the notion of a message and specialized message of a Brauer configuration to give a formula for the number of perfect matchings of a snake graph associated with some Kronecker modules \cite{Fernandez, Canadas4}. Rios \cite{Rios} used such algebras to describe particular classes of Dyck paths and integer friezes. This approach allows Dyck-Brauer messages to obtain alternative formulas for cluster variables associated with cluster algebras of type $\mathbb{A}_{n}$ \cite{Canadas8, Canadas9}.\par\bigskip
 
 Angarita et al. \cite{Canadas2} used BCAs to define alternative methods to protect biometric databases, and F\'uneme et al. \cite{Canadas5} introduced Brauer messages to examine the performance of some provably secure hash functions as those defined by Z\'emor and Tillich \cite{Tillich} and Sosnovski \cite{Sosnovski} based on the Cayley graph of some semigroups. \par\bigskip
 
 Relationships between BCAs and the graph energy theory were given by Agudelo et al. \cite{Canadas3}, who calculated the trace norm of some matrices induced by so-called $\{0,1\}$-Brauer configurations. It is worth noting that Espinosa \cite{Fernandez} also studied relationships between the graph energy theory and Brauer messages.\par\bigskip
 
 Ballester-Bolinches et al. \cite{Canadas1} used BCAs to give solutions to the Yang-Baxter equation. To do that, the authors specialized Brauer messages to define some skew braces in the sense of Rump \cite{Rump}. Espinosa and Ballester-Bolinches et al. \cite{Canadas10} also gave solutions to the Yang-Baxter equation via previous works written by Espinosa regarding the interaction between the Kronecker algebra, snake graphs theory and the theory of Brauer configuration algebras.\par\bigskip
 
Relationships between cryptography and Brauer configuration algebras were given by Marin et al. \cite{Canadas7}, who define the schedule of an AES (Advanced Standard Encryption) key in terms of appropriated mutations of Brauer messages. In this work, such messages allow solutions to some generalizations of the Chicken Mc Nugget problem, particularly the Frobenius problem.\par\bigskip

 On the other hand, it is worth pointing out that throughout history, Bach's work has been a source of a plethora of research, from the study of the symbology used to write his compositions to determining their corresponding fractal dimensions \cite{Rasmussen, Hofstadter, Sylvestre, Tatlow,Milka,Harvard, Shafer, Butler, Barber, Tatlow1, Madden,Hsu,  Tomita, Franklin, Wollny}. For instance, Smend  \cite{Madden} applied gematria to name the number 14 as the Bach number (BACH=2+1+3+8=14), and Tatlow discusses the use of cryptograms in Bach's work in \cite{Tatlow}. Furthermore, Niitsuma and Tomita studied Bach's writing in \cite{Tomita}. Bach journal, published by the Riemenschneider Bach Institute and Understanding Bach journal, edited by Ruth Tatlow, are devoted to Bach works in the Baroque era. Several papers published in these journals regard Bach's canons and his writing stylism, including his famous motif \cite{Wissner, deCoul, Collins1, Greer, Butler, Barber, Wollny, Hacohen, Emery}.\par\bigskip

 Sylvestre and Costa \cite{Sylvestre} and Shafer \cite{Shafer} studied some relationships between Bach's work, Fibonacci numbers and Fourier analysis.
 
 \par\bigskip
  It is worth noting, that the theory of graphs is a helpful tool to analyze symbolic music structure, which is an open problem in Music Information Retrieval (MIR) in this regard Hernandez-Olivan et al. \cite{Hernandez} used adjacency matrices of some graphs to segment symbolic music by its form or structure. According to them, the notes of a music piece $\mathfrak{M}$ can be seen as vertices of an appropriated graph $(V^{\mathfrak{M}},E^{\mathfrak{M}})$ and the set $E^{\mathfrak{M}}$ of edges between two of them is partitioned into three sets, i.e., $E^{\mathfrak{M}}=E^{on}\bigcup E^{cons}\bigcup E^{h}$, where $E^{on}$ consists of the edges connecting notes in the same onset, $E^{cons}$ consists of the edges connecting consecutive notes in time. And $E^{h}$ which is the set of edges linking overlapping notes in time. \par\bigskip
 
 On the other hand, Szeto et al. \cite{Szeto} defined posets (partially ordered sets) to facilitate pattern matching in post-tonal music analysis by using some pitch-class sets. In this work, the authors introduced the notion of stream, defined as the perceptual impression of connected series of musical notes. Whereas, stream segregation is the process of grouping musical notes into streams. According to them, the vertices or points of the corresponding Hasse diagram are given by events, and the main objective of the clustering algorithm is to connect an event to its nearest sequential event. 
 \par\bigskip
 We also remind that Jeong et al. \cite{Jeong} presented a graph neural network to learn note representation from music scores in Western notation.
 
 \par\bigskip
 
In this paper, we combine the techniques introduced by Hernandez Olivan et al. and Szeto et al. to construct appropriated graphs which can be read as symbols commonly used by Bach in his works.\par\bigskip

 In the sequel, we introduce some helpful notation and definitions regarding BCAs and classical cryptography.
  
  \subsection{ Multisets}\label{OSSP}

  A \textit{multiset} is a set with possibly repeated elements. Formally, a multiset is a pair of type $(M, f)$ where $M$ is a set and $f:M\rightarrow \mathbb{N}$ is a map from $M$ to the set of nonnegative integers. In this case, if $m\in M$ then $f(m)$ is said to be the \textit{multiplicity} of $m$ \cite{Andrews, Fontoura}.\par\bigskip
  
  A \textit{permutation} of a multiset $(M, f)$ is a word $w$ whose letters are given by the set $M$ whereas the number of occurrences of a given letter in $w$ is given by the map $f$. In this paper, we assume that the word $w(M)$ of a multiset $(M, f)$ with $M=\{m_{1},m_{2}\dots, m_{s}\}$ is given by a fixed permutation with the form
  \begin{equation}
  w(M)=m^{f(m_{1})}_{1}m^{f(m_{2})}_{2}\dots m^{f(m_{s})}_{s}.
  \end{equation}

  Note that there are $(f(m_{1})+f(m_{2})+\dots +f(m_{s}))!$ permutations associated with the multiset $(M,f)$.\par\bigskip
  
  If $(M_{1}, f_{1})$ and $(M_{2}, f_{2})$ are multisets with $M_{1}=\{m_{1,1}\dots, m_{1,r}\}$ and $M_{2}=\{m_{2,1},\dots, m_{2,s}\}$ then 
  
  \begin{equation}
  \begin{split}
  (M_{1}, f_{1})\cup (M_{2}, f_{2})&=(M_{1}\cup M_{2}, f_{1,2})\hspace{0.1cm}\text{with}\hspace{0.1cm}f_{1,2}(x)=\mathrm{max}\{f_{1}(x), f_{2}(x)\}.\\  
  (M_{1}, g_{1})\cap (M_{2}, g_{2})&=(M_{1}\cap M_{2}, g_{1,2})\hspace{0.1cm}\text{with}\hspace{0.1cm}g_{1,2}(x)=\mathrm{min}\{f_{1}(x), f_{2}(x)\}.\\  
  \end{split}
  \end{equation}
  
  Let $\mathscr{M}=\{(M_{1}, f_{1}), (M_{2}, f_{2}),\dots, (M_{h}, f_{h})\}$ be a collection of multisets such that if
  \par\bigskip
  \begin{centering}
   $w(M_{i})=(m_{i,1})^{f_{i}({m_{i,1}})}(m_{i,2})^{f_{i}({m_{i,2}})}\dots(m_{i, s_{i}})^{f_{i}({m_{i, s_{i}}})}$\hspace{0.1cm} then\hspace{0.1cm} $\underset{j=1}{\overset{s_{i}}{\sum}}f_{i}({m_{i,j}})>1$.\par\bigskip
\end{centering}  
  
In such a case. If $M= \underset{i=1}{\overset{h}{\bigcup}}M_{i}$, $ \underset{x\in I}{\bigcap}M_{x}$ is the interception of all multisets containing an element $y\in M$ with $I\subseteq\{1,2,\dots, h\}$ a fixed set of indices, and $f_{x}(y)$ is the frequency of $y$ in $M_{x}$ then $\underset{x\in I}{\sum} f_{x}(y)$  is said to be the \textit{valency} of $y$ denoted $val(y)$. We endow the set $\mathfrak{I}_{y}=\{M_{x}\mid x\in I\}$ with a linear order $<$. 
\par\bigskip
If  $x, x'\in {I}$, $z\in M_{x}\cap M_{x'}$, $z\notin \underset{M_{x}\in \mathfrak{I}_{y}}{\bigcap}M_{x} $, and $M_{x}<M_{x'}$ then the linear order $\prec$ associated with $\mathfrak{I}_{z}$ is defined in such a way that $M_{x}\prec M_{x'}$. Particularly, for $x\in I$ fixed and $f_{x}(y)>1$, there is a subchain  $\mathfrak{M}_{x, y}=M^{(1)}_{x}<M^{(2)}_{x}<\dots <M^{f_{x}( y)}_{x}$ where $M^{i}_{x}$ is associated with a unique copy of $M_{x}$ named the \textit{expansion} of $M_{x}$ induced by $y$.

\par\bigskip

 A collection of multisets $\mathscr{M}=\{(M_{1}, f_{1}),\dots, (M_{t}, f_{t})\}$ is said to be a \textit{system of multisets of type $M$}, if $M=\underset{j=1}{\overset{t}{\bigcup}}M_{j}$. And it is endowed with a map $\nu: M\rightarrow \mathbb{N}^{+}\times\mathbb{N}^{+}$ such that $\nu(m)=(j, val(m))$, for each $m\in M$. 
 
\par\bigskip 
Green and Schroll \cite{Green} named vertices the elements of $M$. The positive integer $j$ in a pair $(j, val(m))$ associated with a vertex $m$ is the multiplicity value $\mu(m)$ according to them. It is worth pointing out that $\mu(m)$ does not deal with the multiplicity of $m$ as an element of a multiset. In particular, $m\in M$ is said to be non-truncated (truncated) provided that $\mu(m)val(m)>1$ ($\mu(m)val(m)=1$). Thus, the multiplicity function $\mu$ classifies the set of vertices $M$ into the set of truncated and non-truncated vertices. In this work, if $val(m)=1$, it is assumed that $\mu(m)=2$. Therefore, the considered vertices are non-truncated.  
\par\bigskip

If $y\in M_{i_{1}}\cap M_{i_{2}}\cap\dots \cap M_{i_{t}}$ then a \textit{successor sequence} $\mathcal{S}_{y}$ associated with $y$ is a chain with the form

\begin{equation}\label{ss0}
M^{(1)}_{i_{1}}<\dots< M^{(f(i_{1},y))}_{i_{1}}<M^{(1)}_{i_{2}}<\dots <M^{(f(i_{2},y))}_{i_{2}}<\dots<M^{(1)}_{i_{t}}<\dots< M^{(f(i_{t},y))}_{i_{t}}   
\end{equation}

Successor sequence (\ref{ss0}) can be written in the form:

\begin{equation}\label{ess}
\mathfrak{M}_{i_{1},y}<\mathfrak{M}_{i_{2},y}<\dots<\mathfrak{M}_{i_{t},y}.
\end{equation}

  A \textit{Brauer configuration} is a system of multisets of type $M$ endowed with a function $\nu$ and an orientation $\mathcal{O}$ which is obtained by adding a relation $ M^{(f(i_{t},y))}_{i_{t}}<M^{(1)}_{i_{1}}$ ($M_{i_{t}}<M_{i_{1}}$) to each successor sequence $S_{y}$, associating in this fashion a family of equivalent circular orders to any $y\in M$. In this case, the multisets $(M_{i}, f_{i})$ are called \textit{polygons}. We let $\mathscr{M}=(M,\mathscr{M}_{1}, \nu, \mathcal{O})$ denote a Brauer configuration with $\mathscr{M}_{1}=\{(M_{1}, f_{1}),\dots, (M_{h}, f_{h})\}$, map $\nu$ and orientation $\mathcal{O}$ defined as above \cite{Green, Schroll}.
  
  \par\bigskip
   Note that the completed successor sequences or circular orders associated with a given vertex $y$ have the form
  \par\bigskip
  \begin{equation}\label{ss}
M^{(j)}_{i_{s}}<\dots< M^{(f(i_{s},y))}_{i_{s}}<\dots < M^{(f(i_{t},y))}_{i_{t}}<M^{(1)}_{i_{1}}<\dots< M^{(j-1)}_{i_{s}}.  
\end{equation}
  
 The \textit{message} $M(\mathscr{M})$ of a Brauer configuration is the concatenation of the words $w(M_{i})$, i.e.,
 \begin{equation}\label{message}
 M(\mathscr{M})=w(M_{1})w(M_{2})\dots w(M_{h}).
 \end{equation} 

  We say that a Brauer configuration $\mathscr{M}$ is $\mathfrak{S}$-\textit{labeled} if each multiset $(M_{i}, f_{i})$ is labeled by a permutation $\pi_{i}\in \mathfrak{S}_{f_{i}(m_{i, 1})+f_{i}(m_{i, 2})+\dots+f_{i}(m_{i, s_{i} })}$, where $\mathfrak{S}_{n}$ denotes the symmetric group with $n$ elements. Such a labeling means that the permutation $\pi_{i}$ is applied to the word $w(M_{i})$ for each $1\leq i\leq h$ \cite{Fernandez, Canadas7}. In such a case we write
  
\begin{equation}
M(\mathscr{M}^{\pi_{1},\dots, \pi_{h}})=w(M_{1}, \pi_{1})\dots w(M_{h}, \pi_{h})=w(\pi_{1}(M_{1}))\dots w(\pi_{h}(M_{h})).
\end{equation}

  \begin{Remark}\label{Label}
Since Brauer messages are defined by concatenating a fixed set of words given by a set of suitable permutations, we assume that all the Brauer configurations presented in this paper are $\mathfrak{S}$-labeled. i.e., polygons are defined by words obtained after the application of an appropriated permutation.
  \end{Remark}

The \textit{Brauer quiver} $Q_{\mathscr{M}}=(Q_{0}, Q_{1}, s, t)$ (or simply $Q$, if no confusion arises) induced by a Brauer configuration $\mathscr{M}=(M,\mathscr{M}_{1}, \nu, \mathcal{O})$ is defined in such a way that there is a bijective correspondence between its set of vertices $Q_{0}$ and the set of polygons $\mathscr{M}_{1}$. In such a case, each covering $(M_{i}, f_{i})<(M_{j}, f_{j})$ in a circular ordering defines an arrow from $(M_{i}, f_{i})$ to $(M_{j}, f_{j})$ in $Q_{1}$.\par\bigskip

    \subsubsection{Brauer Configuration Algebras}\label{BCAs}

 A Brauer configuration algebra $\Lambda_{\mathscr{M}}$ (or simply $\Lambda$ if no confusion arises) is a bound quiver algebra $\Lambda_{\mathscr{M}}=kQ_{\mathscr{M}}/\langle \rho\rangle$ induced by a Brauer quiver $Q_{\mathscr{M}}$ bounded by an ideal $I=\langle \rho\rangle$ generated by a set of relations $\rho$ of the following types \cite{Green, Schroll}:
 
 \begin{itemize}
 \item $C^{\mu(i)}_{i}-C^{\mu(j)}_{j}$ if $i$ and $j$ are vertices in the same polygon $(M_{i}, f_{i})$ and $C_{x}$ is a cycle defined by circular relations. These cycles are said to be \textit{special cycles}.
 \item $C^{\mu(i)}a$ if $a$ is the first arrow of a special cycle $C_{i}$ associated with a vertex $i$.
 \item $ab$ if $a,b\in Q_{1}$ are arrows of different special cycles and $ab$ is an element of the path algebra induced by $Q_{\mathscr{M}}$. 
 \end{itemize}   
 
 \begin{Remark}\label{Green}
 Green and Schroll \cite{Green} proved that Brauer configuration algebras are symmetric and multiserial and that there is a bijective correspondence between the indecomposable projective $\Lambda_{\mathscr{M}}$-modules and the polygons in $\mathscr{M}_{1}$. Furthermore, the number of summands in the heart of an indecomposable projective $\Lambda_{\mathscr{M}}$-module $P$ with $rad^{2}\hspace{0.1cm}P\neq0$ equals the number of non-truncated vertices in the corresponding polygon. Particularly, it holds that
 
\begin{equation}\label{algebra}
\mathrm{dim}_{k}\hspace{0.1cm}\Lambda_{\mathscr{M}}=2|\mathscr{M}_{1}|+\underset{m\in M}{\sum}val(m)(val(m)\mu(m)-1).
\end{equation}
 Sierra \cite{Sierra} obtained the following formula for the dimension of the center $Z (\Lambda_{\mathscr{M}})$ of a connected Brauer configuration algebra.
  
  \begin{equation}\label{center}
 \mathrm{dim}_{k}\hspace{0.1cm}Z(\Lambda_{\mathscr{M}})=1+|\mathscr{M}_{1}|-|M|+\underset{m\in M}{\sum}\mu(m)+\#(Loops (Q_{\mathscr{M}})) -|\{m\in M\mid val(m)=1\}|.  
  \end{equation}

\end{Remark}

\subsection{Cryptography}\label{Cryptography}
 
 This section reminds the definition of well-known cryptosystems \cite{Stinson}. Some of them related to the theory of Brauer configuration algebras.\par\bigskip 

A cryptosystem or cryptographic system is a quintuple of the form $(\mathcal{P},\mathcal{C},\mathcal{K},\mathcal{E}, \mathcal{D})$, where 

\begin{itemize}
\item $\mathcal{P}$ is a finite set of possible plaintexts.
\item $\mathcal{C}$ is a finite set of possible ciphertexts.
\item $\mathcal{K}$ is a finite set of possible  keys.
\item $\mathcal{E}$ ($\mathcal{D}$) is a set of encryption (decryption) functions such that for each $K\in\mathcal{K}$, it holds that $e_{K}:\mathcal{P}\rightarrow \mathcal{C}\in \mathcal{E}$ ( $d_{K}:\mathcal{C}\rightarrow \mathcal{P}\in \mathcal{D}$), and $d_{K}(e_{K})(x)=x$, for each $x\in\mathcal{P}$.

\end{itemize}

The permutation cryptosystem and the Vigenere cryptosystem are examples of classical cryptosystems associated with Brauer configuration algebras. The following section will be devoted to clarifying these relationships. \par\bigskip

The permutation or transposition cryptosystem is defined in such a way that for fixed integers $n, m>1$, it holds that $\mathcal{P}=\mathcal{C}=\mathbb{Z}^{m}_{n}$, in this case $\mathcal{K}=\mathfrak{S}_{m}$ is the set of all $m$-element permutations. For $\pi\in\mathfrak{S}_{m}$ and $x=(x_{1}, x_{2},\dots, x_{m})\in \mathcal{P}$, $e_{\pi}(x)=(x_{\pi(1)}, x_{\pi(2)},\dots, x_{\pi(m)})$. Whereas, for $y=(y_{1}, y_{2},\dots, y_{m})\in\mathcal{C}$, $d_{\pi}(y_{1}, y_{2},\dots, y_{m})=(y_{\pi^{-1}(1)}, y_{\pi^{-1}(2)},\dots, y_{\pi^{-1}(m)})$. The following is an example of this kind of encryption for $\mathcal{P}=\mathcal{C}=\mathbb{Z}^{4}_{26}$, and $\pi= (3\hspace{0.2cm}  4\hspace{0.2cm}  1\hspace{0.2cm}  2)$.\par\smallskip

\begin{centering}
\begin{equation}\label{examplecr}
\begin{array}{|c|c|c|}\hline
Y & O & H \\\hline
P & T & Y \\\hline
C & R & A\\\hline
R & G & P\\\hline
\end{array}
\end{equation}\par\smallskip
\end{centering}

Applying $\pi^{-1}$ to each block, it is obtained the following array:

\begin{centering}
\begin{equation}\label{decrypt}
\begin{array}{|c|c|c|}\hline
C & R & A \\\hline
R & G & P \\\hline
Y & O & H\\\hline
P & T & Y\\\hline
\end{array}
\end{equation}\par\smallskip
\end{centering}

\begin{figure}[h]
	\begin{equation}
		\xymatrix@=30pt{
			\bullet\ar@[blue][ddd]&\bullet\ar@[blue][r]&\bullet\ar@[blue][ddd]&\\
			\bullet&\bullet&\bullet&\\
			\bullet&\bullet&\bullet&\\
			\bullet\ar@[blue][r]&\bullet\ar@[blue][uuu]&\bullet&\\
		}\notag
	\end{equation}
	\caption{The route (path) used to decrypt the message (\ref{decrypt}) obtained via permutation.}\label{hx1}
\end{figure}
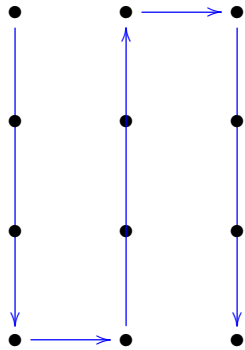

The path shown in figure \ref{hx1} gives rise to the plaintext CRYPTOGRAPHY. So, it is considered as a route cypher. \par\bigskip

\par\bigskip

For fixed integers $m, n>1$ fixed, $\mathbb{Z}^{m}_{n}=\mathcal{P}=\mathcal{C}=\mathcal{K}$ in the Vigenere cryptosystem. In such a case, for a given key $K=(k_{1}, k_{2},\dots, k_{m})\in\mathcal{K}$, $x=(x_{1}, x_{2},\dots, x_{m})\in\mathcal{P}$, and $y=(y_{1}, y_{2},\dots, y_{m})\in\mathcal{C}$, it holds that

\begin{equation}
\begin{split}
e_{K}(x)&=(x_{1}+k_{1}, x_{2}+k_{2}, \dots, x_{m}+k_{m})\mod{n}\\
d_{K}(y)&=(y_{1}-k_{1}, y_{2}-k_{2}, \dots, y_{m}-k_{m})\mod{n}
\end{split}
\end{equation}

The Vigenere cryptosystem is vulnerable to a ciphertext-only attack introduced by Friedman \cite{Stinson}. It is based on the coincidence index of a text. In such a case, it is assumed that $n$ is the size of the underlying alphabet. For instance, $n=26$ if the language is English. The \textit{index of coincidence} $\mathscr{I}(T)$ of a string $T$ of alphabetic characters (which is the probability that two random elements of $T$ be identical) is given by the following formula:

\begin{equation}
\mathscr{I}(T)=\underset{i=0}{\overset{25}{\sum}}\frac{f_{i}(f_{i}-1)}{n(n-1)}.
\end{equation}  

Where $n=|T|$ and for each $i$, $f_{i}$ denotes the frequency of the $i$th character in $T$ ($f_{0}$ for A, $f_{1}$ for B, and so on).\par\bigskip

Suppose $T$ and $T'$ are strings of alphabetic characters. In that case, the \textit{mutual index of coincidence} (which is the probability that a random element of $T$ is identical to an element of $T'$) is given by the following formula.

\begin{equation}
M\mathscr{I}(T, T')=\underset{i=0}{\overset{25}{\sum}}\frac{f_{i}f'_{i}}{nn'}.
\end{equation}

Where $n=|T|$, $n'=|T'|$ and $f_{i}$ ($f'_{i}$) denotes the frequency of the $i$th character in $T$ ($T'$).\par\bigskip

The coincidence index attack is based on the fact that modulo the length of the key, identical characters are encrypted by the same alphabetic character. So, if $y$ is a ciphertext obtained by a Vigenere key $k=(k_{1}, k_{2}, \dots, k_{m})$ of length $m$. Then $y$ can be split into $m$ lists $y_{1}, y_{2},\dots, y_{m}$ whose indices are all close to 0.065. \par\bigskip

Since the lists $y_{1}, y_{2},\dots, y_{m}$ are obtained by applying an appropriated shift $\mathfrak{s}$ and computing the different mutual indices of coincidence $M\mathscr{I}(y_{i}, y^{\mathfrak{s}}_{j})=\frac{\overset{25}{\sum}f_{i}f'_{i-\mathfrak{s}}}{nn'}$. Then a set of equations with the form $k_{i}-k_{j}=\mathfrak{s}_{0}$ is built for each of these comparisons, where $\mathfrak{s}_{0}$ is a value for which $M\mathscr{I}(y_{i}, y^{\mathfrak{s}_{0}}_{j})$ is close to 0.065. The solution of the obtained system of linear equations (modulo 26) gives the Vigenere key used to obtain the ciphertext $y$. Particularly, if $m$ is the length of the key then $\mathscr{I}(y_{i})\sim$ 0.065 for any $1\leq i\leq m$.

\section{Main Results}\label{Main}

This section gives the main results regarding the  dimension of Brauer configurations algebras associated with  polyalphabetic cryptosystems, mutations of Brauer configurations, and iterated cyphers. Subsections \ref{Content} and \ref{Bach} are devoted to the Brauer analysis of musical content, subsection \ref{Bach} analyzes Bach's canons.

\subsection{Brauer Configuration Algebras Associated with the Permutation Cryptosystem}

Suppose that $y=y_{1}y_{2}\dots y_{m}$ is a ciphertext obtained by encrypting  a plaintext $x$ divided into $m$ blocks $x_{1}, x_{2},\dots, x_{m}$ not necessarily of the same size, i.e., $x=x_{1}x_{2}\dots x_{m}$, with $|x_{i}|=s_{i}=|y_{i}|>1$. The encryption and decryption rules are given by the following identities:

\begin{equation}
\begin{split}
e_{K}(x)&=e_{\pi_{1}}(x_{1})e_{\pi_{2}}(x_{2})\dots e_{\pi_{m}}(x_{m})=\pi_{1}(x_{1})\pi_{2}(x_{2})\dots \pi_{m}(x_{m}).\\
d_{K}(y)&=d_{\pi^{-1}_{1}}(x_{1})d_{\pi^{-1}_{2}}(x_{2})\dots d_{\pi^{-1}_{m}}(x_{m})=\pi^{-1}_{1}(x_{1})\pi^{-1}_{2}(x_{2})\dots \pi^{-1}_{m}(x_{m}).
\end{split}
\end{equation}
 where $\pi_{i}\in\mathfrak{S}_{i}$ for $1\leq i\leq m$.\par\bigskip

Plaintexts and ciphertexts in the transposition cryptosystem define labeled Brauer configurations bearing in mind that $x$ can be seen as a Brauer message of a labeled Brauer configuration $\mathscr{M}$ of the form

\begin{equation}
\begin{split}
x&=M(\mathscr{M})=(w(x_{1}), \sigma_{1})(w(x_{2}), \sigma_{2})\dots(w(x_{m}), \sigma_{m}),\quad \sigma_{i}\in\mathfrak{S}_{i}\quad 1\leq i\leq m.\\
x_{i}&=\sigma_{i}(w(x_{i})).
\end{split}
\end{equation}

Where for each $i$, $w(x_{i})=x^{f_{i_{1}}}_{i, 1}x^{f_{i_{2}}}_{i, 2}\dots x^{f_{i_{m_{i}}}}_{i, m_{i}}$ ($f_{i_{1}}+f_{i_{2}}+\dots+f_{i_{m_{i}}}=s_{i})$ is a word determining the polygon $U_{i}\in \mathscr{M}_{1}=\{U_{1}, U_{2},\dots, U_{m}\}$. Note that, $x_{i, j}$ is a character of a fixed alphabet $\mathscr{A}$. Thus,
$\mathscr{M}=(M, \mathscr{M}_{1}, \mu, \mathcal{O})$, where

\begin{itemize}
\item $M=\{x_{i, j}\mid 1\leq i\leq m$, $1\leq j\leq s_{i}\}$.
\item $\mathscr{M}_{1}=\{w(x_{i})\mid 1\leq i\leq m\}$.
\item $\mu(x_{i, j})=1$ ($\mu(x_{i, j})=2$) if $val(x_{i, j})>1$ ($val(x_{i,j})=1$).
\item $(w(x_{1}),\sigma_{1})<(w(x_{2}) ,\sigma_{2})<\dots<(w(x_{m}), \sigma_{m})$ in successor sequences.
\end{itemize}

In the same fashion the ciphertext $y=\pi_{1}(x_{1})\pi_{2}(x_{2})\dots \pi_{m}(x_{m})$ is a labeled Brauer message $M(\mathscr{M}^{\pi_{1},\pi_{2},\dots, \pi_{m}})$ obtained by permuting vertices in $x_{i}$ for $1\leq i\leq m$. We let $\Lambda_{\mathscr{M}}$ ($\Lambda_{\mathscr{M}^{\pi_{1},\dots, \pi_{m}}}$) denote the labeled  Brauer configuration algebra defined by the plaintext $x$ (ciphertext $y$). \par\bigskip

The following result holds:

\begin{Theorem}\label{Permutation}
$\Lambda_{\mathscr{M}}=\Lambda_{\mathscr{M}^{\pi_{1},\dots, \pi_{m}}}$.

\end{Theorem}

\textbf{Proof.} Note that the system of permutations $\pi_{1},\pi_{2},\dots, \pi_{m}$ does not change vertices, polygons nor the orientation of the polygons in successor sequences defining the Brauer configuration $\mathscr{M}$. \hspace{0.5cm}$\square$

\subsection{Brauer configuration Algebras Associated with the Cryptanalysis of the Vigenere Cryptosystem}

Suppose that $y$ is a captured ciphertext of a Vigenere cryptosystem. We can assume that $y$ contains all the characters of an alphabet $\mathscr{A}$ of a natural language $\mathscr{L}$ and the frequency $f_{j}$ of any $j\in\mathscr{A}$ is greatest than 1. If $m$ is the length of the key $K=(k_{1}, k_{2},\dots, k_{m})$ then $y$ is a concatenation of $m$ lists $y_{1}, y_{2},\dots, y_{m}$. Thus $y$ can be seen as a Brauer message $M(\mathfrak{C})$ such that:

\begin{itemize}
\item $\mathfrak{C}_{0}=\mathscr{A}$.
\item $\mathfrak{C}_{1}=\{y_{1}, y_{2},\dots, y_{m}\}$.
\item $\mu(i)=1$, for any $i\in\mathscr{A}$.
\item $y_{1}<y_{2}<\dots<y_{m}$ in successor sequences.
\end{itemize}

We have the following result:

\begin{Theorem}\label{Vigenere1}
Let $\Lambda_{\mathfrak{C}}$ be a Brauer configuration algebra induced by a Vigenere ciphertext $y$ obtained with a key of length $m$. Then 
\begin{equation}
\mathrm{dim}_{k}\hspace{0.1cm}\Lambda_{\mathfrak{C}}=2m+|y|(|y|-1)\mathscr{I}(y).
\end{equation}

\end{Theorem}

\textbf{Proof.} Since the size of the key is $m$. Then $|\mathfrak{C}_{1}|=m$ and $\mathscr{I}(y)=\underset{i\in\mathscr{A}}{\sum}\frac{f_{i}(f_{i}-1)}{|y|(|y|-1|)}$.\hspace{0.5cm}$\square$
\par\bigskip

The following result regards the dimension $\mathrm{dim}_{k}\hspace{0.1cm}Z(\Lambda_{\mathfrak{C}})$ of the center of the algebra $\Lambda_{\mathfrak{C}}$ induced by a Vigenere ciphertext $y$.

\begin{Theorem}\label{Vigenere2}
Let $\Lambda_{\mathfrak{C}}$ be the Brauer configuration algebra induced by a Vigenere ciphertext $y$ and $y_{1}, y_{2},\dots y_{m}$ are its corresponding lists obtained with a key of length $m$. Then if $f_{i ,j}$ denotes the frequency of the $j$th alphabetic character in the list $y_{i}$, it holds that
\begin{equation}
\mathrm{dim}_{k}\hspace{0.1cm}Z(\Lambda_{\mathfrak{C}})=1+m+\underset{i=1}{\overset{m}\sum}\underset{j\in\mathscr{A}}{\sum}(f_{i ,j}-1).
\end{equation}

\end{Theorem}

\textbf{Proof.} Since the length of the key is $m$ then $|\mathfrak{C}_{1}|=m$. Furthermore, the number of loops associated with the $j$th character-vertex of the alphabet $\mathscr{A}$ in the list-polygon $y_{i}$  is $f_{i, j}-1$.\hspace{0.5cm}$\square$

\par\bigskip

As an example, the Vigenere ciphertext 

\begin{equation}
{C}=\mathrm{OOPAELRIXFGGBWDODDEPK}
\end{equation}

is obtained from the plaintext

\begin{equation}
\mathrm{classicalcryptography}
\end{equation}
 with the key
 
 \begin{equation}
 \mathrm{MDPI}
 \end{equation}

This scheme induces the Brauer configuration $\mathfrak{C}=(\mathfrak{C}_{0},\mathfrak{C}_{1},\mu,\mathcal{O})$ such that

\begin{itemize}
\item $\mathfrak{C}_{0}=\{\mathrm{O,P,A,E,L,R,I,X,F,G,B,W,D,K}\}$.
\item$\mathfrak{C}_{1}=\{y_{1}=\{\mathrm{O,E,X,B,D,K}\}, y_{2}=\{\mathrm{O,L,F,W,D}\}, y_{3}=\{\mathrm{P,R,G,D,E}\}, y_{4}=\{\mathrm{A,I,G,O,P}\}\}$.
\item $\mu(O)=\mu(P)=\mu(E)=\mu(G)=\mu(D)=1$, \quad $\mu(j)=2$, for the remaining vertices $j\in M$.
\item $y_{1}<y_{2}<y_{3}<y_{4}$ in successor sequences.
\end{itemize}

\begin{equation}
\begin{split}
S_{\mathrm{O}}&=y_{1}<y_{2}<y_{4},\\
S_{\mathrm{P}}&=y_{3}<y_{4},\\
S_{\mathrm{A}}&=y_{4},\\
S_{\mathrm{E}}&=y_{1}<y_{3},\\
S_{\mathrm{L}}&=y_{2},\\
S_{\mathrm{R}}&=y_{3},\\
S_{\mathrm{I}}&=y_{4},\\
S_{\mathrm{X}}&=y_{1},\\
S_{\mathrm{F}}&=y_{2},\\
S_{\mathrm{G}}&=y_{3}<y_{4},\\
S_{\mathrm{B}}&=y_{1},\\
S_{\mathrm{W}}&=y_{2},\\
S_{\mathrm{D}}&=y_{1}<y_{2}<y_{3},\\
S_{\mathrm{K}}&=y_{1}.
\end{split}
\end{equation}

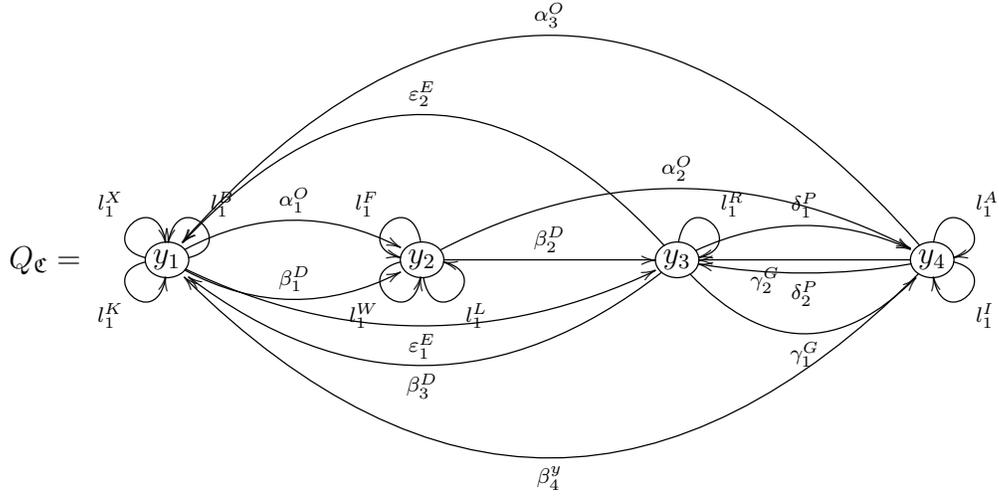
\begin{figure}[H]
 \begin{equation}
Q_{\mathfrak{C}}=\xymatrix@=80pt{*+[o][F-]{y_{1}}\ar@/^15pt/@[\blue][r]^{\alpha^{O}_{1}}\ar@/_15pt/@[\blue][r]^{\beta^{D}_{1}}\ar@(l,d)@[\blue]_{l_1^K}\ar@(l,u)@[\blue]^{l_1^X}\ar@(r,u)@[\blue]_{l_1^B}\ar@/_25pt/@[\blue][rr]_{\varepsilon^{E}_{1}}&*+[o][F-]{y_{2}}\ar@/^27pt/@[\blue][rr]^{\alpha^{O}_{2}}\ar@/_0pt/@[\blue][r]^{\beta^{D}_{2}}\ar@(u,l)@[\blue]_{l_1^F}\ar@(l,d)@[\blue]_{l_1^W}\ar@(d,r)@[\blue]_{l_1^L}&*+[o][F-]{y_{3}}\ar@/^40pt/@[\blue][ll]^{\beta^{D}_{3}}\ar@/_28pt/@[\blue][r]_{\gamma^{G}_{1}}\ar@(u,r)@[\blue]^{l_1^R}\ar@/^13pt/@[\blue][r]^{\delta^{P}_{1}}
\ar@/_55pt/@[\blue][ll]_{\varepsilon^{E}_{2}}&*+[o][F-]{y_{4}}\ar@(u,r)@[\blue]^{l_1^A}\ar@(r,d)@[\blue]^{l_1^I}\ar@/_85pt/@[\blue][lll]_{\alpha^{O}_{3}}\ar@/^75pt/@[\red][lll]^{\beta^{y}_{4}}\ar@/15pt/@[\blue][l]^{\gamma^{G}_{2}}\ar@/^5pt/@[\blue][l]^{\delta^{P}_{2}}&
}\notag
\end{equation}

\caption{Brauer quiver of a Vigenere ciphertext.}\label{mutfinal0}
\end{figure}

 The following relations define the admissible ideal $I_{\mathfrak{C}}$ with $M=\{\alpha^{j_{1}}_{i_{1}},\beta^{j_{2}}_{i_{2}},\gamma^{j_{3}}_{i_{3}},\delta^{j_{4}}_{i_{4}},\varepsilon^{j_{5}}_{i_{5}}\}$.
 
 \begin{itemize}
\item  $m^{{j_{s}}}_{i_{s}}n^{{j_{t}}}_{i_{t}}$, \quad $m^{{j_{s}}}_{i_{s}}\neq n^{{j_{t}}}_{i_{t}}$,\quad $m^{{j_{s}}}_{i_{s}},n^{{j_{t}}}_{i_{t}}\in M$.
\item $(l^{j}_{1})^{2}$, for all the possible values of $j$.
\item $\alpha^{O}_{1}\alpha^{O}_{2}\alpha^{O}_{3}=\beta^{D}_{1}\beta^{D}_{2}\beta^{D}_{3}=\varepsilon^{E}_{1}\varepsilon^{E}_{2}$.
\item $\alpha^{O}_{2}\alpha^{O}_{3}\alpha^{O}_{1}=\beta^{D}_{2}\beta^{D}_{3}\beta^{D}_{1}$.
\item $\gamma^{G}_{1}\gamma^{G}_{2}=\beta^{D}_{3}\beta^{D}_{1}\beta^{D}_{2}=\varepsilon^{E}_{2}\varepsilon^{E}_{1}=\delta^{P}_{1}\delta^{P}_{2}$.
\item $\gamma^{G}_{2}\gamma^{G}_{1}=\alpha^{O}_{3}\alpha^{O}_{1}\alpha^{O}_{2}=\delta^{P}_{2}\delta^{P}_{1}$.
\item $C^{j}_{i}f$, if $C^{j}_{i}$ is a special cycle and $f$ is its first arrow. 
 \end{itemize}
 It holds that
 
 \begin{equation}
 \begin{split}
 \mathrm{dim}_{k}\hspace{0.1cm}\Lambda_{\mathfrak{C}}&=35,\\
  \mathrm{dim}_{k}\hspace{0.1cm}Z(\Lambda_{\mathfrak{C}})&=1+18+5-14+4+9-9=14,\\
  \mathscr{I}(C)&=\frac{ \mathrm{dim}_{k}\hspace{0.1cm}\Lambda_{\mathfrak{C}}-2(4)}{21(20)}=\frac{27}{420}=0,0642. 
 \end{split}
 \end{equation}

\begin{Proposition}\label{Vigenere3}
Let $\Lambda_{\mathscr{M}}$ be a Brauer configuration algebra induced by a Brauer configuration $\mathscr{M}=(M,\mathscr{M}_{1},\mu,\mathcal{O})$ with
$\mathscr{M}_{1}=\{(U_{1}, f_{1}),\dots, (U_{m}, f_{m})\}$ and for any  $u_{i, j}\in U_{i}$, $1\leq i\leq m$ it holds that $f_{i}({u_{i, j}})=1$. And there are $n$ vertices $\alpha$ such that  $val(\alpha)=1$. Then

\begin{equation}
\mathrm{dim}_{k}\hspace{0.1cm}Z(\Lambda_{\mathscr{M}})=m+n+1.
\end{equation}

\end{Proposition}

\textbf{Proof.} Note that $\#(Loops (Q_{\mathscr{M}}))=|\{i\in M\mid val(i)=1\}|$, and $\underset{i\in M}{\sum}\mu(i)=2n+|M|-n=n+|M|$. Thus, $\mathrm{dim}_{k}\hspace{0.1cm}Z(\Lambda_{\mathscr{M}})=1+m+n+|M|-|M|=1+m+n$. We are done.\hspace{0.5cm}$\square$

\subsection{Brauer Configuration Algebras Associated with Musical Content}\label{Content}

This section proves that some musical contents give rise to Brauer configuration algebras.\par\bigskip

Let $\mathscr{A}=\mathscr{X}\bigcup\mathscr{S}_{0}$ be an alphabet of a language whose letters or vertices are partitioned into two sets $\mathscr{X}$ and $\mathscr{S}_{0}$, $\mathscr{X}$ is endowed with a partial order $\preceq$ and a subset of real numbers  $\mathscr{T}=\{t_{1},t_{2},\dots, t_{k}\}$. Elements $s\in\mathscr{S}_{0}$ are said to be constants. Each element $x_{i}\in\mathscr{X}$ is labeled by a unique element $t_{i}\in \mathscr{T}$. Therefore, each element $x_{i}\in \mathscr{X}$ can be seen as an ordered pair of the form $(x_{i}, t_{i})$ and

\begin{equation}
\mathscr{X}=\{(x_{1},t_{1}),(x_{2}, t_{2}),\dots, (x_{k}, t_{k})\}
\end{equation}

$(x_{i}, t_{i})\unlhd (x_{j}, t_{j})$ if and only if  $x_{i}\preceq x_{j}$ or $t_{i}\leq t_{j}$ (the usual real numbers order). \par\bigskip

In this case,

\begin{equation}
\mathrm{min}(\mathscr{X})=x_{1}\unlhd x_{2}\unlhd \dots\unlhd x_{k-1}\unlhd x_{k}=\mathrm{max}(\mathscr{X})
\end{equation}

$\mathscr{X}$ is endowed with a circular order of the form\par\bigskip

$\mathfrak{c}_{h}=(x_{h}, t_{h})\unlhd (x_{h+1}, t_{h+1})\unlhd \dots\unlhd(x_{k}, t_{k})\unlhd (x_{1}, t_{1})\unlhd \dots\unlhd (x_{h-1}, t_{h-1})\unlhd(x_{h}, t_{h})$, $1\leq h\leq k$.\par\bigskip

Numbers $t_{j}$ in pairs $(x_{j},t_{j})$ are given by coverings $x_{j-1}\preceq x_{j}$, i.e.,

\begin{equation}
x_{j-1}\preceq x_{j}\longrightarrow t_{j}.
\end{equation}

If numbers $t_{j}$ are predefined, then (if no confusion arises), we will omit their use in the notation of the pairs constituting a set $\mathscr{X}$.\par\bigskip

The words $\mathfrak{w}\in\mathscr{X}^{*}$ defined by the set $\mathscr{X}$ have the form:

\begin{equation}\label{word}
\mathfrak{w}=\mathfrak{z}^{F_{1}}_{1}\mathfrak{z}^{F_{2}}_{2}\dots\mathfrak{z}^{F_{s}}_{s}.
\end{equation}
\par\bigskip

For each $1\leq i\leq s $, $\mathfrak{z}_{i}\in\mathscr{X}_{j}\subset\mathscr{X}^{*}$ for some integer $j\geq 1$. Where, $\mathscr{X}_{j}$ consists of length $j$ words with the shape:

\begin{equation}
(x_{i_{h,1}},t_{i_{h,1}} )(x_{i_{h,2}}, t_{i_{h,2}})\dots(x_{i_{h, j}}, t_{i_{h, j}}), \quad 1\leq h\leq n_{j}.
\end{equation}
To define operators, we will assume the notation $(x_{i_{h,j}},t_{i_{h,j}})=\mathfrak{y}_{i_{h, j}}$.
\par\bigskip

We let $F_{i}$ denote an integral vector $F_{i}=(f_{i,1}, f_{i,2},\dots, f_{i,j})$ such that

\begin{equation}
\mathfrak{z}^{F_{i}}_{i}=(x_{i_{h,1}},t_{i_{h,1}} )^{f_{i,1}}(x_{i_{h,2}}, t_{i_{h,2}})^{f_{i,2}}\dots(x_{i_{h, j}}, t_{i_{h, j}})^{f_{i,j}}.
\end{equation}

Where $(x_{i_{h,1}},t_{i_{h,1}} )^{f_{i,1}}=\underset{f_{i,1}-\mathrm{times}}{\underbrace{{(x_{i_{h,1}},t_{i_{h,1}})\dots(x_{i_{h,1}},t_{i_{h,1}} )}}}$.

\begin{center}
\textbf{Operators}

\end{center}

Two maps $\sigma_{r_{\mathfrak{y}_{i_{h, j}}}}$ and $\varphi_{h}$ are defined in such a way that if $\mathfrak{w}$ is given as in (\ref{word}) then

 \[\sigma_{r_{\mathfrak{y}_{i_{h, j}}}}(\mathfrak{y}_{i_{h', j'}})=
\begin{cases}
(x_{i_{h,1}},t_{i_{h,1}}+\frac{(t_{i_{h,1}})(r_{\mathfrak{y}_{i_{h, j}}})}{2}),  & \hspace{0.2cm}\text{if}\hspace{0.2cm}h=h'\quad\text{and}\quad j=j', \\
\mathfrak{y}_{i_{h', j'}},  & \hspace{0.2cm}\text{otherwise.}
\end{cases}\]
\par\bigskip
${r_{\mathfrak{y}_{i_{h, j}}}}\in\mathbb{R}$. 

\par\bigskip
\begin{equation}
\varphi_{h}((x_{i_{p, j}}, t_{i_{p,j}})^{f_{p,j}})=(x_{i_{p, j}}, t_{i_{p ,j}})^{f_{p, j}+hf_{p, j}}
\end{equation}

\par\bigskip

\begin{equation}
\begin{split}
\sigma_{2}(x_{i-1}, t_{i-1})&=(x_{i}, 2t_{i-1})=(x_{i}, t_{i})\\
\sigma_{-2}(x_{i}, t_{i})&=(x_{i-1}, 0)=(x_{i-1}, t_{i-1})\\
\sigma_{x}(\sigma_{-x})((x_{i}, t_{i}))&=(x_{i}, t_{i}).
\end{split}
\end{equation}

Note that in the cycle $\mathfrak{c}_{h}$, it holds that $(x_{h+k}, t_{h+k})=\sigma_{2k}(x_{h}, t_{h})$.

\begin{center}
\textbf{Grouping Symbols}
\end{center}

Letters within a word $\mathfrak{w}$  or in different words can be grouped using parentheses, brackets and braces. Often, parentheses are used to group letters with the same coordinates. Whereas brackets and braces group letters within a word with the same frequency.

\begin{equation}\label{word1}
\underset{n-\text{times}}{{(\underbrace{(x_{1}, t_{1})^{f_{1}}(x_{1}, t_{1})^{f_{1}}\dots(x_{1}, t_{1})^{f_{1}}}}(x_{1}, t_{1})^{f_{i}})}=(x_{1}, t_{1})^{\underset{n-\text{veces}}{\underbrace{f_{1}+f_{1}+\dots+f_{1}}+f_{i}}}\
\end{equation}

\begin{equation}\label{word2}
\mathfrak{w}=\underset{w_{1}}{\underbrace{(x_{i_{1}}, t_{i_{1}})^{f_{i_{1}}}(x_{i_{2}}, t_{i_{2}})^{f_{i_{2}}}\dots((x_{i_{s}}, t_{i_{s}})^{f_{i_{s}}})}}\underset{w_{2}}{\underbrace{(x_{i_{s}}, t_{i_{s}})^{f_{i_{s}}})(x_{j_{1}}, t_{j_{1}})^{f_{j_{1}}}}}\underset{w_{3}}{\underbrace{(x_{j_{2}}, t_{j_{2}})^{f_{j_{2}}}}}\dots\underset{w_{r}}{\underbrace{(x_{j_{r}}, t_{j_{r}})^{f_{j_{r}}}}}
\end{equation}

\begin{equation}\label{word3}
\mathfrak{w}=\underset{w_{1}}{\underbrace{\{(x_{i_{1}}, t_{i_{1}})(x_{i_{2}}, t_{i_{2}})\dots(x_{i_{s}}, t_{i_{s}})\}^{f_{j_{0}}}}}\underset{w_{2}}{\underbrace{(x_{j_{1}}, t_{j_{1}})^{f_{j_{1}}}}}\underset{w_{3}}{\underbrace{(x_{j_{2}}, t_{j_{2}})^{f_{j_{2}}}}}\dots\underset{w_{r}}{\underbrace{(x_{j_{r}}, t_{j_{r}})^{f_{j_{r}}}}}
\end{equation}

$(x_{i_{1}}, t_{i_{1}})\unlhd(x_{i_{2}}, t_{i_{2}})\unlhd\dots\unlhd(x_{i_{s}}, t_{i_{s}})$.\par\bigskip

We let $\sigma_{r, x}^{\Gamma}$ denote a function which applies $\sigma_{r}$ to each occurrence of a vertex $x\in \Gamma_{0}$, where $\Gamma_{0}$ is the set of vertices of a Brauer configuration whose polygons are given by words of type $\mathfrak{w}$ (see (\ref{word}), (\ref{word1})-(\ref{word3}) ).

\begin{equation}
\mathfrak{a}^{\Gamma}=\{\sigma^{\Gamma}_{r, x}\mid x\in\Gamma_{0}, r\in\mathbb{R}\}
\end{equation}

$\mathfrak{a}^{\Gamma}$ can be considered part of a $\Gamma$ \textit{labeling}, which contains symbols giving information on the ordering and number of vertices in the polygons.\par\bigskip

Identities (\ref{ex1}) gives the words determining polygons of a Brauer configuration $\Gamma=(\Gamma_{0},\Gamma_{1},\mu, \mathcal{O})$ with

\begin{itemize}

\item $\Gamma_{0}=\{a, e, g, \sigma_{-1}(b),  \sigma_{-1}(g), c_{0}, d_{0}\}$.
\item $\varphi_{\frac{n}{2}}(x)^{n}=x^{n+\frac{n}{2}}$.
\item $\sigma_{-1}(x^{n})=(\sigma_{-1}(x))^{n}$.
\item $c_{0}$ and $d_{0}$ are constants.
\item $\pi_{j}$ are labeling denoting appropriated permutations.
\item $w_{i}<w_{j}$ if $i<j$ for successor sequences.
\item $\mu(\alpha)=2$ or $\mu(\alpha)=1$ if the valency of a vertex $\alpha$ is whether 1 or greatest than 1.
\end{itemize}

\begin{equation}\label{ex1}
\begin{split}
w_{1}&=(\varphi_{8}(e^{16})\sigma_{-1}(\varphi_{8}(g^{16}))a^{16}d^{8}_{0}e^{8},\pi_{1})\\
w_{2}&=(d^{8}_{0}\sigma_{-1}(g^{8})d^{8}_{0}\sigma_{-1}(b^{8})\varphi_{8}(a^{16})d^{8}_{0},\pi_{2})\\
w_{3}&=(\varphi_{8}(e^{16})\sigma_{-1}(\varphi_{8}(g^{16}))a^{16}d^{8}_{0}g^{8},\pi_{3})\\
w_{4}&=(d^{8}_{0}\varphi_{16}(e^{32})c^{16}_{0},\pi_{4})\\
w_{5}&=(\varphi_{8}(e^{16})\sigma_{-1}(g^{16})a^{16}d^{8}_{0}e^{8},\pi_{5})\\
w_{6}&=(d^{8}_{0}\sigma_{-1}(g^{8})d^{8}_{0}\sigma_{-1}(b^{8})\varphi_{8}(a^{16})d^{8}_{0},\pi_{6})\\
w_{7}&=(\varphi_{8}(e^{16})\sigma_{-1}(\varphi_{8}(g^{16}))a^{16}d^{8}_{0}g^{8},\pi_{7})\\
w_{8}&=(\varphi_{8}(e^{16})\sigma_{-1}(\varphi_{8}(g^{16}))a^{16}d^{8}_{0}g^{8},\pi_{8})\\
w_{9}&=(d^{8}_{0}\varphi_{16}(e^{32})c^{16}_{0},\pi_{9})\\
\mathrm{dim}_{k}\hspace{0.1cm}\Lambda_{\Gamma}&=90426\\
\mathrm{dim}_{k}\hspace{0.1cm}Z(\Lambda_{\Gamma})&=612
\end{split}
\end{equation}

The following result regards Western musical writing without considering ornaments or interpretation symbols (e.g. fermata, coda, crescendo-decrescendo symbols or any other dynamic symbol, which can be included as label polygons).
\begin{Theorem}\label{Musical}
Any musical piece $\mathscr{M}$ written with the standard staff notation can be seen as a transposition ciphertext of the message ${M}(\mathfrak{M})$ of a Brauer configuration $\mathfrak{M}$ with the following properties: 
$\mathfrak{M}=(\mathfrak{M}_{0}, \mathfrak{M}_{1}, \mu, \mathcal{O})$.

\begin{itemize}
\item $\mathfrak{M}_{0}= \{a<\sigma_{-1}(b)<b<\sigma_{-1}(c)<c<d<\sigma_{-1}(e)<e<\sigma_{-1}(f)<f<\sigma_{-1}(g)<g<a\}\bigcup\mathcal{H}\bigcup \mathscr{S}_{0}$.
\item $\mathcal{H}=\{\varphi_{\frac{k}{2}}(x^{k})\mid x=a,b,\dots, g, k=2^{t}\hspace{0.1cm}\text{for some}\hspace{0.1cm}t>0\}$.
\item $\mathscr{S}_{0}= \{a_{0}, b_{0}, c_{0}, d_{0}, e_{0}, f_{0}, g_{0}\}$.
\item Any labeled polygon $U\in \mathfrak{M}_{1}$ (by a key $(\mathcal{K}_{i},\mathfrak{s}_{i})$) is determined by a word of type (\ref{word}) and (\ref{word1}, \ref{word}, \ref{word3}) with $f_{i,j}=2^{n}+k2^{n/2}$, $k, n\geq0$. We assume that if $\mathfrak{M}_{1}=\{U_{1}, U_{2},\dots, U_{k}\}$ then $\mathcal{K}_{i}$ is a permutation of the vertices $x_{i, j}$ in polygon $U_{i}$ and  $\mathfrak{s}_{i}\leq 64$ is the total number of its vertices counting repetitions. 
\item  \[\mu(\alpha)=
\begin{cases}
2,  & \hspace{0.2cm}\text{if}\hspace{0.2cm}val(\alpha)=1, \\
1,  & \hspace{0.2cm}\text{Otherwise.}
\end{cases}\]

\item The orientation $\mathcal{O}$ is given by enumerating the words $w_{1}<w_{2}<\dots<w_{n}$ associated with the polygons of $\Gamma_{1}$ ($|\Gamma_{1}|=n$) and assuming this order to build the successor sequences. 
\end{itemize}

\end{Theorem}

\textbf{Proof} Consider the following equivalences:

  \begin{enumerate}
\item 

$x^{64} $ $\longleftrightarrow$ ~~ (\Ganz, \hspace{0.1cm}$x$)   \\
$x^{32}$ $\longleftrightarrow$ ~~ (\Halb, \hspace{0.1cm}$x$)  \\
$x^{16}$ $\longleftrightarrow$ ~~  (\Vier,\hspace{0.1cm}$x$)  \\
$x^{8}$ $\longleftrightarrow$ ~~~  (\Acht,\hspace{0.1cm}$x$ )\\
$x^{4}$ $\longleftrightarrow$  ~~~ (\Sech,\hspace{0.1cm}$x$)    \\
$x^{2} $ $\longleftrightarrow$ ~~~  (\Zwdr,\hspace{0.1cm}$x$)   \\
Finally, the sixty fourth note and its rest have associated the vertices $x^{1}$ and $g_{0}$ respectively.
\item 
\begin{equation}
\begin{split}
\sigma_{-1}(x)&=\flat x.\\
\sigma_{1}(x)&=\sharp x.\\
\sigma_{-n}\sigma_{n}(x)&=\natural x.
\end{split}
\end{equation}

\item

$\varphi_{32}(x^{64}) $ $\longleftrightarrow$ ~~ ($\underset{\cdot}{\Ganz}$, \hspace{0.1cm}$x$), \quad $\varphi_{16}(x^{32}) $ $\longleftrightarrow$ ~~ ~ ($\underset{\cdot}{\Halb}$, \hspace{0.1cm}$x$),\quad $\varphi_{8}(x^{16}) $ $\longleftrightarrow$ ~~~  ($\underset{\cdot}{\Vier}$, \hspace{0.1cm}$x$)  $\dots$
\item $a^{64}_{0}$ $\longleftrightarrow$ ~~~  \GaPa\\
 $b^{32}_{0}$ $\longleftrightarrow$ ~~~  \HaPa\\
  $c^{16}_{0}$ $\longleftrightarrow$ ~~~  \ViPa\\
 $d^{8}_{0}$ $\longleftrightarrow$ ~~~  \AcPa\\
  $e^{4}_{0}$ $\longleftrightarrow$ ~~~  \SePa\\
   $f^{2}_{0}$ $\longleftrightarrow$ ~~~  \ZwPa\\
   
   \item $a$ $\longleftrightarrow$ C\\
    $b$ $\longleftrightarrow$ D\\
     $c$ $\longleftrightarrow$ E\\
      $d$ $\longleftrightarrow$ F\\
       $e$ $\longleftrightarrow$ G\\
        $f$ $\longleftrightarrow$ A\\
         $g$ $\longleftrightarrow$ B\\
         \item Standard scales are given by the following cycles.
         \begin{enumerate}
         \item $a<b<c<d<e<f<g<a$.
          \item $d<\sigma_{-1}(e)<e<\sigma_{-1}(f)<f<\sigma_{-1}(g)<g<a<\sigma_{-1}(b)<b<\sigma_{-1}(c)<c<d$ 
         \item $a<\sigma_{-1}(b)<b<\sigma_{-1}(c)<c<d<\sigma_{-1}(e)<e<\sigma_{-1}(f)<f<\sigma_{-1}(g)<g<a$
         \item $\frac{x+y}{2}=\sigma_{-1}(y)=\sigma_{1}(x)$, if $x<y$.
         \end{enumerate}
 \item Measures are given by the polygons in $\mathfrak{M}_{1}$.       
\item Grouping symbols (parentheses, brackets, and braces) in polygons define slurs, ties, tuples and chords.  

\end{enumerate}

Therefore, the musical piece $\mathscr{M}$ is nothing but a ciphertext of a transposition cryptosystem $S=(\mathcal{P}, \mathcal{C},\mathcal{K},\mathcal{E},\mathcal{D})$, where

\begin{equation}
\begin{split}
\mathcal{P}&=\{a,b,\dots, g,\sigma_{-1}(b),\sigma_{\pm1}(c), \sigma_{\pm1}(e), \sigma_{\pm1}(f),\sigma_{\pm1}(g)\}\bigcup \mathcal{H}\bigcup \mathscr{S}_{0},\\
\mathcal{C}&=\{(\mathfrak{n},x)\mid \mathfrak{n}\hspace{0.1cm}\text{is a musical note and}\hspace{0.1cm}x\in\mathcal{P}\}.\\
\mathcal{K}&=\{n\in\mathbb{N}\mid n=2^{k}, 1\leq k\leq 6\}.\\
e_{2^{k}}(x)&=x^{2^{k}}=(\mathfrak{n},x),\hspace{0.1cm}\text{for any}\hspace{0.1cm}x\in\mathcal{P}\hspace{0.1cm}\text{the duration of}\hspace{0.1cm}\mathfrak{n}\hspace{0.1cm}\text{is}\hspace{0.1cm}\frac{\mathfrak{d}}{2^{k}} (see (\ref{duration})).\\
\end{split}
\end{equation}
\begin{equation}
\begin{split}
d_{n}(\mathfrak{n},y)&=y^{2^{k}}\hspace{0.1cm}\text{if}\hspace{0.1cm}\text{the duration of}\hspace{0.1cm}\mathfrak{n}=\frac{\mathfrak{d}}{2^{k}}.\hspace{0.1cm} \text{Furthermore},\\
d_{n}(\flat\mathfrak{n},y)&=\sigma_{-1}(y^{2^{k}}).\\
d_{n}(\sharp\mathfrak{n},y)&=\sigma_{1}(y^{2^{k}}).\\
d_{n}(\natural\mathfrak{n},y)&=y^{2^{k}}.\\
d_{n}(\underset{\cdot}{\mathfrak{n}},y)&=\varphi_{2^{k-1}}(y^{2^{k}}).\\
d_{n}({\mathfrak{r}})&=x^{2^{k}}_{0},\hspace{0.1cm}\text{where}\hspace{0.1cm}\mathfrak{r}\hspace{0.1cm}\text{is a}\hspace{0.1cm}\text{rest}\hspace{0.1cm}\text{of duration}\hspace{0.1cm}\frac{\mathfrak{d}}{2^{k}}.\hspace{0.5cm}\square
\end{split}
\end{equation}

As an example the following is the musical notation of the words (\ref{ex1}) defined by the equivalences described in Theorem \ref{Musical}, the corresponding labeled Brauer configuration induced by these words is said to be \textit{$\mathfrak{M}$-reduced}:

\begin{equation}\label{SLYm}
\begin{split}
(\omega_{1},\mathcal{K})&= [b^{8}f^{8}e^{8}b^{8}][d^{8}\sigma_{-1}g^{8}e^{8}b^{8}]\longleftrightarrow [(\eighthnote, b)(\eighthnote, f)(\eighthnote, e)(\eighthnote, b)] [(\eighthnote, d)(\flat\eighthnote, g)(\eighthnote, e)(\eighthnote, b)]\\
(\omega_{2},\mathcal{K})&= \sigma_{-1}c^{16}[\sigma_{-1}g^{8}e^{8}]b^{16}f^{16}\longleftrightarrow (\flat\Vier, c)[(\flat\eighthnote, g)(\eighthnote, e)](\Vier, b)(\Vier, f)\\
(\omega_{3},\mathcal{K})&= b^{16}[f^{8}\sigma_{-1} c^{8}]a^{16}f^{16}\longleftrightarrow (\Vier, b)[(\eighthnote, f)(\flat\eighthnote, c)](\Vier, a)(\Vier, f)\\
(\omega_{4},\mathcal{K})&= [b^{8}f^{8}e^{8}b^{8}][d^{8}\sigma_{-1}g^{8}e^{8}b^{8}]\longleftrightarrow [(\eighthnote, b)(\eighthnote, f)(\eighthnote, e)(\eighthnote, b)]\\& [(\eighthnote, d)(\flat\eighthnote, g)(\eighthnote, e)(\eighthnote, b)]\\
(\omega_{5},\mathcal{K})&= [\sigma_{-1}c^{8}\sigma_{-1}g^{8}e^{8}b^{8}][\sigma_{-1}b^{8}\sigma_{-1}b^{8}\sigma_{-1}g^{8}e^{8}]\longleftrightarrow [(\flat\eighthnote, c)(\flat\eighthnote, g)(\eighthnote, e)(\eighthnote, b)]\\& [(\flat\eighthnote, b)(\flat\eighthnote, b)(\flat\eighthnote, g)(\eighthnote, e)]\\
(\omega_{6},\mathcal{K})&= b^{16}[f^{8}\sigma_{-1}c^{8}]a^{32}\longleftrightarrow (\Vier, b)[(\eighthnote, f)(\flat\eighthnote, c)](\Halb, a)\\
(\omega_{7},\mathcal{K})&= b^{16}[f^{8}\sigma_{-1}c^{8}]a^{16}f^{16}\longleftrightarrow (\Vier, b)[(\eighthnote, f)(\flat\eighthnote, c)](\Vier, a)(\Vier, f)\\
\end{split}
\end{equation}

The $\mathfrak{M}$-reduced Brauer configuration is defined as follows:

\begin{itemize}
\item $\mathfrak{M}_{0}=\{(\eighthnote, b), (\flat\eighthnote, b), (\flat\eighthnote, c),(\eighthnote, d),(\eighthnote, e), (\flat\eighthnote, f),(\flat\eighthnote, g),(\Vier, a), (\Halb, a),(\Vier, b)$, \\\hspace*{1.2cm}$(\flat\Vier, c), (\Vier, f), (\eighthnote, f)\}$.
\item $\Gamma_{1}=\{(\omega_{1}, \mathcal{K}),\dots, (\omega_{7},\mathcal{K})\}$.
\item The successor sequences associated with vertices are defined as follows (if no confusion arises, we use the same notation for copies of the same polygon):
\begin{equation}
\begin{split}
S_{(\eighthnote, b)}&=\omega_{1}<\omega_{1}<\omega_{1}<\omega_{4}<\omega_{4}<\omega_{4}<\omega_{5},\quad val((\eighthnote, b))=7,\\
S_{(\flat\eighthnote, b)}&=\omega_{5}<\omega_{5},\quad val((\flat\eighthnote, b))=2,\\
S_{(\flat\eighthnote, c)}&=\omega_{3}<\omega_{5}<\omega_{6}<\omega_{7},\quad val((\eighthnote, c))=4,\\
S_{(\eighthnote, d)}&=\omega_{1}<\omega_{4},\quad val((\eighthnote, d))=2,\\
S_{(\eighthnote, e)}&=\omega_{1}<\omega_{1}<\omega_{2}<\omega_{4}<\omega_{4}<\omega_{5}<\omega_{5},\quad val((\eighthnote, e))=7,\\
S_{(\eighthnote, f)}&=\omega_{1}<\omega_{3}<\omega_{4}<\omega_{6}<\omega_{7},\quad val((\eighthnote, f))=5,\\
S_{(\flat\eighthnote, g)}&=\omega_{1}<\omega_{2}<\omega_{4}<\omega_{5}<\omega_{5},\quad val((\flat\eighthnote, g))=5,\\
S_{(\Vier, a)}&=\omega_{3}<\omega_{7},\quad val((\Vier, a))=2,\\
S_{(\Halb, a)}&=\omega_{6},\quad val((\Halb, a))=1,\\
S_{(\Vier, b)}&=\omega_{2}<\omega_{3}<\omega_{6}<\omega_{7},\quad val((\Vier, b))=4,\\
S_{(\flat\Vier, c)}&=\omega_{2},\quad val((\flat\Vier, c))=1,\\
S_{(\Vier, f)}&=\omega_{2}<\omega_{3}<\omega_{7},\quad val((\Vier, f))=3.\\
\end{split}
\end{equation}

\item $\mu((\flat\Vier, c))=\mu((\Halb, a))=2$, $\mu(\alpha)=1$, for the remaining vertices $\alpha\in\Gamma_{0}$.

\item The signature of the key $\mathcal{K}$ is given by the pair $(\trebleclef,\begin{array}{c}n \\2^m\end{array})$, meaning that each polygon $\omega_{i}$, $1\leq i\leq 7$ has $\mathfrak{s}_{i}=n\times 2^{6-m}$ vertices (or its equivalent number).

\item $\mathrm{dim}_{k}\hspace{0.1cm}\Lambda_{\mathfrak{M}}=14+162=176$.
\item  $\mathrm{dim}_{k}\hspace{0.1cm}Z(\Lambda_{\mathfrak{M}})=1+10+4-12+7+12-2=20$ (note that, there are 2 vertices whose valency equals 1. Furthermore, the Brauer quiver $Q_{\mathfrak{M}}$ has 12 loops and 7 vertices).

\end{itemize}
\par\bigskip

\begin{figure}[H]
		\centering
	\includegraphics[scale=0.7]{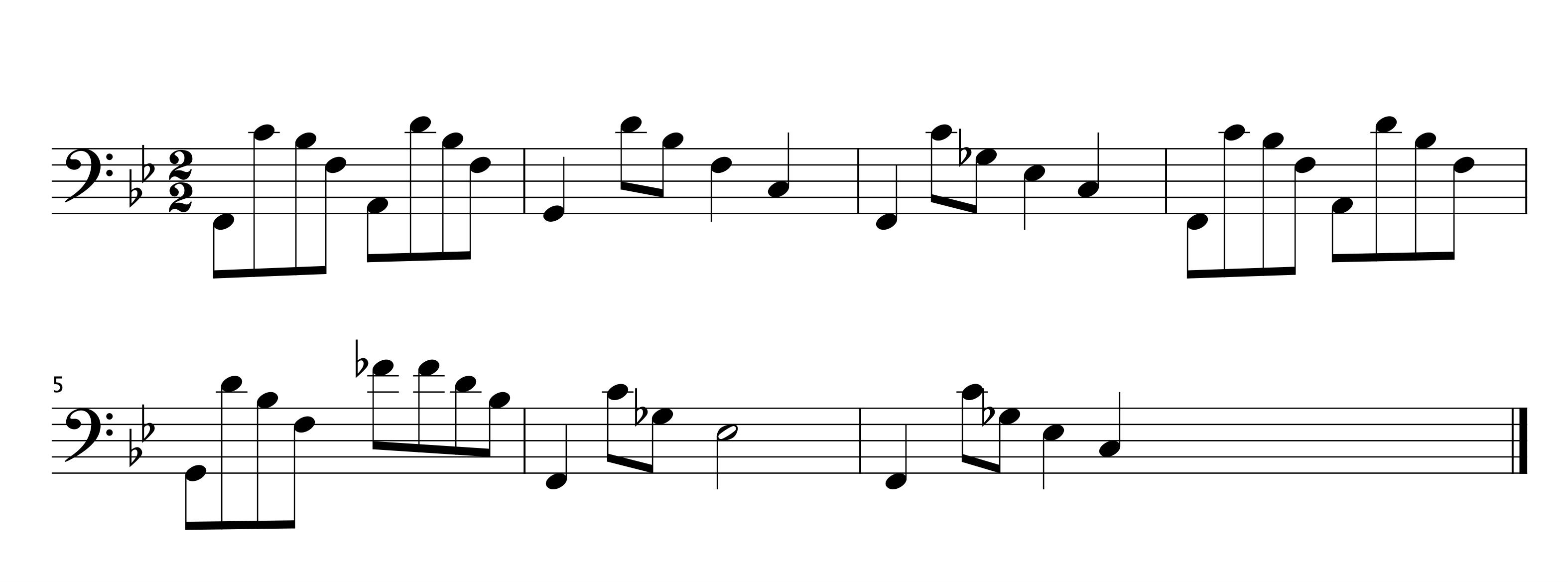}
	\caption{Ciphertext defined by the Brauer configuration (\ref{SLYm}). }
	\label{SLYmod1}
\end{figure}

Figure \ref{SLYmod1} shows the labeled Brauer message $M(\mathfrak{M},\mathcal{K})$ written using the staff notation. Here, the label of the Brauer configuration is given by a signature of the form  $\mathcal{K}=(\bassclef, \begin{array}{c}2  \\2 \end{array}, \mathfrak{a}^{\mathfrak{M}}=\{\sigma_{-1}(c),\sigma_{-1}(g)\})$, which means that the circular order associated with the vertices is given by the following sequence (\ref{seq}) and that each polygon has two half notes or their equivalents.

  \begin{equation}\label{seq}
 d<\sigma_{-1}(e)<e<\sigma_{-1}(f)<f<\sigma_{-1}(g)<g<a<\sigma_{-1}(b)<b<\sigma_{-1}(c)<c<d 
  \end{equation}

 \vspace{-6pt} 

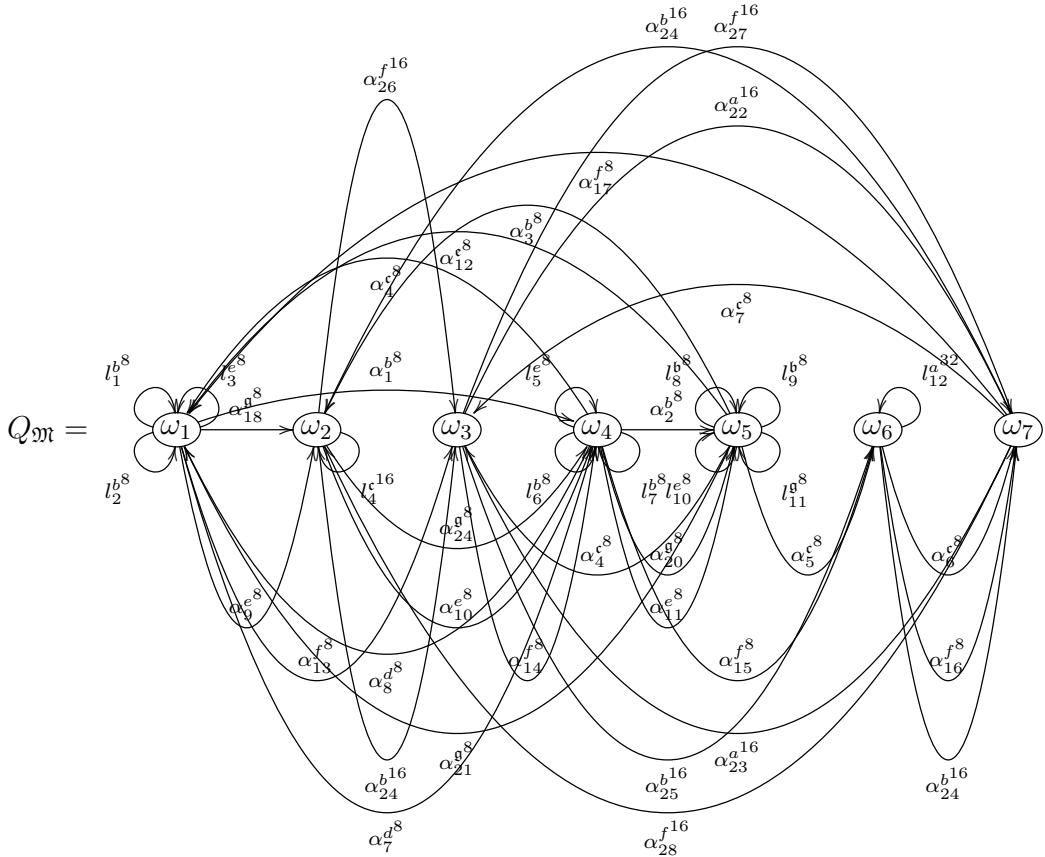
\begin{figure}[H]
 \begin{equation}
Q_{\mathfrak{M}}=\xymatrix@=35pt{*+[o][F-]{\omega_{1}}\ar@/^15pt/@[\blue][rrr]^{\alpha^{b^{8}}_{1}}\ar@(l,d)@[\blue]_{l_2^{b^{8}}}\ar@(l,u)@[\blue]^{l_1^{b^{8}}}\ar@(r,u)@[\blue]_{l_3^{e^{8}}}\ar@/_145pt/@[\blue][rrr]_{\alpha^{{d}^{8}}_{7}}\ar@/_75pt/@[\blue][r]^{\alpha^{{e}^{8}}_{9}}\ar@/_95pt/@[\blue][rr]^{\alpha^{{f}^{8}}_{13}}\ar@/_0pt/@[\blue][r]^{\alpha^{\mathfrak{g}^{8}}_{18}}&*+[o][F-]{\omega_{2}}\ar@/_75pt/@[\blue][rr]^{\alpha^{{e}^{8}}_{10}}\ar@/_45pt/@[\blue][rr]^{\alpha^{\mathfrak{g}^{8}}_{24}}\ar@/_125pt/@[\blue][r]_{\alpha^{{b}^{16}}_{24}}\ar@(r,d)@[\blue]^{l_4^{\mathfrak{c}^{16}}}\ar@/^125pt/@[\blue][r]^{\alpha^{{f}^{16}}_{26}}&*+[o][F-]{\omega_{3}}
\ar@/_55pt/@[\blue][rr]^{\alpha^{\mathfrak{c}^{8}}_{4}}\ar@/_95pt/@[\blue][r]^{\alpha^{{f}^{8}}_{14}}\ar@/^115pt/@[\blue][rrrr]^{\alpha^{a^{16}}_{22}}\ar@/_125pt/@[\blue][rrr]_{\alpha^{{b}^{16}}_{25}}\ar@/^145pt/@[\blue][rrrr]^{\alpha^{{f}^{16}}_{27}}
&*+[o][F-]{\omega_{4}}\ar@(l,d)@[\blue]_{l_6^{b^{8}}}\ar@(r,d)@[\blue]^{l_7^{b^{8}}}\ar@/_0pt/@[\blue][r]^{\alpha^{b^{8}}_{2}}\ar@/^85pt/@[\blue][lll]^{\alpha^{d^{8}}_{8}}\ar@/_65pt/@[\blue][lll]^{\alpha^{\mathfrak{c}^{8}}_{4}}\ar@(l,u)@[\blue]^{l_5^{e^{8}}}\ar@/_75pt/@[\blue][r]^{\alpha^{{e}^{8}}_{11}}\ar@/_95pt/@[\blue][rr]^{\alpha^{{f}^{8}}_{15}}\ar@/_55pt/@[\blue][r]^{\alpha^{\mathfrak{g}^{8}}_{20}}&*+[o][F-]{\omega_{5}}\ar@(l,u)@[\blue]^{l_8^{\mathfrak{b}^{8}}}\ar@(r,u)@[\blue]_{l_{9}^{\mathfrak{b}^{8}}}\ar@(d, l)@[\blue]^{l_{10}^{{e}^{8}}}\ar@/_85pt/@[\blue][lll]^{\alpha^{b^{8}}_{3}}\ar@/_55pt/@[\blue][r]^{\alpha^{\mathfrak{c}^{8}}_{5}}\ar@(r,d)@[\blue]^{l_{11}^{\mathfrak{g}^{8}}}\ar@/_75pt/@[\blue][llll]^{\alpha^{\mathfrak{e}^{8}}_{12}}\ar@/^115pt/@[\blue][llll]^{\alpha^{\mathfrak{g}^{8}}_{21}}&*+[o][F-]{\omega_{6}}\ar@/_55pt/@[\blue][r]^{\alpha^{\mathfrak{c}^{8}}_{6}}\ar@/_95pt/@[\blue][r]^{\alpha^{{f}^{8}}_{16}}\ar@(r,u)@[\blue]_{l_{12}^{a^{32}}}\ar@/_125pt/@[\blue][r]_{\alpha^{{b}^{16}}_{24}}&*+[o][F-]{\omega_{7}}\ar@/_55pt/@[\blue][llll]^{\alpha^{\mathfrak{c}^{8}}_{7}}\ar@/_105pt/@[\blue][llllll]^{\alpha^{{f}^{8}}_{17}}\ar@/^115pt/@[\blue][llll]^{\alpha^{a^{16}}_{23}}\ar@/_145pt/@[\blue][lllll]_{\alpha^{{b}^{16}}_{24}}\ar@/^145pt/@[\blue][lllll]^{\alpha^{{f}^{16}}_{28}}
}\notag
\end{equation}

\caption{Brauer quiver associated with the $\mathfrak{M}$-reduced Brauer configuration (\ref{SLYm}).}\label{mutfinal10}
\end{figure}

Figure \ref{mutfinal10} shows the Brauer quiver $Q_{\mathfrak{M}}$ induced by the $\mathfrak{M}$-reduced Brauer configuration (\ref{SLYm}). Here $\mathfrak{z}^{j}$ corresponds with a vertex $\sigma_{-1}z^{j}$. \par\bigskip

The Brauer configuration algebra $\Lambda_{\mathfrak{M}}$ is bounded by an ideal $I$ generated by relations of the following types:

\begin{itemize}
\item $C^{\mu(j)}_{j}-C^{\mu(i)}_{i}$, if $i, j$ belong to the same polygon. Whereas, $C_{i}$ and $C_{j}$ are corresponding special cycles.
\item $\alpha^{s'}_{s}\alpha^{t'}_{t}$ if $t\neq t'$.
\item $(l^{s'}_{s})^{2}$, if $l^{s'}_{s}$ is a loop.
\item $C^{\mu(j)}_{j}f$, if $C_{j}$ is a special cycle and $f$ is its first arrow.
\end{itemize}

\subsection{Brauer Analysis of Bach's Canons}\label{Bach}

This section gives some properties of the Brauer configuration algebras induced by Bach's canons known as  canon \^a  6 Voc, canon \^a 4. Voc: Perpetuus (BWV 1073), canon 1, \^a 2 (the crab canon), canon \^a4 Quaerendo Invenietis. 
The musical notes are represented in the Euclidean plane to build some of the most common symbols used by Bach in his manuscripts. In this line, Bach's canons are interpreted as ciphertexts of a transposition cryptosystem whose plaintext consists of such symbols, which give the structure or form of the studied Bach's canons. Algorithm \ref{Alg} describes the procedure: 

\begin{algorithm}
\begin{itemize}
 \item Choose a vertex circular order of the set $\{a,b,\dots,g\}$ according to the signature clef. The order associated with a treble, (bass, C) clef starts with the letter $e$, ($d, a$, respectively).\\
 \item Segregate the notes into classes according to the equivalence relation $\sim$, such that $(\mathfrak{n},x)\sim(\mathfrak{n}',x')$ if and only if $\mathfrak{n}=\mathfrak{n}'$ and $x=x'$.
\item  Associate with the first occurrence of the $i$th note $\mathfrak{n}_{i}$ a point $(\varepsilon^{i,1},\varepsilon^{i,2})$ in the Euclidean plane. Where, the first coordinate $\varepsilon^{i,1}$ is a pair $(\mathfrak{n}^{1}_{i,1}, \mathfrak{n}^{1}_{i,2})$, $\mathfrak{n}^{1}_{i,1}\in\{\Ganz,\Halb,\Vier,\dots\}$, $\mathfrak{n}^{2}_{i,2}\in\{a,b, c,\dots, g\}$ is given by the appearance order of the note in the musical piece. \\
\item $\varepsilon^{i, 2}$ is an integer number assigned to the vertical coordinate of the note. The clef note has assigned the number $\varepsilon^{1,2}=0$ and the remaining notes have assigned the number $i$ if the vertices order is $a_{0}<a_{1}<\dots $, $a_{i}=\mathfrak{n}^{2}_{i,2}\in\{a, b, c,\dots, g\}$. It is positive above (below) the clef note $(\varepsilon^{1,1},\varepsilon^{1,2})$, if it is written in standard (reversed) orientation.
\item Each letter or number is built by connecting consecutive points. If there are repeated notes only the first of them is represented in the diagram. In such a way that points $(\varepsilon^{i,1},\varepsilon^{i,2})$ define a discrete function. Generally, two consecutive points with the same value $\varepsilon^{i,2}$ are not connected, such a connection is drawn
only for aesthetic reasons as well as additional lines can be drawn to represent a symbol containing a cycle (e.g., A, B, $\alpha$). If necessary, lines of the previous symbol can be used to construct a new one. In particular, vertical lines can be used to represent rests.
\end{itemize}
 \caption{This algorithm gives the steps to built a graph induced by the musical notes of Bach's canons.}
 \label{Alg}
\end{algorithm}
\par\bigskip

\begin{figure}[H]
		\centering
	\includegraphics[scale=0.5]{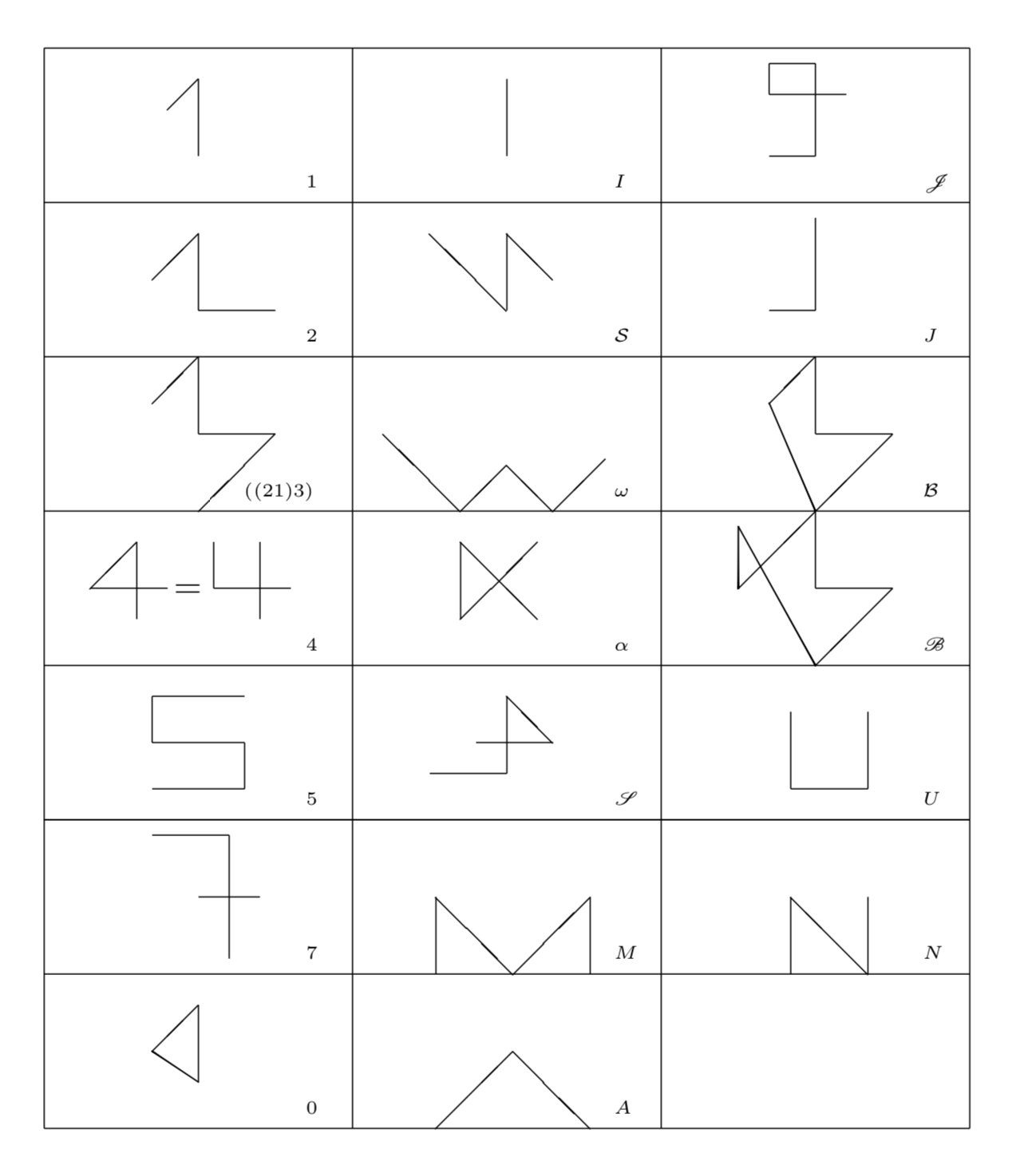}\\
	\includegraphics[scale=0.5]{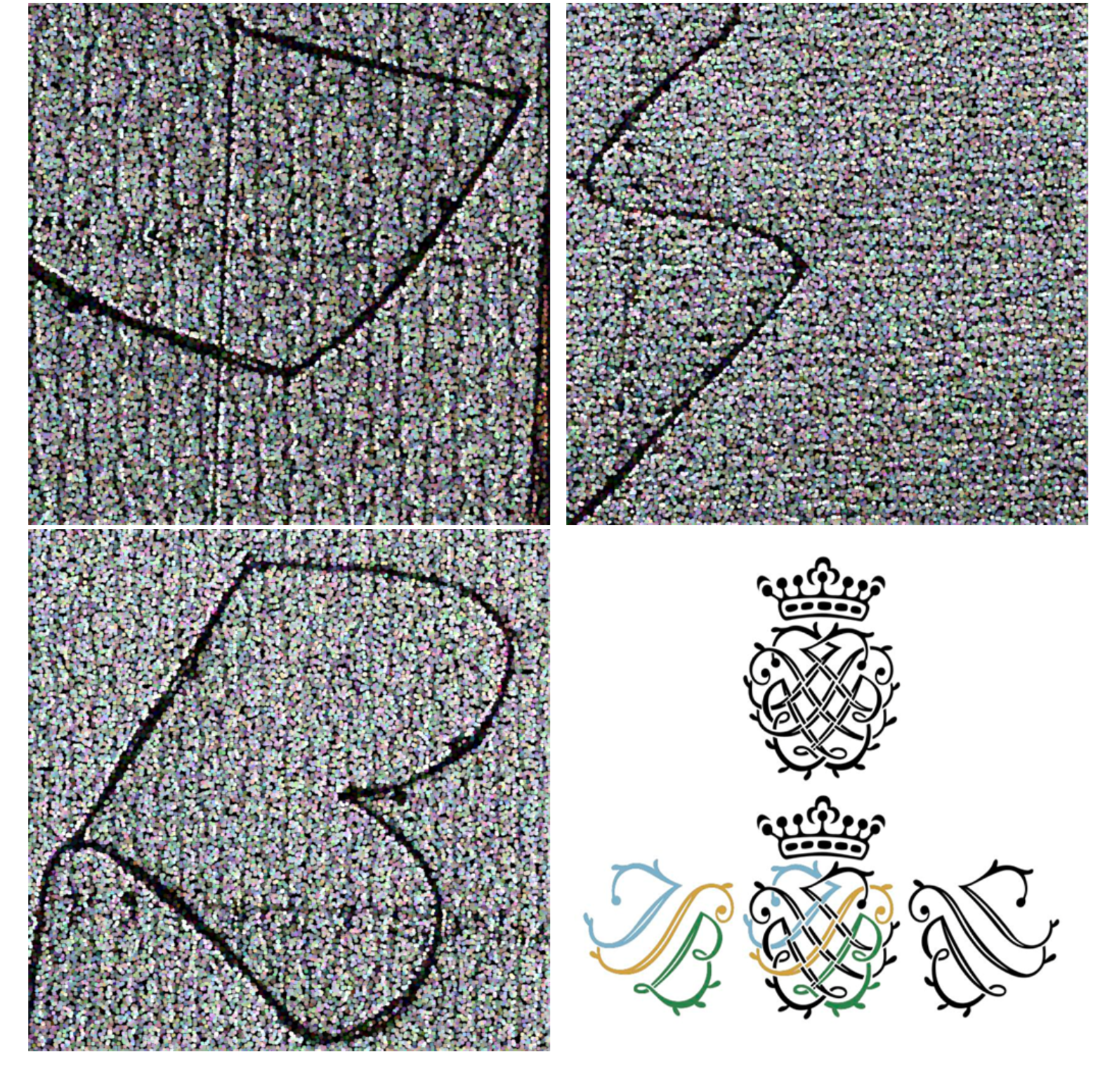}
	\caption{Symbols that are constructed using points and lines. And examples of symbols written or mentioned in Bach's manuscripts. Bach's monogram on his wax seal based on his initials and different versions of the Greek letters $\alpha$ and $\omega$ \cite{Seal, Chris}. }
	\label{dictionary}
\end{figure}

\begin{figure}[H]
  \centering
		\includegraphics[scale=0.5]{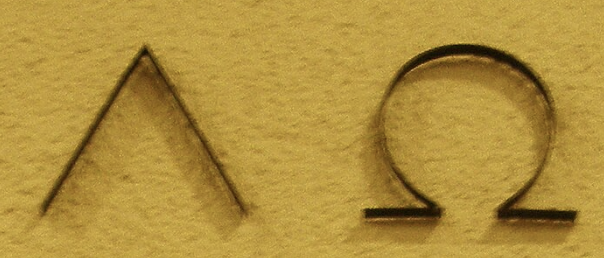}\includegraphics[scale=0.15]{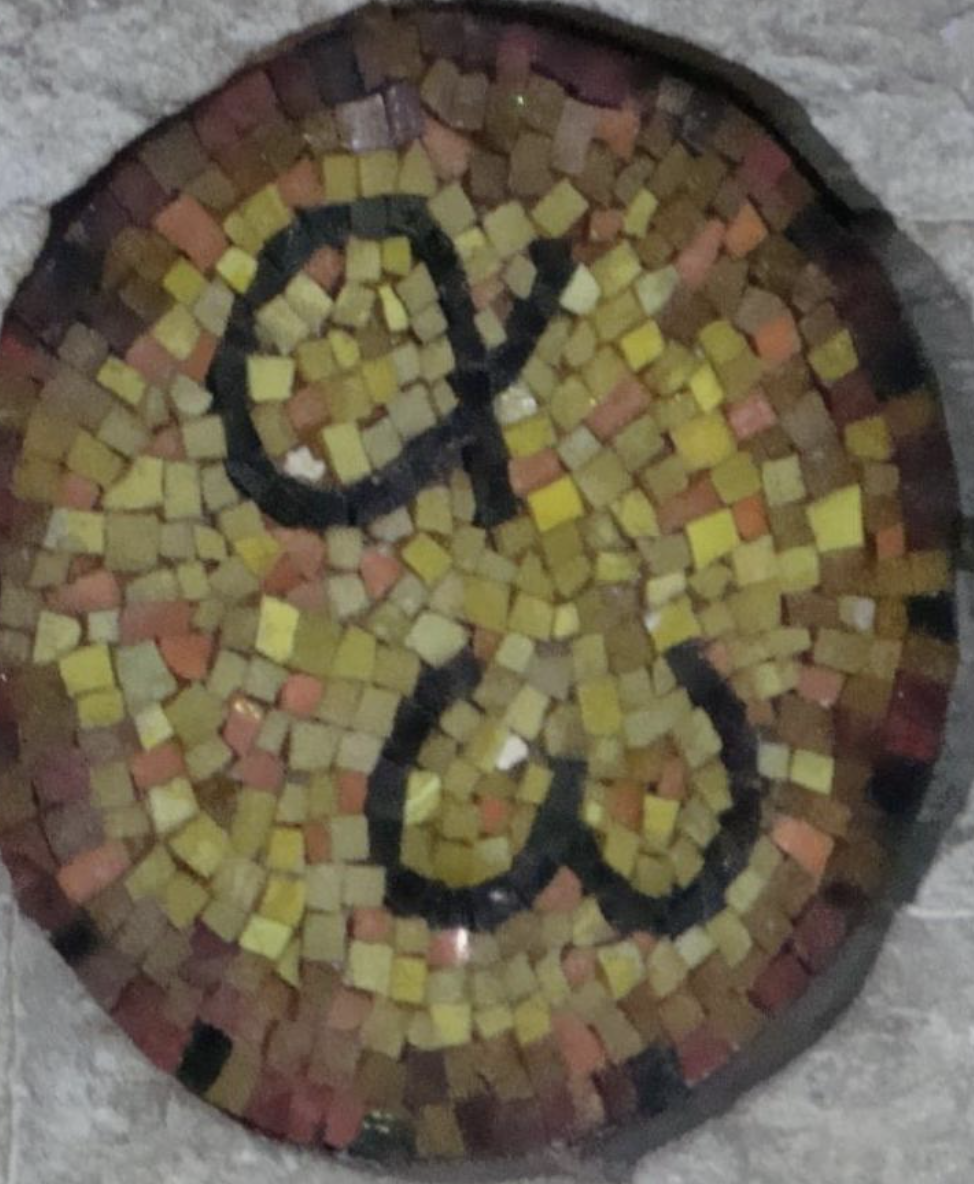}
	
	\label{chrismon(1)}
\end{figure}

\begin{figure}[H]
		\centering
	\includegraphics[scale=0.5]{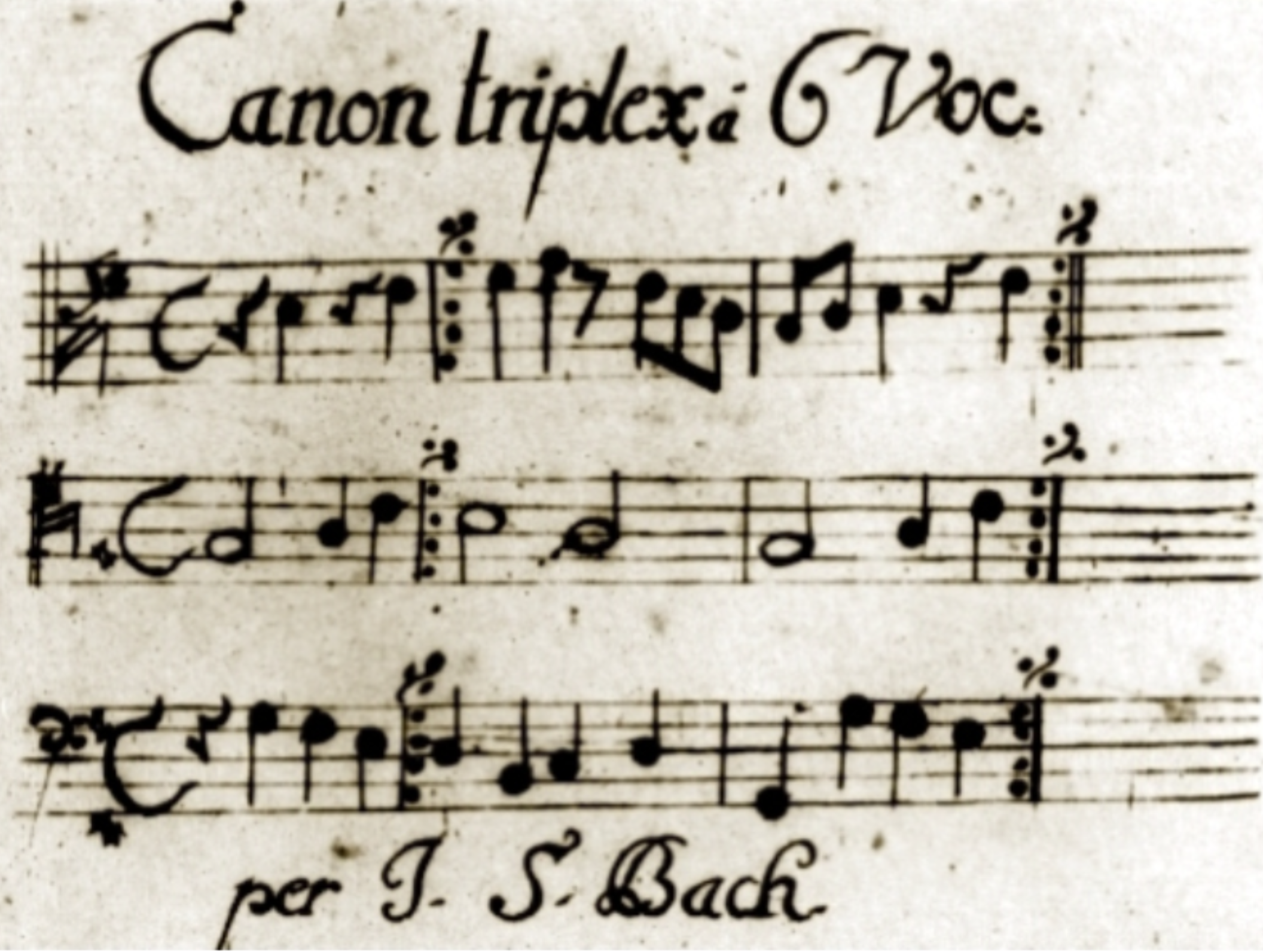}
	\caption{The reduced Brauer message of the Brauer configuration (\ref{A6Br}) gives the Bach's canon \^a6 Voc \cite{a6V}. }
	\label{a6V1}
\end{figure}

Figure \ref{dictionary} shows examples of symbols and letters that can be built using Algorithm \ref{Alg}. Bach was used to inverting some of these symbols in his manuscripts (see Figures  \ref{chrismon(1)} and \ref{a6V1}, where Bach's monogram in his famous wax seal is built upon his initials. And a bass symbol is used to write the number 6).

\vspace{-6pt} 

Bach's canon Triplex \^a6 Voc is obtained by reducing the Brauer message (\ref{MA6}). It allows defining the labeled Brauer configuration $\mathfrak{M}^{A6}=(\mathfrak{M}^{A6}_{0},\mathfrak{M}^{A6}_{1},\mu,\mathcal{O})$  whose words or polygons  (\ref{WA6}) are obtained after an appropriate segmentation.

\begin{equation}\label{MA6}
\begin{split}
c^{16}_{0}&e^{16}c^{16}_{0}d^{16}c^{16}b^{16}d^{8}_{0}[d^{8}e^{8}f^{8}][g^{8}f^{8}]e^{16}c^{16}_{0}d^{16}g^{32}f^{16}d^{16}e^{32}f^{32}g^{32}f^{16}d^{16}c^{16}_{0}c^{16}d^{16}e^{16}f^{16}a^{16}g^{16}\\
f^{16}&c^{16}c^{16}d^{16}e^{16}.
\end{split}
\end{equation}
\begin{equation}\label{WA6}
\begin{split}
(\omega_{1},\mathcal{K})&=c^{16}_{0}e^{16}c^{16}_{0}d^{16}\longleftrightarrow \ViPa(\Vier, e)\ViPa(\Vier, d)\\
(\omega_{2},\mathcal{K})&=c^{16}b^{16}d^{8}_{0}[d^{8}e^{8}f^{8}]\longleftrightarrow(\Vier, c)(\Vier, b)\AcPa[(\eighthnote, d)(\eighthnote, e)(\eighthnote, f)]\\
\omega_{3}, \mathcal{K})&=[g^{8}f^{8}]e^{16}c^{16}_{0}d^{16}\longleftrightarrow[(\eighthnote, g)(\eighthnote, f)](\Vier, e)\ViPa(\Vier, d)\\
(\omega_{4},\mathcal{K})&=g^{32}f^{16}d^{16}\longleftrightarrow(\Halb, g)(\Vier, f)(\Vier, d)\\
(\omega_{5,\mathcal{K})}&=e^{32}f^{32}\longleftrightarrow(\Halb, e)(\Halb, f)\\
(\omega_{6},\mathcal{K})&=g^{32}f^{16}d^{16}\longleftrightarrow(\Halb, g)(\Vier, f)(\Vier, d)\\
(\omega_{7},\mathcal{K})&=c^{16}_{0}c^{16}d^{16}e^{16}\longleftrightarrow \ViPa(\Vier, c)(\Vier, d)(\Vier, e)\\
(\omega_{8},\mathcal{K})&=f^{16}a^{16}g^{16}f^{16}\longleftrightarrow (\Vier, f)(\Vier, a)(\Vier, g)(\Vier, f)\\
(\omega_{9},\mathcal{K})&=c^{16}c^{16}d^{16}e^{16}\longleftrightarrow (\Vier, c)(\Vier, c)(\Vier, d)(\Vier, e)\\
\end{split}
\end{equation}

The successor sequences are given by the following identities:
 \begin{equation}
  \begin{split}
  S_{\ViPa}&=\omega_{1}<\omega_{1}<\omega_{3}<\omega_{7},\\
  S_{(\Vier, e)}&=\omega_{1}<\omega_{3}<\omega_{7}<\omega_{9},\\
   S_{(\Vier, d)}&=\omega_{1}<\omega_{3}<\omega_{4}<\omega_{6}<\omega_{7}<\omega_{9},\\
    S_{(\Vier, c)}&=\omega_{2}<\omega_{7}<\omega_{9}<\omega_{9},\\
     S_{(\Vier, b)}&=\omega_{2},\\
     S_{\AcPa}&=\omega_{2},\\
 S_{(\eighthnote, d)}&=\omega_{2},\\
S_{(\eighthnote, e)}&=\omega_{2},\\
 S_{(\eighthnote, f)}&=\omega_{2}<\omega_{3},\\
 S_{(\eighthnote, g)}&=\omega_{3},\\
 S_{(\Halb, g)}&=\omega_{4}<\omega_{6},\\
  S_{(\Vier, f)}&=\omega_{4}<\omega_{6}<\omega_{8}<\omega_{8},\\
   S_{(\Halb, e)}&=\omega_{5},\\
    S_{(\Halb, f)}&=\omega_{5},\\
    S_{(\Vier, a)}&=\omega_{8},\\
    S_{(\Vier, g)}&=\omega_{8}.  
  \end{split}
  \end{equation}

 The following identities give the main properties of the Brauer configuration algebra $\Lambda_{\mathfrak{M}^{A6}}$ induced by the Bach's canon \^a6 Voc. 
 \begin{equation}\label{A6Br}
 \begin{split}
 \mathfrak{M}^{A6}_{0}&=\{\ViPa,\AcPa, (\Vier, a),(\Vier, b),(\Vier, c), (\Vier, d),(\Vier, e),(\Vier, f),(\Vier, g),(\eighthnote, d), (\eighthnote, e),(\eighthnote, f),\\&
 (\eighthnote, g), (\Halb, e), (\Halb, f), (\Halb,g)\},\\
  \mathfrak{M}^{A6}_{1}&=\{\omega_{1},\dots, \omega_{9}\},\quad \text{(see identities (\ref{WA6}))},\\ 
  \mathcal{K}&=(\bassclef, \begin{array}{c}4  \\4 \end{array},d<c<b<\dots<f<e<d, \text{4th line})\\
  \{\alpha\in  \mathfrak{M}^{A6}_{0}\mid val(\alpha)=1\}&=\{(\Vier, b), \AcPa,(\eighthnote, d), (\eighthnote, e), (\eighthnote, g), (\Halb, e), (\Halb, f), (\Vier, a), (\Vier, g)\},\\
    \{\alpha\in  \mathfrak{M}^{A6}_{0}\mid val(\alpha)=2\}&=\{(\eighthnote, f), (\Halb, g)\},\\
       \{\alpha\in  \mathfrak{M}^{A6}_{0}\mid val(\alpha)=4\}&=\{\ViPa, (\Vier, e), (\Vier, c), (\Vier, f)\},\\
          \{\alpha\in  \mathfrak{M}^{A6}_{0}\mid val(\alpha)=6\}&=\{(\Vier, d)\},\\
          |\mathfrak{M}^{A6}_{0}|&=16,\quad\text{there are 14 vertices, excluding the rests.}\\
          \#Loops\hspace{0.1cm}Q_{\mathfrak{M}^{A6}}&=12,\\
          \mathrm{dim}_{k}\hspace{0.1cm}\Lambda_{\mathfrak{M}^{A6}}&=109,\\
           \mathrm{dim}_{k}\hspace{0.1cm}Z(\Lambda_{\mathfrak{M}^{A6}})&=22.
 \end{split}
 \end{equation}

\begin{figure}[H]
		\centering
		\includegraphics[scale=0.3]{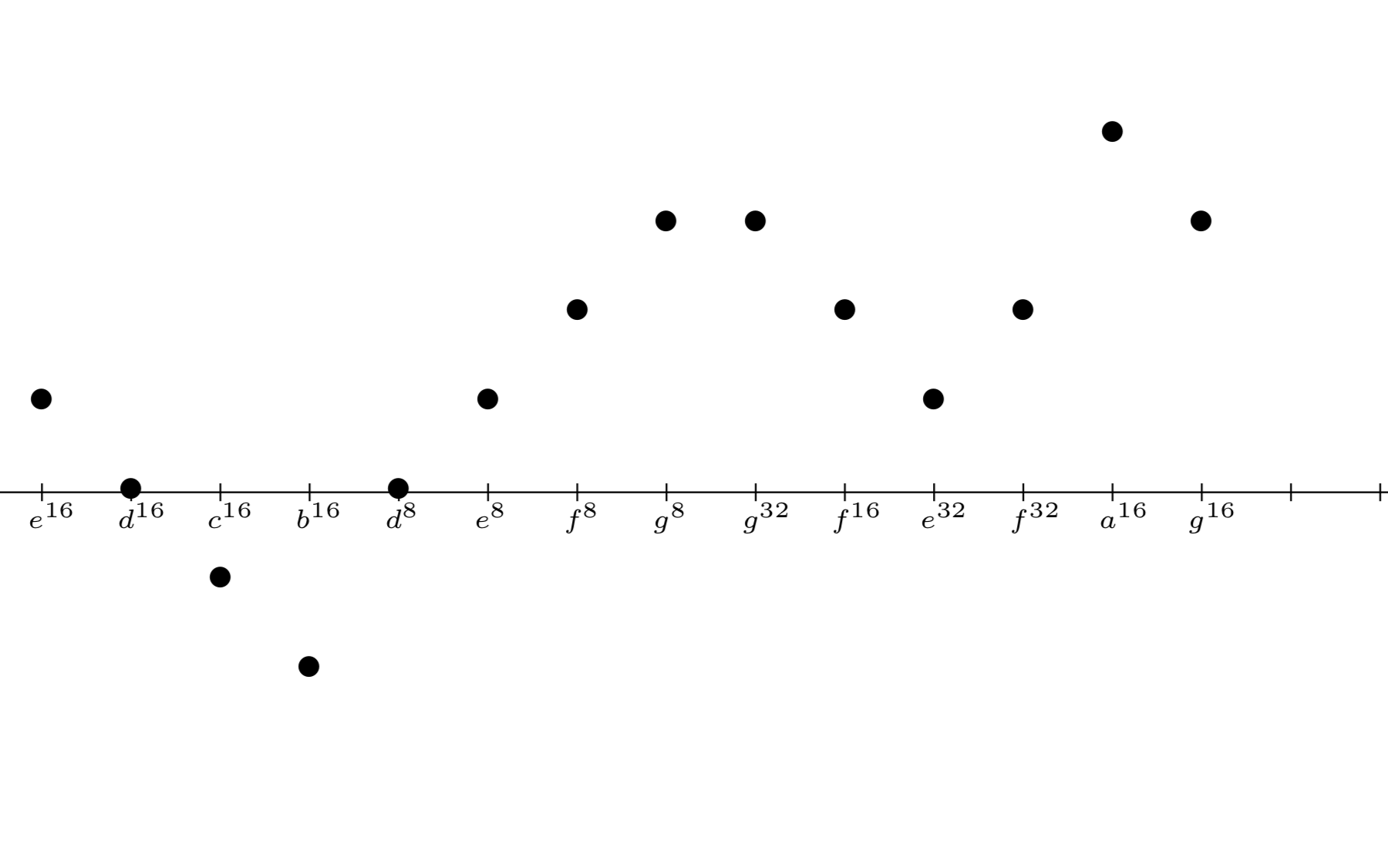}\\
		\includegraphics[scale=0.3]{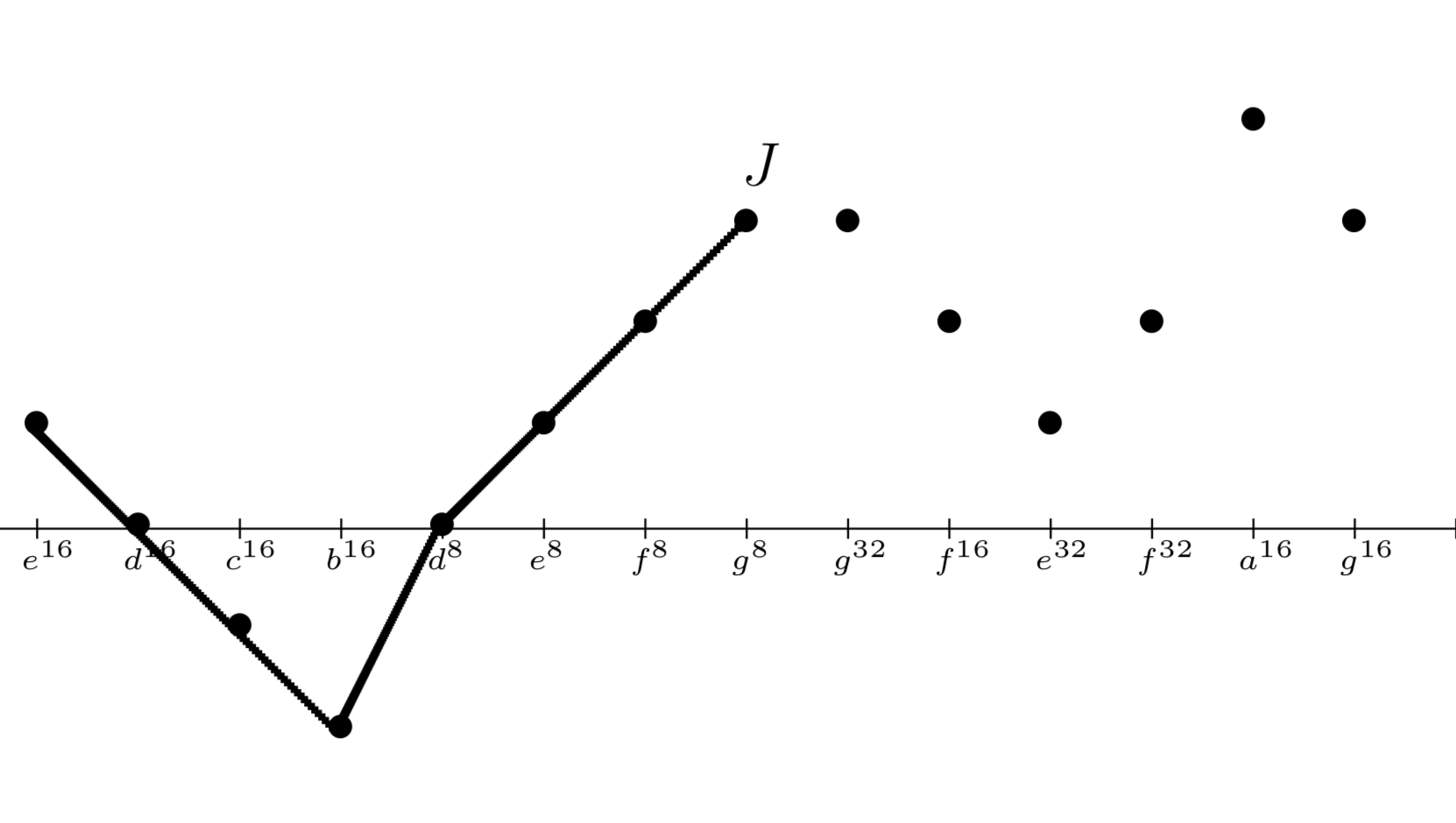}\\
		\includegraphics[scale=0.3]{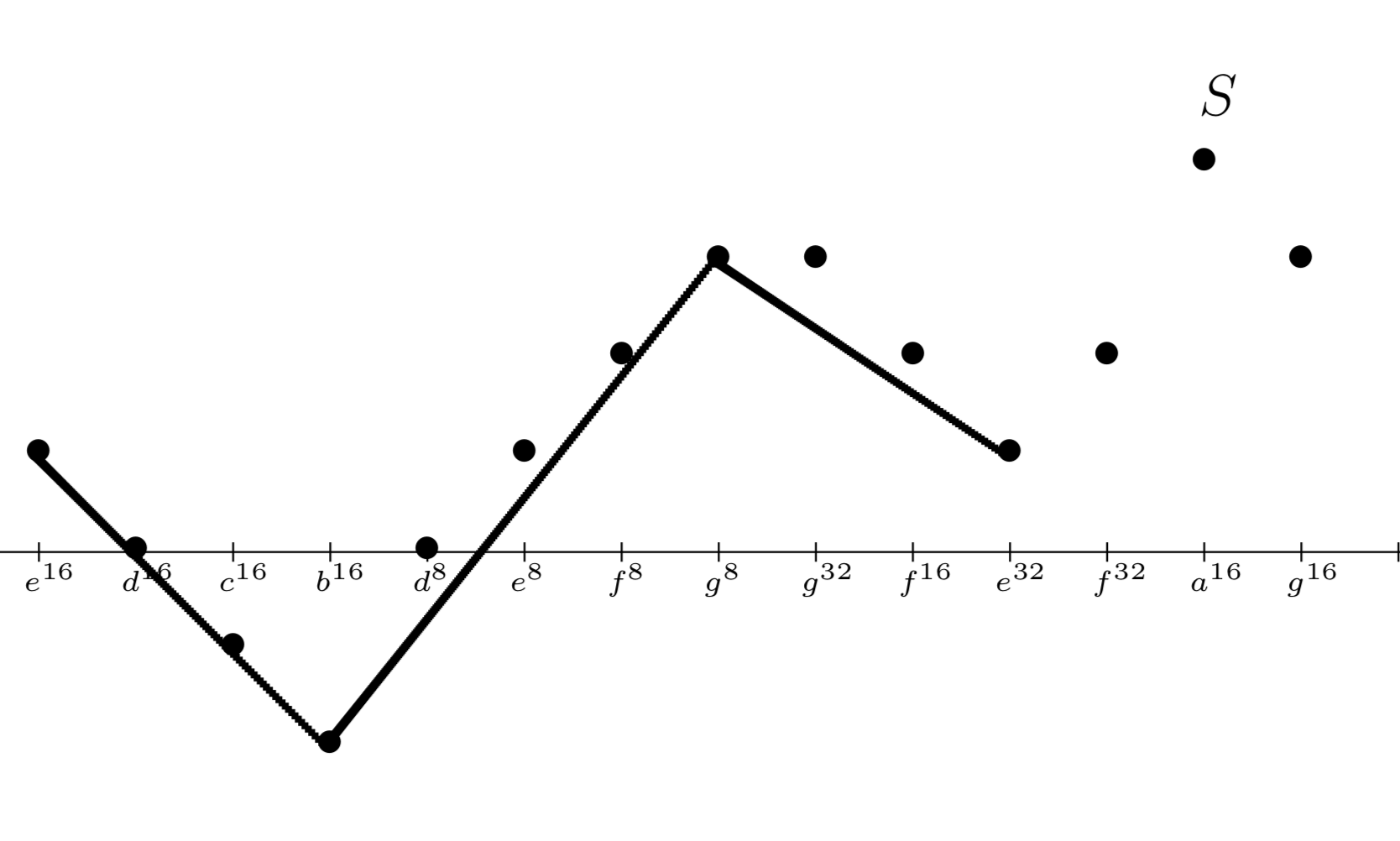}\\
		\includegraphics[scale=0.3]{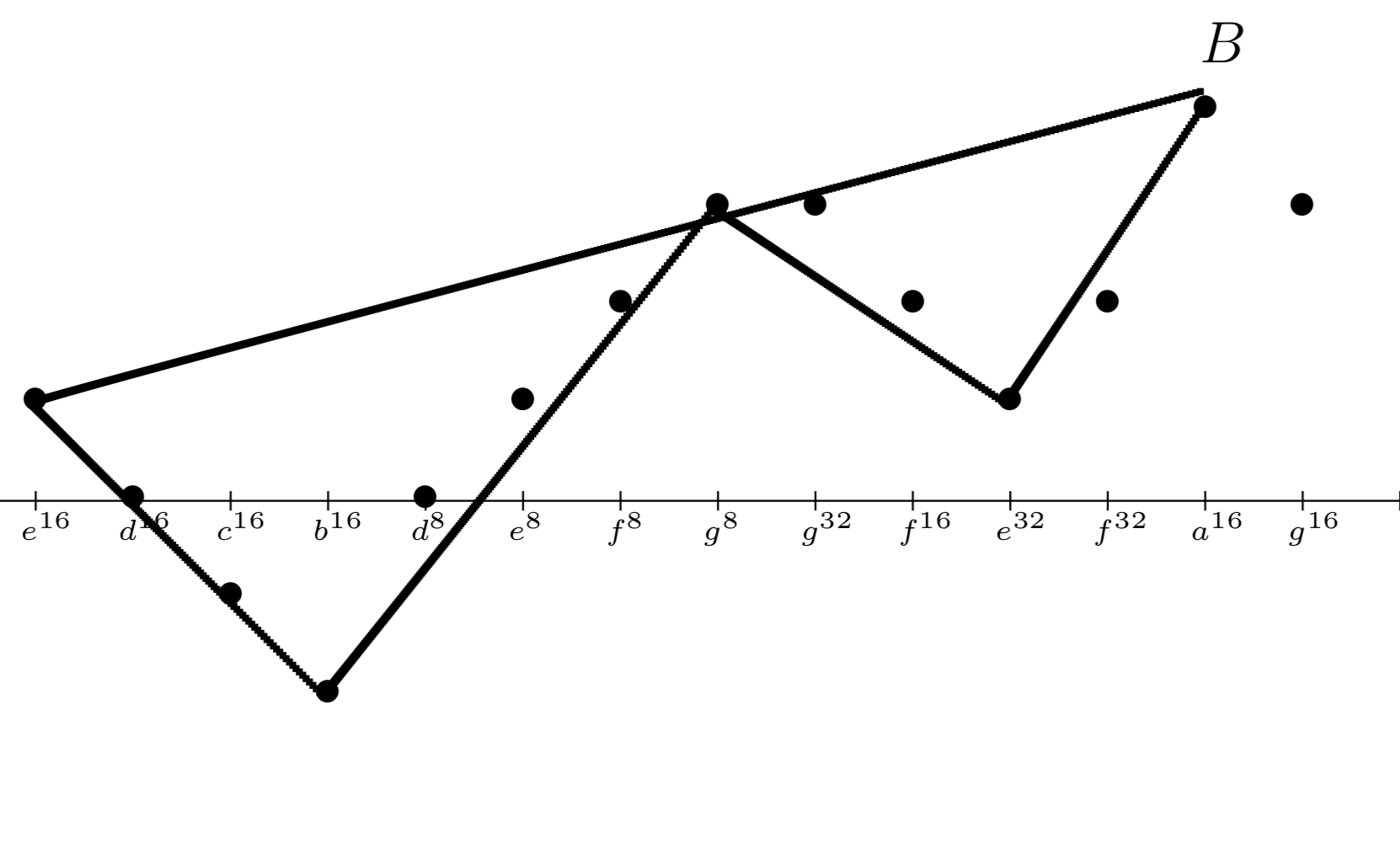}\\		
	\includegraphics[scale=0.6]{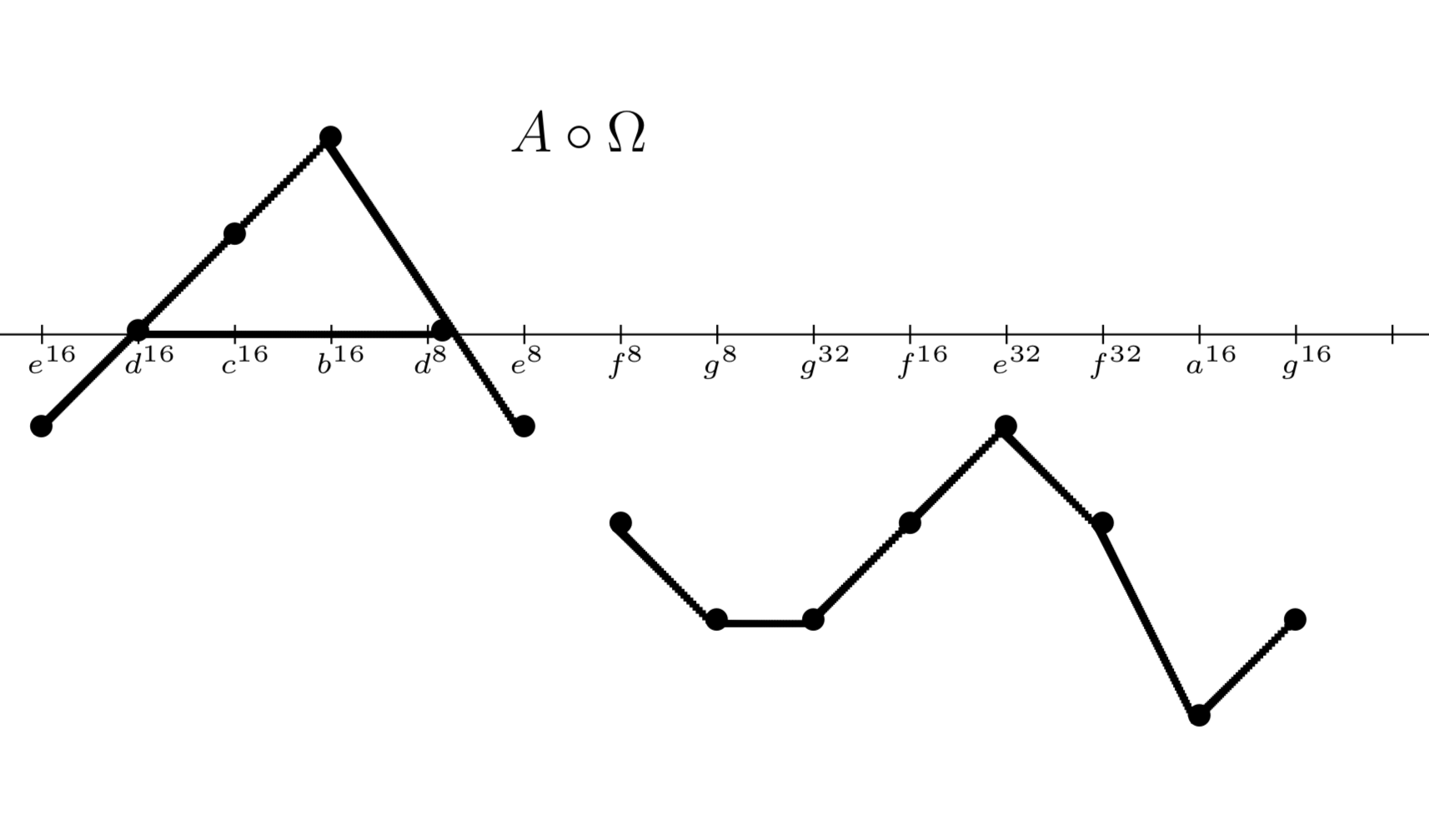}
	
	\caption{Fourteen vertices are enough to build up Bach's initials. Symbols alpha and omega arise by applying a reflection to these points (see Algorithm \ref{Alg}). }
	\label{a6VS}
\end{figure}

The following list gives the 14 points (without rests) induced by canon \^a6 Voc. Such points are represented in the diagram shown in Figure  \ref{a6VS}. The initials of Johann Sebastian Bach and the Greek letters alpha and omega can be built up via an appropriate connection of the points.

\begin{equation}
\begin{split}
 ( {(\Vier, e);\varepsilon^{1}})&, ({(\Vier, d);\varepsilon^{0}}),
({(\Vier, c);\varepsilon^{-1}})),
((\Vier, b);\varepsilon^{-2}), ({(\eighthnote, d);\varepsilon^{0}}),
({(\eighthnote, e);\varepsilon^{1}}),
({(\eighthnote, f);\varepsilon^{2}}),\\
({(\eighthnote, g);\varepsilon^{3}})&, ({(\Halb, g);\varepsilon^{3}}),
  ({(\Vier, f);\varepsilon^{2}}), ({(\Halb, e);\varepsilon^{1}}), ({(\Halb, f);\varepsilon^{2}}), ({(\Vier, a);\varepsilon^{4}}),
    ({(\Vier, g);\varepsilon^{3}}).
\end{split}
\end{equation}

 \vspace{-6pt} 

\begin{figure}[H]
		\centering
	\includegraphics[scale=0.7]{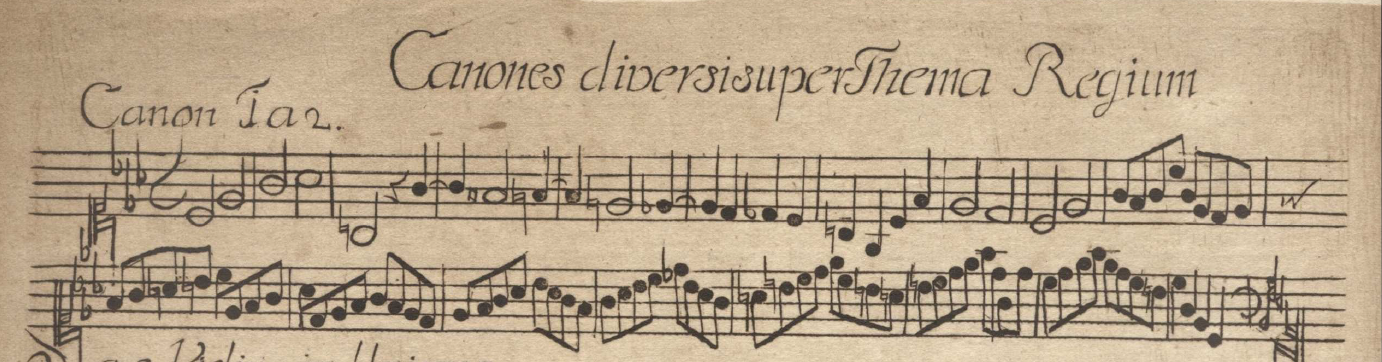}
	\caption{Brauer message of the Brauer configuration (\ref{A2Br}) gives Bach's canon 1 \^a2 \cite{Crab}. }
	\label{SLYmod}
\end{figure}

Reductions and segmentations of the Brauer message (\ref{MMA2}) give rise to  Bach's \textit{Canon 1, a2} shown in Figure \ref{SLYmod}. Words $\omega_{1}-\omega_{18}$ (see identities $(\ref{crB})$) define the $\mathfrak{M}$-reduced Brauer configuration $\mathfrak{M}^{A2}$.

\begin{equation}\label{MMA2}
\begin{split}
f^{32}\sigma_{-1}d^{32}b^{32}\sigma_{-1}a^{32}&\sigma_{1}\sigma_{-1}g^{32}c^{16}_{0}(b^{16}b^{16})\sigma_{-1} c^{32}(\sigma_{1}\sigma_{-1}c^{16}c^{16})\sigma_{1}\sigma_{-1}d^{32}(\sigma_{-1}d^{16}\\
\sigma_{-1} d^{16})e^{16}\sigma_{-1}e^{16}&f^{16}\sigma_{1}\sigma_{-1} g^{16}b^{16}f^{16}c^{16}\sigma_{-1} d^{32}e^{32}f^{32}\sigma_{-1}d^{32}[b^{8}c^{8}b^{8}f^{8}][b^{8}\sigma_{-1} d^{8}e^{8}\sigma_{-1}d^{8}] [c^{8}b^{8}\\
\sigma_{1}\sigma_{-1}a^{8}\sigma_{1}\sigma_{-1} g^{8}]&[f^{8}\sigma_{-1} d^{8}c^{8}b^{8}] [a^{8}e^{8}\sigma_{-1} d^{8}c^{8}][b^{8}c^{8}\sigma_{-1}d^{8}e^{8}] [\sigma_{-1}d^{8}c^{8}b^{8}\sigma_{-1}a^{8}][\sigma_{-1}g^{8}\sigma_{-1}a^{8}b^{8}c^{8}]\\
 [b^{8}\sigma_{-1}a^{8}\sigma_{-1}g^{8}f^{8}]&[\sigma_{-1}e^{8}\sigma_{-1}g^{8}\sigma_{-1}a^{8}b^{8}][\sigma_{1}\sigma_{-1}a^{8}\sigma_{1}\sigma_{-1}g^{8}f^{8}e^{8}][\sigma_{-1}d^{8}f^{8}\sigma_{1}\sigma_{-1}g^{8} \sigma_{1}\sigma_{-1}a^{8}]\\
 [\sigma_{1}\sigma_{-1}g^{8}f^{8}e^{8}\sigma_{-1}d^{8}]&[c^{8}e^{8}b^{8}e^{8}][f^{8}e^{8}\sigma_{-1}d^{8}c^{8}][\sigma_{-1}d^{8}e^{8}f^{8}\sigma_{1}\sigma_{-1}g^{8}] f^{8}b^{8}\sigma_{-1}d^{8}f^{8}.
 \end{split}
 \end{equation}

\begin{equation}\label{crB}
\begin{split}
(\omega_{1},\mathcal{K})&=f^{32}\sigma_{-1}d^{32}\longleftrightarrow(\Halb, f)(\flat\Halb, d)\\
(\omega_{2},\mathcal{K})&=b^{32}\sigma_{-1}a^{32}\longleftrightarrow(\Halb, b)(\flat\Halb, a)\\
(\omega_{3},\mathcal{K})&=\sigma_{1}\sigma_{-1}g^{32}c^{16}_{0}(b^{16}\longleftrightarrow(\natural\Halb, g)\ViPa((\Vier, b)\\
(\omega_{4},\mathcal{K})&=b^{16})\sigma_{-1} c^{32}(\sigma_{1}\sigma_{-1}c^{16}\longleftrightarrow (\Vier, b))(\flat\Halb, c)((\natural\Vier, c)\\
(\omega_{5},\mathcal{K})&=c^{16})\sigma_{1}\sigma_{-1}d^{32}(\sigma_{-1}d^{16}\longleftrightarrow (\Vier, c))(\natural\Halb, d)((\flat\Vier, d)\\
(\omega_{6},\mathcal{K})&=\sigma_{-1} d^{16})e^{16}\sigma_{-1}e^{16}f^{16}\longleftrightarrow (\flat\Vier, d))(\Vier, e)(\flat\Vier, e)(\Vier, f)\\
(\omega_{7},\mathcal{K})&=\sigma_{1}\sigma_{-1} g^{16}b^{16}f^{16}c^{16}\longleftrightarrow (\natural\Vier, g)(\Vier, b)(\Vier, f)(\Vier, c)\\
(\omega_{8},\mathcal{K})&=\sigma_{-1} d^{32}e^{32}\longleftrightarrow (\flat\Halb, d)(\Halb, e)\\
(\omega_{9},\mathcal{K})&=f^{32}\sigma_{-1}d^{32}\longleftrightarrow (\Halb, f)(\flat\Halb, d)\\
(\omega_{10},\mathcal{K})&= [b^{8}c^{8}b^{8}f^{8}][b^{8}\sigma_{-1} d^{8}e^{8}\sigma_{-1}d^{8}]\longleftrightarrow [(\eighthnote, b)(\eighthnote, c)(\eighthnote, b)(\eighthnote, f)] [(\eighthnote, b)(\flat\eighthnote, d)(\eighthnote, e)(\flat\eighthnote, d)]\\
(\omega_{11},\mathcal{K})&= [c^{8}b^{8}\sigma_{1}\sigma_{-1}a^{8}\sigma_{1}\sigma_{-1} g^{8}][f^{8}\sigma_{-1} d^{8}c^{8}b^{8}]\longleftrightarrow [(\eighthnote, c)(\eighthnote, b)(\natural\eighthnote, a)(\natural\eighthnote, g)]\\& [(\eighthnote, f)(\flat\eighthnote, d)(\eighthnote, c)(\flat\eighthnote, b)]\\
(\omega_{12},\mathcal{K})&= [a^{8}e^{8}\sigma_{-1} d^{8}c^{8}][b^{8}c^{8}\sigma_{-1}d^{8}e^{8}]\longleftrightarrow [(\eighthnote, a)(\eighthnote, e)(\flat\eighthnote, d)(\eighthnote, c)]\\& [(\eighthnote, b)(\eighthnote, c)(\flat\eighthnote, d)(\eighthnote, e)]\\
(\omega_{13},\mathcal{K})&= [\sigma_{-1}d^{8}c^{8}b^{8}\sigma_{-1}a^{8}][\sigma_{-1}g^{8}\sigma_{-1}a^{8}b^{8}c^{8}]\longleftrightarrow [(\flat\eighthnote, d)(\eighthnote, c)(\eighthnote, b)(\flat\eighthnote, a)]\\& [(\flat\eighthnote, g)(\flat\eighthnote, a)(\eighthnote, b)(\eighthnote, c)]\\
(\omega_{14},\mathcal{K})&= [b^{8}\sigma_{-1}a^{8}\sigma_{-1}g^{8}f^{8}][\sigma_{-1}e^{8}\sigma_{-1}g^{8}\sigma_{-1}a^{8}b^{8}]\longleftrightarrow [(\eighthnote, b)(\flat\eighthnote, a)(\flat\eighthnote, g)(\flat\eighthnote, g)(\eighthnote, f)]\\& [(\flat\eighthnote, e)(\flat\eighthnote, g)(\flat\eighthnote, a)(\eighthnote, b)]\\
(\omega_{15},\mathcal{K})&= [\sigma_{1}\sigma_{-1}a^{8}\sigma_{1}\sigma_{-1}g^{8}f^{8}e^{8}][\sigma_{-1}d^{8}f^{8}\sigma_{1}\sigma_{-1}g^{8}\sigma_{1}\sigma_{-1}a^{8}]\longleftrightarrow [(\natural\eighthnote, a)(\natural\eighthnote, g)(\eighthnote, f)(\eighthnote, e)]\\& [(\flat\eighthnote, d)(\eighthnote, f)(\natural\eighthnote, g)(\natural\eighthnote, a)]\\
(\omega_{16},\mathcal{K})&= [\sigma_{1}\sigma_{-1}g^{8}f^{8}e^{8}\sigma_{-1}d^{8}][c^{8}e^{8}b^{8}e^{8}]\longleftrightarrow [(\natural\eighthnote, g)(\eighthnote, f)(\eighthnote, e)(\flat\eighthnote, d)]\\& [(\eighthnote, c)(\eighthnote, e)(\eighthnote, b)(\eighthnote, e)]\\
(\omega_{17},\mathcal{K})&= [f^{8}e^{8}\sigma_{-1}d^{8}c^{8}][\sigma_{-1}d^{8}e^{8}f^{8}\sigma_{1}\sigma_{-1}g^{8}]\longleftrightarrow [(\eighthnote, f)(\eighthnote, e)(\flat\eighthnote, d)(\flat\eighthnote, c)]\\& [(\flat\eighthnote, d)(\eighthnote, e)(\eighthnote, f)(\natural\eighthnote, g)]\\
(\omega_{18},\mathcal{K})&= f^{8}b^{8}\sigma_{-1}d^{8}f^{8}\longleftrightarrow (\eighthnote, f)(\eighthnote, b)(\flat\eighthnote, d)(\eighthnote, f)
\end{split}
\end{equation}

The following data (\ref{A2SBr}) and (\ref{A2Br}) give the dimension of the Brauer configuration algebra (and its center) induced by Canon's Bach 1 a2 known as the \textit{crab canon} (see Figure \ref{SLYmod}).
\begin{equation}\label{A2SBr}
\begin{split}
S_{(\flat\Halb, a)}&=\omega_{1},\\
S_{(\Halb, b)}&=\omega_{1},\\
S_{(\flat\Halb, c)}&=\omega_{3},\\
S_{(\flat\Halb, d)}&=\omega_{1}<\omega_{7}<\omega_{8},\\
S_{(\Halb, e)}&=\omega_{7},\\
S_{(\flat\Halb, f)}&=\omega_{1}<\omega_{8},\\
S_{(\flat\Halb, f)}&=\omega_{1},\\
S_{(\Halb, g)}&=\omega_{2},\\
S_{\ViPa}&=\omega_{2},\\
S_{(\Vier, b)}&=\omega_{2}<\omega_{6},\\
S_{(\flat\Vier, d)}&=\omega_{4},\\
S_{(\Vier, e)}&=\omega_{5},\\
S_{(\flat\Vier, e)}&=\omega_{5},\\
S_{(\Vier, f)}&=\omega_{5}<\omega_{6},\\
S_{(\Vier, g)}&=\omega_{6},\\
S_{(\eighthnote, a)}&=\omega_{10}<\omega_{11}<\omega_{14}<\omega_{14}<\omega_{16},\\
S_{(\eighthnote, b)}&=\omega_{9}<\omega_{9}<\omega_{9}<\omega_{10}<\omega_{11}<\omega_{12}<\omega_{12}<\omega_{13}<\omega_{13}<\omega_{15}<\omega_{17},\\
S_{(\eighthnote, c)}&=\omega_{9}<\omega_{10}<\omega_{10}<\omega_{11}<\omega_{11}<\omega_{12}<\omega_{12}<\omega_{15}<\omega_{16},\\
S_{(\eighthnote, e)}&=\omega_{9}<\omega_{11}<\omega_{11}<\omega_{14}<\omega_{15}<\omega_{15}<\omega_{15}<\omega_{16}<\omega_{16},\\
S_{(\eighthnote, f)}&=\omega_{9}<\omega_{13}<\omega_{14}<\omega_{14}<\omega_{15}<\omega_{16}<\omega_{16}<\omega_{17}<\omega_{17},\\
S_{(\eighthnote, g)}&=\omega_{10}<\omega_{11}<\omega_{14}<\omega_{14}<\omega_{15},\\
S_{(\flat\eighthnote, a)}&=\omega_{12}<\omega_{12}<\omega_{13}<\omega_{13},\\
S_{(\flat\eighthnote, d)}&=\omega_{9}<\omega_{9}<\omega_{10}<\omega_{11}<\omega_{11}<\omega_{12}<\omega_{14}<\omega_{15}<\omega_{16}<\omega_{16}<\omega_{17},\\
S_{(\flat\eighthnote, e)}&=\omega_{13},\\
S_{(\flat\eighthnote, g)}&=\omega_{12}<\omega_{13}<\omega_{13}.\\
\end{split}
\end{equation}

\begin{equation}\label{A2Br}
\begin{split}
\mathfrak{M}^{A2}_{0}&=\{(\flat\Halb, a),(\Halb,b), (\flat\Halb, c), (\flat\Halb,d),(\Halb,d),(\Halb, e),(\Halb,f), (\Halb, g), \ViPa, (\Vier, b),\\&(\flat\Vier, d), (\Vier, e), (\flat\Vier, e),(\Vier, f),(\Vier, g),(\eighthnote, a),(\eighthnote, b),(\eighthnote, c),(\eighthnote, e),(\eighthnote, f),\\
&(\eighthnote, g),(\flat\eighthnote, a),\flat\eighthnote, d),(\flat\eighthnote, e),(\flat\eighthnote, g)\},\\
\mathfrak{M}^{A2}_{1}&=\{\omega_{1},\dots, \omega_{13}\},\\
\mathcal{K}&=(\bassclef,d<e<f<\dots<c<d, \text{4th line}),\quad\text{labels polygons in}\hspace{0.1cm}\mathfrak{M}^{A2},\\
\omega_{1}&<\omega_{2}<\dots<\omega_{13}\quad\text{for successor sequences,}\\
\mathrm{dim}_{k}\hspace{0.1cm}{\Lambda_{\mathfrak{M}^{A2}}}&=582,\\
|\{\alpha\in\mathfrak{M}^{A2}_{0}\mid val(\alpha)=1\}|&=12,\\
|\{\alpha\in\mathfrak{M}^{A2}_{0}\mid val(\alpha)=2\}|&=3,\\
|\{\alpha\in\mathfrak{M}^{A2}_{0}\mid val(\alpha)=3\}|&=2,\\
|\{\alpha\in\mathfrak{M}^{A2}_{0}\mid val(\alpha)=4\}|&=1,\\
|\{\alpha\in\mathfrak{M}^{A2}_{0}\mid val(\alpha)=5\}|&=2,\\
|\{\alpha\in\mathfrak{M}^{A2}_{0}\mid val(\alpha)=9\}|&=2,\\
|\{\alpha\in\mathfrak{M}^{A2}_{0}\mid val(\alpha)=11\}|&=3,\\
\#Loops(Q_{\mathfrak{M}^{A2}})&=32,\\
\mathrm{dim}_{k}\hspace{0.1cm}Z({\Lambda_{\mathfrak{M}^{A2}}})&=46.
\end{split}
\end{equation}

 \vspace{-6pt} 

\begin{equation}
\begin{split}
((\Halb, f), \varepsilon^{2}),&((\flat\Halb, d), \varepsilon^{0}),((\Halb, b), \varepsilon^{-2}),((\flat\Halb, a), \varepsilon^{-3}),((\natural\Halb, g), \varepsilon^{3}),\ViPa,((\Vier, b), \varepsilon^{-2}),\\
((\sharp\Halb,c),\varepsilon^{-1}),&((\natural\Vier, c), \varepsilon^{-1}),((\natural\Halb, d), \varepsilon^{0}),((\flat\Vier, d), \varepsilon^{0}),((\Vier, e), \varepsilon^{1}),((\flat\Vier, e), \varepsilon^{1}),((\Vier, f), \varepsilon^{2}),\\
((\natural\Vier, g), \varepsilon^{3}),&((\Halb, e), \varepsilon^{1}),((\eighthnote, b), \varepsilon^{-2}),((\eighthnote, c), \varepsilon^{-1}),((\eighthnote, f), \varepsilon^{-5}),((\flat\eighthnote, d), \varepsilon^{0}),((\eighthnote, e), \varepsilon^{1}),\\
((\natural\eighthnote, a), \varepsilon^{-3}),&((\natural\eighthnote, g), \varepsilon^{-4}),((\flat\eighthnote, a), \varepsilon^{-3}),((\flat\eighthnote, g), \varepsilon^{-4}),((\flat\eighthnote, e), \varepsilon^{-7}).
\end{split}
\end{equation}

 \vspace{-6pt} 
  \begin{figure}[H]
		\centering
		\includegraphics[scale=0.4]{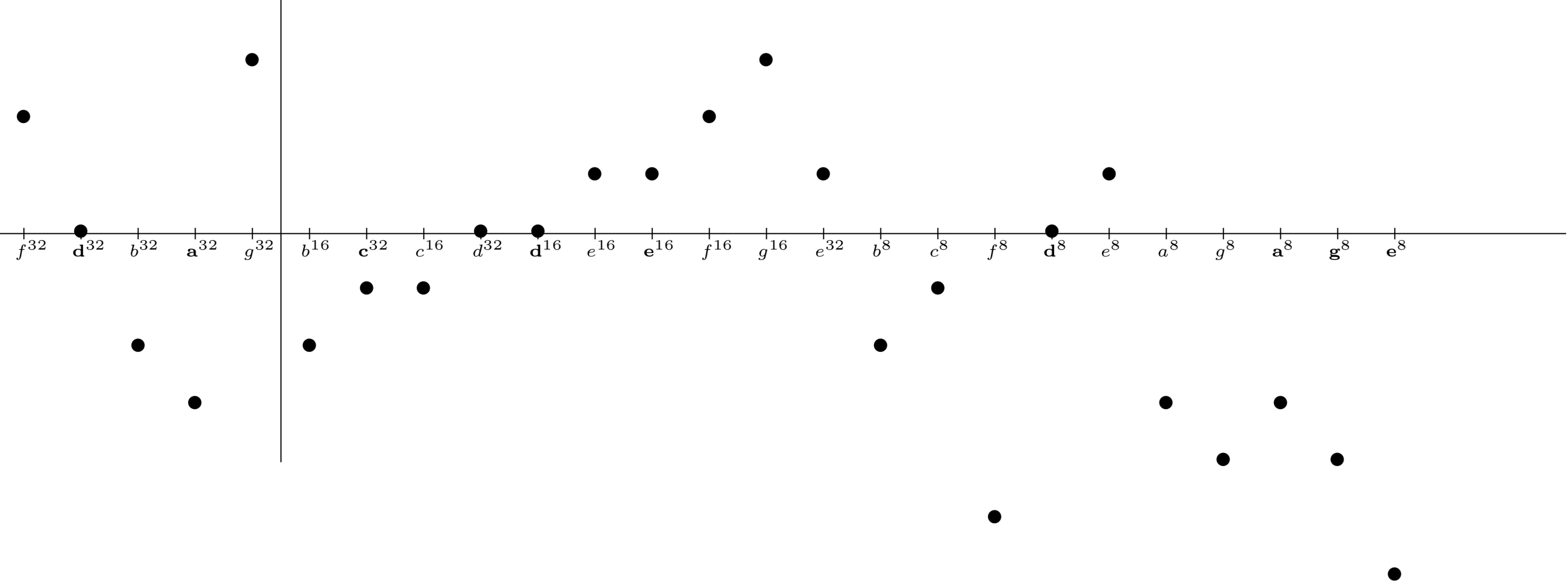}\\
	\includegraphics[scale=0.4]{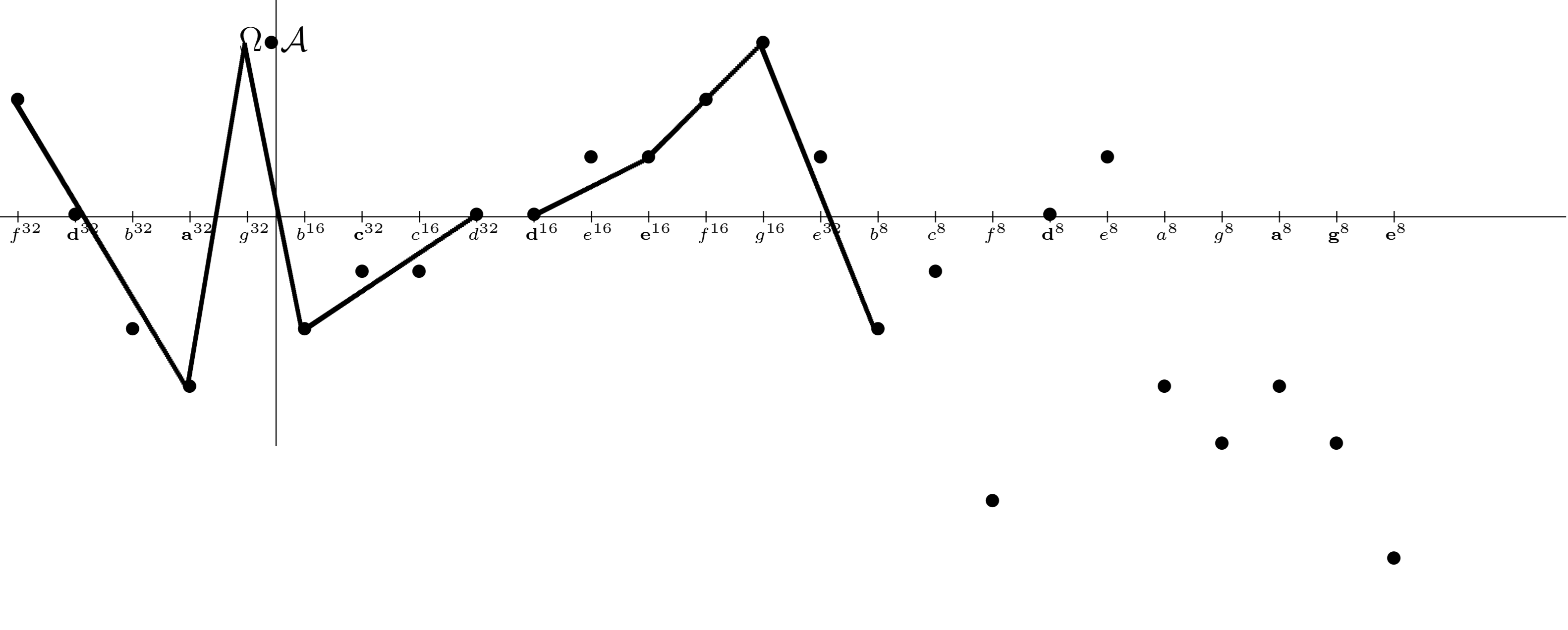}	\\
	\includegraphics[scale=0.4]{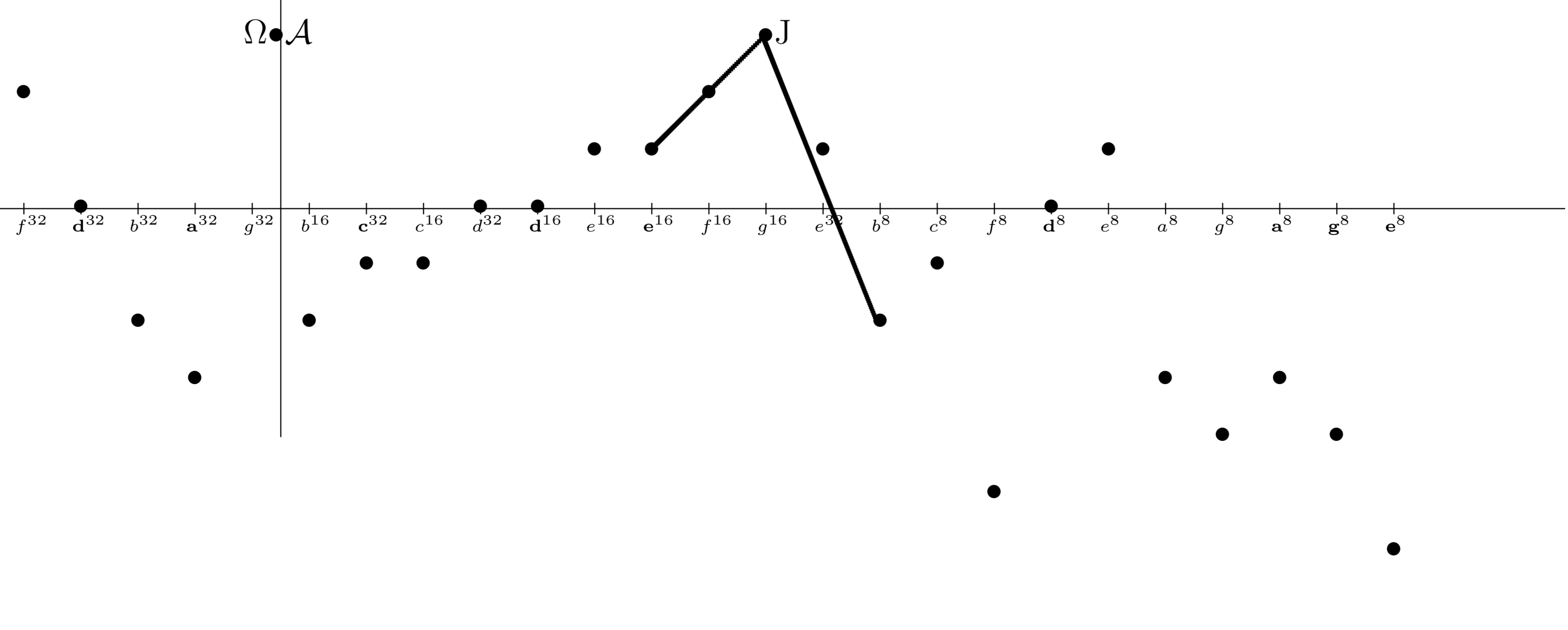}\\
	\includegraphics[scale=0.4]{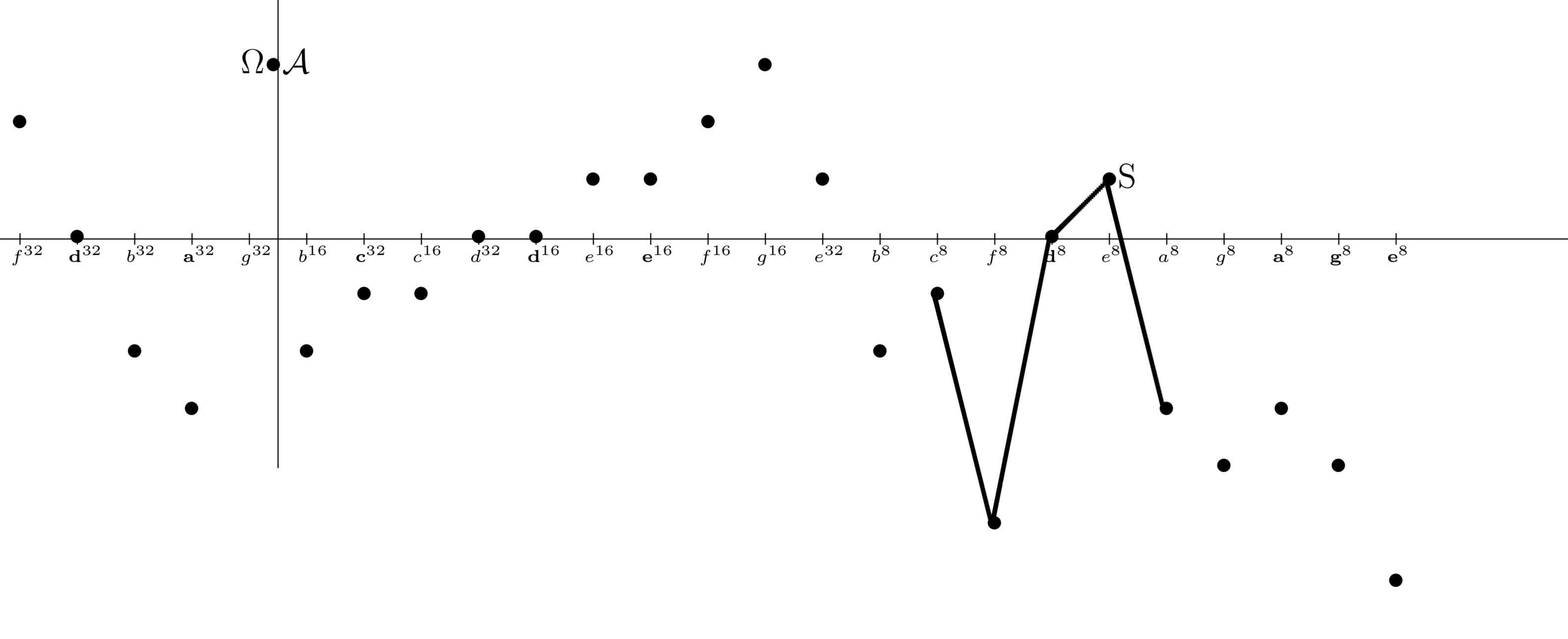}\\	
	\includegraphics[scale=0.4]{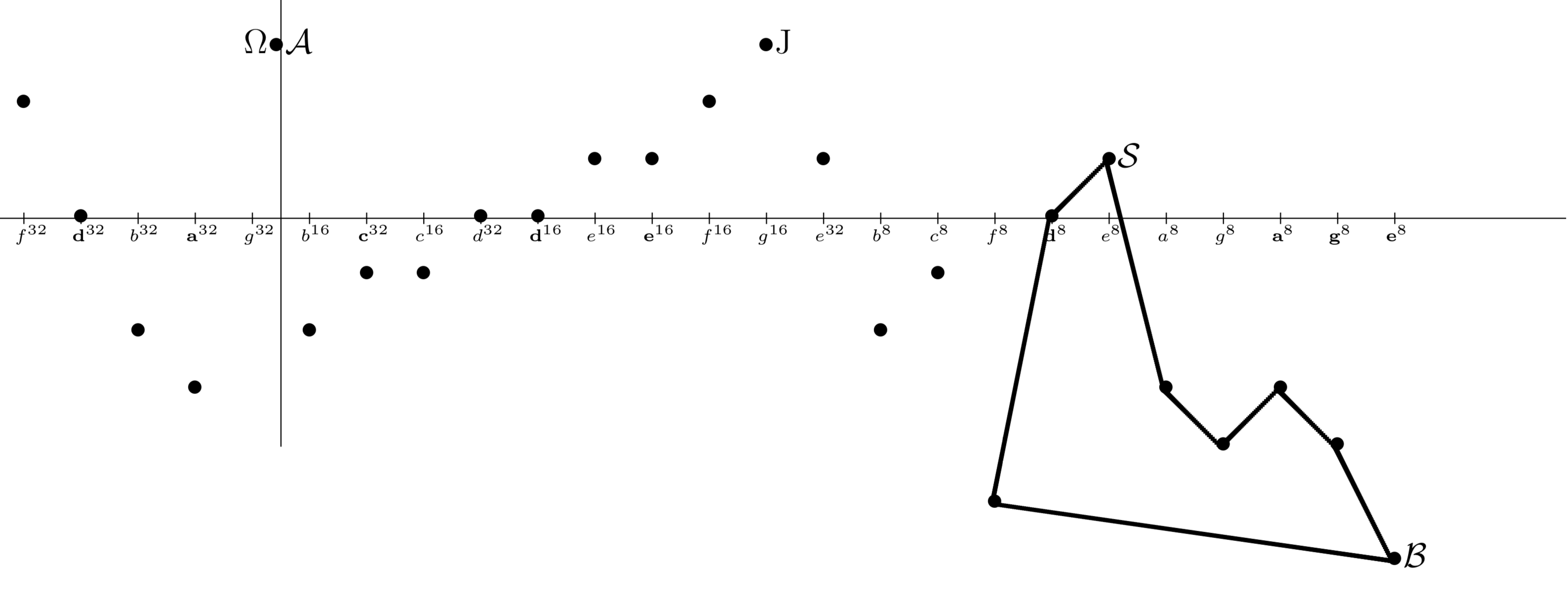}\\
	\caption{Letters, W, A, J,S , and B as they arise by  connecting consecutive points defined by Bach's canon 1 \^a 2 (see Algorithm \ref{Alg}). }
	\label{a6V}
\end{figure}

  Canon's Bach a4 \textit{Quaerendo Invenietis} is given by the $\mathfrak{M}$-reduced Brauer message of the Brauer configuration $\mathfrak{M}^{QI}=(\mathfrak{M}^{QI}_{0},\mathfrak{M}^{QI}_{1},\mu,\mathcal{O})$ defined by words (\ref{QIW}) which were obtained after appropriate segmentations of the Brauer message (\ref{MMQ0I}).
  
  \newpage
  \begin{equation}\label{QIW}
  \begin{split}
  \omega_{1}&=[e^{8}f^{8}][\sigma_{-1}g^{8}a^{8}b^{8}\sigma_{-1}c^{8}]\longleftrightarrow [(\eighthnote, e)(\eighthnote, f)][(\flat\eighthnote, g)(\eighthnote, a)(\eighthnote, b)(\flat\eighthnote, c)]\\  
  \omega_{2}&=\sigma_{1}d^{16}c^{16}_{0}c^{16}_{0}(b^{16}\longleftrightarrow (\sharp\Vier, d)\ViPa\ViPa((\Vier, b)\\  
  \omega_{3}&=[b^{8})\sigma_{1}a^{8}f^{8}\sigma_{1}\sigma_{-1}g^{8}](\sigma_{-1}\sigma_{1}a^{32}\longleftrightarrow [(\eighthnote, b))(\sharp\eighthnote, a)(\eighthnote, f)(\natural\eighthnote, g)](\natural\Halb, a)\\  
  \omega_{4}&=[a^{8})\sigma_{1}\sigma_{-1}g^{8}e^{8}f^{8}](\sigma_{-1}g^{32}\longleftrightarrow [(\eighthnote, a))(\sharp\eighthnote, g)(\eighthnote, e)(\eighthnote, f)](\flat\Halb, g)\\  
   \omega_{5}&=[(\sigma_{-1}g^{8})f^{8}\sigma_{1}e^{8}f^{8})][(\sigma_{-1}f^{8}e^{8}\sigma_{1}d^{8}e^{8})]\longleftrightarrow [((\flat\eighthnote, g))(\eighthnote, f)(\sharp\eighthnote, e)(\eighthnote, f))]\\&[((\flat\eighthnote, f)(\eighthnote,e)(\eighthnote,d)(\sharp\eighthnote,g)(\eighthnote,e))]\\  
    \omega_{6}&=[(\sigma_{1}d^{8}\sigma_{1}\sigma_{-1}c^{8}b^{8}\sigma_{-1}c^{8})][(b^{8}\sigma_{1}d^{8}f^{8}a^{8}]\longleftrightarrow [((\sharp\eighthnote, d)(\natural\eighthnote, c)(\eighthnote, b)(\flat\eighthnote,c))]\\&[((\eighthnote, b)(\sharp\eighthnote,d)(\eighthnote, f)(\eighthnote,a))]\\  
    \omega_{7}&=[(\sigma_{-1}g^{8}b^{8}\sigma_{1}d^{8}e^{8})][f^{8}\sigma_{-1}c^{8}b^{8}a^{8}]\longleftrightarrow [((\flat\eighthnote, g)(\eighthnote, b)(\sharp\eighthnote, d)(\eighthnote, e))]\\&[(\eighthnote, f)(\flat\eighthnote, c)(\eighthnote, b)(\eighthnote,a)]\\ 
   \omega_{8}&=[\sigma_{-1}g^{4}f^{4}e^{4}\sigma_{1}d^{4}][e^{4}e^{4}\sigma_{1}d^{4}\sigma_{-1}c^{4}][b^{4}a^{4}\sigma_{-1}g^{4}f^{4}]e^{16}\longleftrightarrow[(\flat\text{\Sech},g)(\text{\Sech},f)(\text{\Sech},e)(\sharp\text{\Sech},d)]\\&[(\text{\Sech},e)(\text{\Sech},e)(\sharp\text{\Sech},d)(\flat\text{\Sech},c)][(\text{\Sech},b)(\text{\Sech},a)(\text{\flat\Sech},g)(\text{\Sech},f)](\Vier, e)\\
   \omega_{9}&=d^{8}_{0}[f^{8}\sigma_{-1}g^{8}a^{8}][b^{8}d^{8}\sigma_{1}\sigma_{-1}c^{8}e^{8}]\longleftrightarrow \ViPa[(\eighthnote, f)(\flat\eighthnote, g)(\eighthnote, a)]\\&[(\eighthnote, b)(\eighthnote, d)(\natural\eighthnote, c)(\eighthnote, e)]\\ 
   \omega_{10}&=[d^{8}\sigma_{1}\sigma_{-1}c^{4}b^{4}]\sigma_{1}a^{16}d^{8}_{0}[\sigma_{-1}c^{8}b^{8}d^{8}]\longleftrightarrow [(\eighthnote, d)(\natural\eighthnote, c)(\text{\Sech}, b)](\sharp\Vier, a)\AcPa\\&[(\flat\eighthnote, c)(\eighthnote, b)(\eighthnote, d)]\\ 
   \omega_{11}&=[\sigma_{-1}c^{8}b^{4}a^{4}]\sigma_{1}\sigma_{-1}g^{16}d^{8}_{0}[b^{8}\sigma_{1}a^{8}\sigma_{1}\sigma_{-1}c^{8}]\longleftrightarrow [(\flat\eighthnote, c)(\text{\Sech}, b)(\text{\Sech}, a)](\natural\Vier, g)\AcPa\\&[(\eighthnote, b)(\sharp\eighthnote, a)(\natural\eighthnote, c)]\\ 
  \omega_{12}&=[b^{8}\sigma_{1}\sigma_{-1}a^{8}\sigma_{1}\sigma_{-1}g^{8}b^{8}][a^{8}\sigma_{1}\sigma_{-1}g^{8}](a^{16}\longleftrightarrow [(\eighthnote, b)(\natural\eighthnote, a)(\natural\eighthnote, g)(\eighthnote, b)]\\&[(\eighthnote, a)(\natural\eighthnote, g)]((\Vier, a)\\
  \omega_{13}&=[a^{8})e^{8}\sigma_{1}d^{8}e^{8}][f^{8}b^{8}](b^{16}\longleftrightarrow [(\eighthnote, a))(\eighthnote, e)(\sharp\eighthnote, d)(\eighthnote, e)][(\eighthnote, f))(\eighthnote, b)]((\Vier,b)\\    
  \end{split}
  \end{equation}

  \begin{equation}
  \begin{split}
   \omega_{14}&=b^{16})[a^{8}\sigma_{-1}g^{8}]a^{16}[\sigma_{-1}g^{8}f^{8}]\longleftrightarrow (\Vier, b))[(\eighthnote, a)(\flat\eighthnote, g)](\Vier, a)[(\flat\eighthnote, g)(\eighthnote, f)]\\  
    \omega_{15}&=[\sigma_{-1}g^{4}a^{4}\sigma_{-1}g^{4}f^{4}](e^{16}[e^{8})f^{8}\sigma_{-1}g^{8}a^{8}]\longleftrightarrow [(\flat\text{\Sech}, g)(\text{\Sech}, a)(\flat\text{\Sech}, g)(\text{\Sech}, f)]((\Vier, e)\\&[(\eighthnote, e))(\eighthnote, f)(\flat\eighthnote, g)(\eighthnote, a)]\\  
    \omega_{16}&=[b^{8}\sigma_{-1}c^{8}b^{8}a^{8}][\sigma_{-1}g^{8}f^{8}e^{8}\sigma_{-1}g^{8}]\longleftrightarrow [(\eighthnote, b)(\flat\eighthnote, c)(\eighthnote, b)(\eighthnote, a)]\\&[(\flat\eighthnote, g))(\eighthnote, f)(\eighthnote, e)(\flat\eighthnote, g)]\\  
     \omega_{17}&=(f^{32}[f^{8})e^{8}\sigma_{1}d^{8}f^{8}]\longleftrightarrow((\Halb, f) [(\eighthnote, f))(\eighthnote, e)(\sharp\eighthnote, d)(\eighthnote, f)]\\  
     \omega_{18}&=(e^{32}[e^{8})d^{8}\sigma_{1}\sigma_{-1}c^{8}e^{8}]\longleftrightarrow((\Halb, e) [(\eighthnote, e))(\eighthnote, d)(\natural\eighthnote, c)(\eighthnote, e)]\\  
      \omega_{19}&=d^{16}d^{8}_{0}d^{8}a^{16}d^{8}_{0}\sigma_{-1}c^{8}\longleftrightarrow(\Vier, d) \AcPa(\eighthnote, d)(\Vier, a)\AcPa(\flat\eighthnote, c)\\  
        \omega_{20}&=b^{16}d^{8}_{0}e^{8}\sigma_{1}d^{8}d^{8}_{0}b^{8}\longleftrightarrow(\Vier, b) \AcPa(\eighthnote, e)(\sharp\eighthnote, d)\AcPa(\eighthnote, b)\\  
         \omega_{21}&=e^{16}c^{16}a^{16}b^{16}\longleftrightarrow(\Vier, e) (\Vier, c)(\Vier, a)(\Vier, b)\\  
         \omega_{22}&=[e^{8}e^{4}f^{4}]\varphi_{8}b^{16}[a^{4}\sigma_{-1}g^{4}]\longleftrightarrow[(\eighthnote, e) (\text{\Sech}, e)(\text{\Sech}, f)](\underset{\cdot}{\Vier}, b)[(\text{\Sech}, a)(\flat\text{\Sech}, g)]\\  
         \omega_{23}&=[f^{8}\sigma_{1}d^{8}e^{8}f^{8}]b^{16}c^{16}_{0}\longleftrightarrow[(\eighthnote, f) (\sharp\eighthnote, d)(\eighthnote, e)(\eighthnote, f)](\Vier, b)\ViPa\\  
         \omega_{24}&=d^{8}_{0}[\sigma_{1}\sigma_{-1}c^{8}\sigma_{1}\sigma_{-1}d^{8}e^{8}]f^{32}\longleftrightarrow\AcPa[(\natural\eighthnote, c) (\natural\eighthnote, d)(\eighthnote, e)](\Halb, f)\\  
         \omega_{25}&=d^{8}_{0}[b^{8}\sigma_{1}\sigma_{-1}c^{8}\sigma_{1}d^{8}]e^{32}\longleftrightarrow\AcPa[(\eighthnote, b) (\natural\eighthnote, c)(\sharp\eighthnote, d)](\Halb, e)\\ 
          \omega_{26}&=c^{16}_{0}b^{16}\sigma_{-1}c^{16}f^{16}\longleftrightarrow\ViPa(\Vier, b)(\flat\Vier, c)(\Vier, f)\\   
          \omega_{27}&=[f^{4}\sigma_{-1}g^{4}]a^{16}\sigma_{-1}g^{8}[a^{8}f^{8}]\sigma_{1}d^{16}\longleftrightarrow [(\text{\Sech}, f)(\flat\text{\Sech}, g)](\Vier, a)(\flat\eighthnote,g)\\&[(\eighthnote, a)(\eighthnote, f)](\sharp\Vier, d)\\ 
           \omega_{28}&=e^{16}d_{0}^{8}\sigma_{-1}g^{8}\sigma_{-1}c^{8}d^{8}_{0}\sigma_{1}d^{8}\longleftrightarrow(\Vier, e) \AcPa(\flat\eighthnote, g)(\flat\eighthnote, c)\AcPa(\sharp\eighthnote, d)\\                  
 \end{split}
  \end{equation}

  \begin{equation}
  \begin{split}
S_{(\eighthnote,e)}&=\omega_{1}<\omega_{4}<\omega_{5}<\omega_{5}<\omega_{7}<\omega_{9}<\omega_{13}<\omega_{13}<\omega_{15}<\omega_{16}<\omega_{17}<\\&\omega_{18}<\omega_{18}<\omega_{20}<\omega_{22}<\omega_{23}<\omega_{24},\quad val((\eighthnote, e))=17.\\
S_{(\eighthnote,f)}&=\omega_{1}<\omega_{3}<\omega_{4}<\omega_{5}<\omega_{5}<\omega_{6}<\omega_{7}<\omega_{9}<\omega_{13}<\omega_{14}<\omega_{15}\\&\omega_{16}<\omega_{17}<\omega_{17},\quad val((\eighthnote, f))=14.\\
S_{(\flat\eighthnote,g)}&=\omega_{1}<\omega_{5}<\omega_{7}<\omega_{9}<\omega_{13}<\omega_{14}<\omega_{14}<\omega_{15}<\omega_{16}<\omega_{16}<\\&\omega_{22}<\omega_{27}<\omega_{28},\quad val((\eighthnote, e))=13.\\
S_{(\eighthnote,a)}&=\omega_{1}<\omega_{4}<\omega_{6}<\omega_{7}<\omega_{9}<\omega_{12}<\omega_{12}<\omega_{13}<\omega_{14}<\omega_{15}<\omega_{16}<\\&\omega_{22}<\omega_{27},\quad val((\eighthnote, a))=13.\\
S_{(\eighthnote,b)}&=\omega_{1}<\omega_{3}<\omega_{6}<\omega_{6}<\omega_{7}<\omega_{7}<\omega_{9}<\omega_{10}<\omega_{11}<\omega_{12}<\omega_{12}\\&<\omega_{13}<\omega_{16}<\omega_{16}<\omega_{20}<\omega_{25},\quad val((\eighthnote, b))=16.\\
S_{(\flat\eighthnote,c)}&=\omega_{1}<\omega_{5} <\omega_{6}<\omega_{7}<\omega_{11} <\omega_{16}<\omega_{19},\quad val((\flat\eighthnote, c))=8.\\
S_{(\sharp\Vier,d)}&=\omega_{2}<\omega_{20}<\omega_{27},\quad val((\sharp\Vier, d))=3.\\
S_{\ViPa}&=\omega_{2}<\omega_{2} <\omega_{23}<\omega_{26},\quad val(\ViPa)=4.\\
S_{(\Vier,b)}&=\omega_{2}<\omega_{13}<\omega_{14}<\omega_{20}<\omega_{21}<\omega_{23}<\omega_{26},\quad val((\Vier, b))=7.\\
S_{(\sharp\eighthnote,a)}&=\omega_{3},\quad val((\sharp\Vier, d))=1.\\
S_{(\natural\eighthnote,g)}&=\omega_{3}<\omega_{4}<\omega_{12}<\omega_{12}, \quad val((\natural\eighthnote, g))=4.\\
S_{(\Halb,a)}&=\omega_{3},\quad val((\Halb, a))=1.\\
S_{(\flat\Halb, g)}&=\omega_{4},\quad val((\flat\Halb, g))=1.\\
S_{(\sharp\eighthnote,e)}&=\omega_{5},\quad val((\sharp\eighthnote, e))=1.\\
S_{(\flat\eighthnote, f)}&=\omega_{5},\quad val((\flat\eighthnote, f))=1.\\
S_{(\sharp\eighthnote,d)}&=\omega_{5}<\omega_{6}<\omega_{6}<\omega_{7}<\omega_{13}<\omega_{17}<\omega_{23}<\omega_{25}<\omega_{28},\quad val((\sharp\eighthnote, d))=9.
\end{split}
\end{equation}
\begin{equation}
\begin{split}
S_{(\eighthnote,c)}&=\omega_{6}<\omega_{9}<\omega_{11}<\omega_{18}<\omega_{24}<\omega_{25},\quad val((\eighthnote, c))=6.\\
S_{(\flat\text{\Sech},g)}&=\omega_{8}<\omega_{8}<\omega_{15}<\omega_{15}<\omega_{22}<\omega_{27},\quad val((\flat\text{\Sech},g))=6.\\
S_{(\text{\Sech},f)}&=\omega_{8}<\omega_{8}<\omega_{15}<\omega_{22}<\omega_{27},\quad val((\flat\text{\Sech},g))=5.\\
S_{(\text{\Sech},e)}&=\omega_{8}<\omega_{8}  <\omega_{8}<\omega_{22},\quad val((\flat\text{\Sech},g))=4.\\
S_{(\sharp\text{\Sech},d)}&=\omega_{8}<\omega_{8},\quad val((\text{\sharp\Sech},d))=2.\\
S_{(\flat\text{\Sech},c)}&=\omega_{8},\quad val((\flat\text{\Sech},d))=1.\\
S_{(\text{\Sech},b)}&=\omega_{8}<\omega_{10}<\omega_{11},\quad val((\text{\Sech},b))=3.\\
S_{(\text{\Sech},a)}&=\omega_{8}<\omega_{11}<\omega_{15}<\omega_{22},\quad val((\text{\Sech},a))=4.\\
S_{(\Vier,e)}&=\omega_{8}<\omega_{15}<\omega_{21}<\omega_{28},\quad val((\Vier, b))=4.\\
S_{\AcPa}&=\omega_{9}<\omega_{10}<\omega_{11}<\omega_{19}<\omega_{19}<\omega_{20}<\omega_{20}<\omega_{24}<\omega_{25}<\\&
\omega_{28}<\omega_{28},\quad val(\AcPa)=11.\\
S_{(\eighthnote, d)}&=\omega_{9}<\omega_{10} <\omega_{10}<\omega_{18}<\omega_{19}<\omega_{24},\quad val((\flat\eighthnote, c))=6.\\
S_{(\text{\Sech},c)}&=\omega_{10},\quad val((\text{\Sech},c))=1.\\
S_{(\sharp\Vier, a)}&=\omega_{10},\quad val((\sharp\Vier, a))=1.\\
S_{(\Vier, g)}&=\omega_{11},\quad val((\Vier, g))=1.\\
S_{(\sharp\eighthnote, a)}&=\omega_{11},\quad val((\sharp\eighthnote, a))=1.\\
S_{(\Vier, a)}&=\omega_{12}<\omega_{14}<\omega_{19}<\omega_{21}<\omega_{27},\quad val((\Vier, a))=5.\\
S_{(\Halb,f)}&=\omega_{17}<\omega_{24},\quad val((\Halb, f))=2.\\
S_{(\Halb,e)}&=\omega_{18}<\omega_{25},\quad val((\Halb, e))=2.\\
S_{(\Vier, d)}&=\omega_{19},\quad val((\Vier, d))=1.\\
S_{(\flat\Vier, c)}&=\omega_{21}<\omega_{28},\quad val((\flat\Vier, c))=2.\\
S_{(\underset{\cdot}{\Vier}, b)}&=\omega_{22},\quad val((\underset{\cdot}{\Vier}, b))=1.\\
S_{(\Vier, f)}&=\omega_{26},\quad val((\Vier, f))=1.\\
\end{split}
\end{equation}

The following Brauer message (\ref{MMQ0I}) gives rise to Bach's canon Quaerendo Invenietis.

\begin{equation}\label{MMQ0I}
\begin{split}
[e^{8}f^{8}]&[\sigma_{-1}g^{8}a^{8}b^{8}\sigma_{-1}c^{8}]\sigma_{1}d^{16}c^{16}_{0}c^{16}_{0}(b^{16}[b^{8})\sigma_{1}a^{8}f^{8}\sigma_{1}\sigma_{-1}g^{8}](\sigma_{-1}\sigma_{1}a^{32}\\
[a^{8})\sigma_{1}\sigma_{-1}g^{8}e^{8}f^{8}]&(\sigma_{-1}g^{32}[(\sigma_{-1}g^{8})f^{8}\sigma_{1}e^{8}f^{8})][(\sigma_{-1}f^{8}e^{8}\sigma_{1}d^{8}e^{8})][(\sigma_{1}d^{8}\sigma_{1}\sigma_{-1}c^{8}b^{8}\sigma_{-1}c^{8})]\\
[(b^{8}\sigma_{1}d^{8}f^{8}a^{8}]&[(\sigma_{-1}g^{8}b^{8}\sigma_{1}d^{8}e^{8})][f^{8}\sigma_{-1}c^{8}b^{8}a^{8}][\sigma_{-1}g^{4}f^{4}e^{4}\sigma_{1}d^{4}][e^{4}e^{4}\sigma_{1}d^{4}\sigma_{-1}c^{4}]\\
[b^{4}a^{4}\sigma_{-1}g^{4}f^{4}]&e^{16}d^{8}_{0}[f^{8}\sigma_{-1}g^{8}a^{8}][b^{8}d^{8}\sigma_{1}\sigma_{-1}c^{8}e^{8}][d^{8}\sigma_{1}\sigma_{-1}c^{4}b^{4}]\sigma_{1}a^{16}d^{8}_{0}[\sigma_{-1}c^{8}b^{8}d^{8}]\\
[\sigma_{-1}c^{8}b^{4}a^{4}]&\sigma_{1}\sigma_{-1}g^{16}d^{8}_{0}[b^{8}\sigma_{1}a^{8}\sigma_{1}\sigma_{-1}c^{8}[b^{8}\sigma_{1}\sigma_{-1}a^{8}\sigma_{1}\sigma_{-1}g^{8}b^{8}][a^{8}\sigma_{1}\sigma_{-1}g^{8}](a^{16}\\
[a^{8})e^{8}\sigma_{1}d^{8}e^{8}]&[f^{8}b^{8}](b^{16}b^{16})[a^{8}\sigma_{-1}g^{8}]a^{16}[\sigma_{-1}g^{8}f^{8}][\sigma_{-1}g^{4}a^{4}\sigma_{-1}g^{4}f^{4}](e^{16}[e^{8})f^{8}\\
\sigma_{-1}g^{8}a^{8}]&[b^{8}\sigma_{-1}c^{8}b^{8}a^{8}][\sigma_{-1}g^{8}f^{8}e^{8}\sigma_{-1}g^{8}](f^{32}[f^{8})e^{8}\sigma_{1}d^{8}f^{8}](e^{32}[e^{8})d^{8}\sigma_{1}\sigma_{-1}c^{8}e^{8}]\\
d^{16}d^{8}_{0}d^{8}a^{16}&d^{8}_{0}\sigma_{-1}c^{8}b^{16}d^{8}_{0}e^{8}\sigma_{1}d^{8}d^{8}_{0}b^{8}e^{16}c^{16}a^{16}b^{16}[e^{8}e^{4}f^{4}]\varphi_{8}b^{16}[a^{4}\sigma_{-1}g^{4}]\\
[f^{8}\sigma_{1}d^{8}e^{8}f^{8}]&b^{16}c^{16}_{0}d^{8}_{0}[\sigma_{1}\sigma_{-1}c^{8}\sigma_{1}\sigma_{-1}d^{8}e^{8}]f^{32}d^{8}_{0}[b^{8}\sigma_{1}\sigma_{-1}c^{8}\sigma_{1}d^{8}]e^{32}c^{16}_{0}b^{16}\sigma_{-1}c^{16}f^{16}\\
[f^{4}\sigma_{-1}g^{4}]&a^{16}\sigma_{-1}g^{8}[a^{8}f^{8}]\sigma_{1}d^{16}e^{16}d_{0}^{8}\sigma_{-1}g^{8}\sigma_{-1}c^{8}d^{8}_{0}\sigma_{1}d^{8}.
 \end{split}
 \end{equation}   
  
Note that,

\begin{equation}\label{QIBr}
\begin{split}
\mathfrak{M}^{QI}_{0}&=\mathfrak{M}^{QI,4}_{0}\cup\mathfrak{M}^{QI,8}_{0}\cup\mathfrak{M}^{QI,16}_{0}\cup\mathfrak{M}^{QI,32}_{0}\cup\mathfrak{M}^{QI,R}_{0},\\
\mathfrak{M}^{QI,4}_{0}&=\{(\text{\Sech}, a), (\text{\Sech}, c),(\text{\Sech}, d), (\text{\Sech}, f),(\text{\Sech}, g),(\sharp\text{\Sech}, b),(\flat\text{\Sech}, e) \},\\
\mathfrak{M}^{QI,8}_{0}&=\{(\eighthnote, a), (\eighthnote, b), (\eighthnote, c), (\eighthnote, d),(\eighthnote, e),(\eighthnote, f),(\eighthnote, g),(\flat\eighthnote, a),(\sharp\eighthnote, b),\\&(\sharp\eighthnote, c),(\flat\eighthnote, d), (\flat\eighthnote, e),  (\sharp\eighthnote, f), (\flat\eighthnote, f) \},\\
\mathfrak{M}^{QI,16}_{0}&=\{(\Vier, a), (\Vier, b),(\Vier, c),(\Vier, d), (\Vier, e),(\Vier, f),(\Vier, g),\\&(\flat\Vier, a), (\sharp\Vier, b),(\sharp\Vier, f),(\underset{\cdot}{\Vier}, g)\},\\
\mathfrak{M}^{QI,32}_{0}&=\{(\Halb, c),(\Halb, d),(\Halb, e),(\Halb, f) \},\\
\mathfrak{M}^{QI,R}_{0}&=\{\AcPa, \ViPa\},\\
\mathfrak{M}^{QI}_{1}&=\{\omega_{1},\dots, \omega_{28}\} \quad\text{(see, identities (\ref{QIW}))},\\
\mathcal{K}&=(\trebleclef, e<f<\dots<d<e,\text{1st line}),\hspace{0.1cm}\text{labels polygons in}\hspace{0.1cm}\mathfrak{M}^{QI},\\
\end{split}
\end{equation}
\begin{equation}
\begin{split}
\{\alpha\in\mathfrak{M}^{QI}_{0}\mid val(\alpha)=1\}&=13\\
|\{\alpha\in\mathfrak{M}^{A2}_{0}\mid val(\alpha)=2\}|&=4,\\
|\{\alpha\in\mathfrak{M}^{A2}_{0}\mid val(\alpha)=3\}|&=2,\\
|\{\alpha\in\mathfrak{M}^{A2}_{0}\mid val(\alpha)=4\}|&=5,\\
|\{\alpha\in\mathfrak{M}^{A2}_{0}\mid val(\alpha)=5\}|&=2,\\
|\{\alpha\in\mathfrak{M}^{A2}_{0}\mid val(\alpha)=6\}|&=3,\\
|\{\alpha\in\mathfrak{M}^{A2}_{0}\mid val(\alpha)=7\}|&=1,\\
|\{\alpha\in\mathfrak{M}^{A2}_{0}\mid val(\alpha)=8\}|&=1,\\
|\{\alpha\in\mathfrak{M}^{A2}_{0}\mid val(\alpha)=9\}|&=1,\\
|\{\alpha\in\mathfrak{M}^{A2}_{0}\mid val(\alpha)=11\}|&=1,\\
|\{\alpha\in\mathfrak{M}^{A2}_{0}\mid val(\alpha)=13\}|&=2,\\
|\{\alpha\in\mathfrak{M}^{A2}_{0}\mid val(\alpha)=14\}|&=1,\\
|\{\alpha\in\mathfrak{M}^{A2}_{0}\mid val(\alpha)=16\}|&=1,\\
|\{\alpha\in\mathfrak{M}^{A2}_{0}\mid val(\alpha)=17\}|&=1,\\
\#Loops\hspace{0.1cm}Q_{\mathfrak{M}^{QI}}&=38,\\
\mathrm{dim}_{k}\hspace{0.1cm}\Lambda_{\mathfrak{M}^{QI}}&=1565,\\
\mathrm{dim}_{k}\hspace{0.1cm}Z(\Lambda_{\mathfrak{M}^{QI}})&=67.
\end{split}
\end{equation}

    \vspace{-6pt} 
  \begin{figure}[H]
	\includegraphics[scale=1.0]{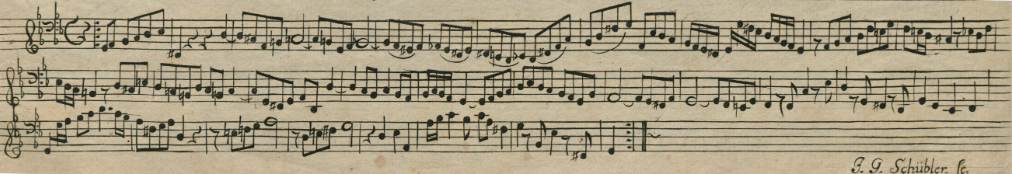}
	\caption{Brauer message (\ref{MMQ0I}) gives rise to Bach's canon \^a4 Quaerendo Invenietis. }
	\label{QI}
\end{figure}

\begin{figure}[H]
		\centering
		\includegraphics[scale=0.35]{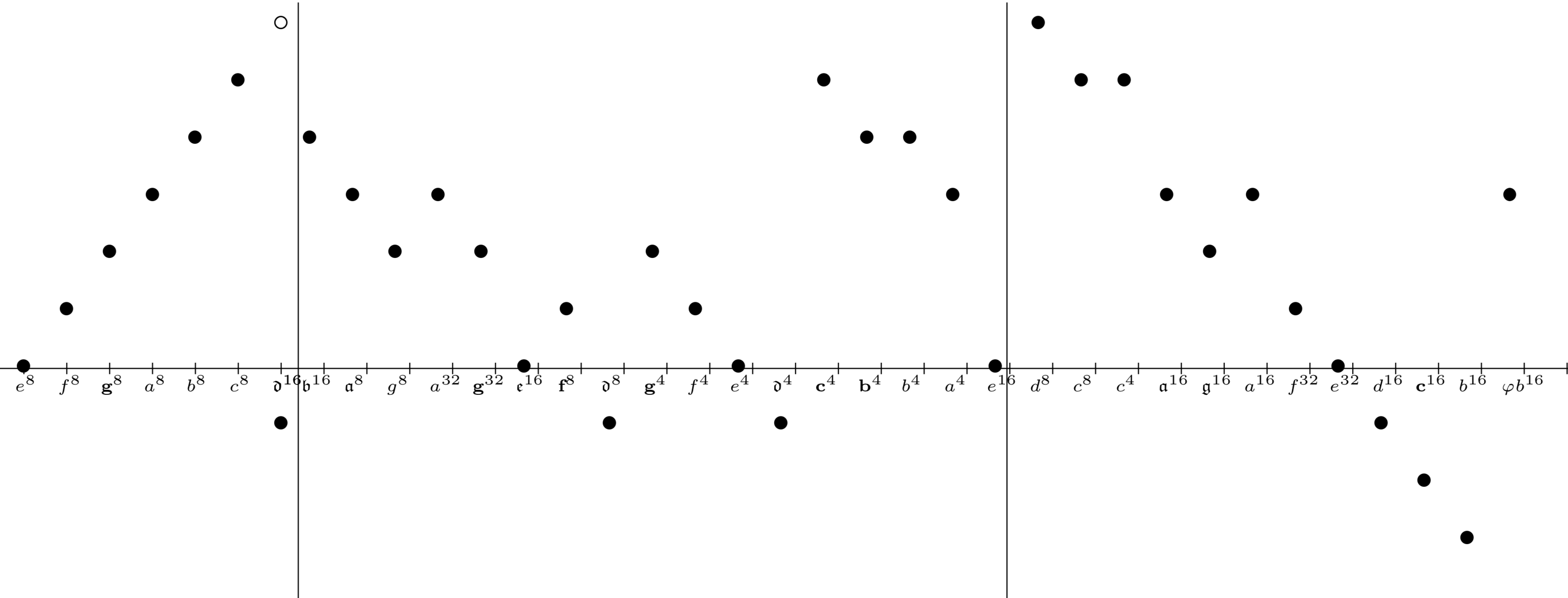}\\
			\includegraphics[scale=0.35]{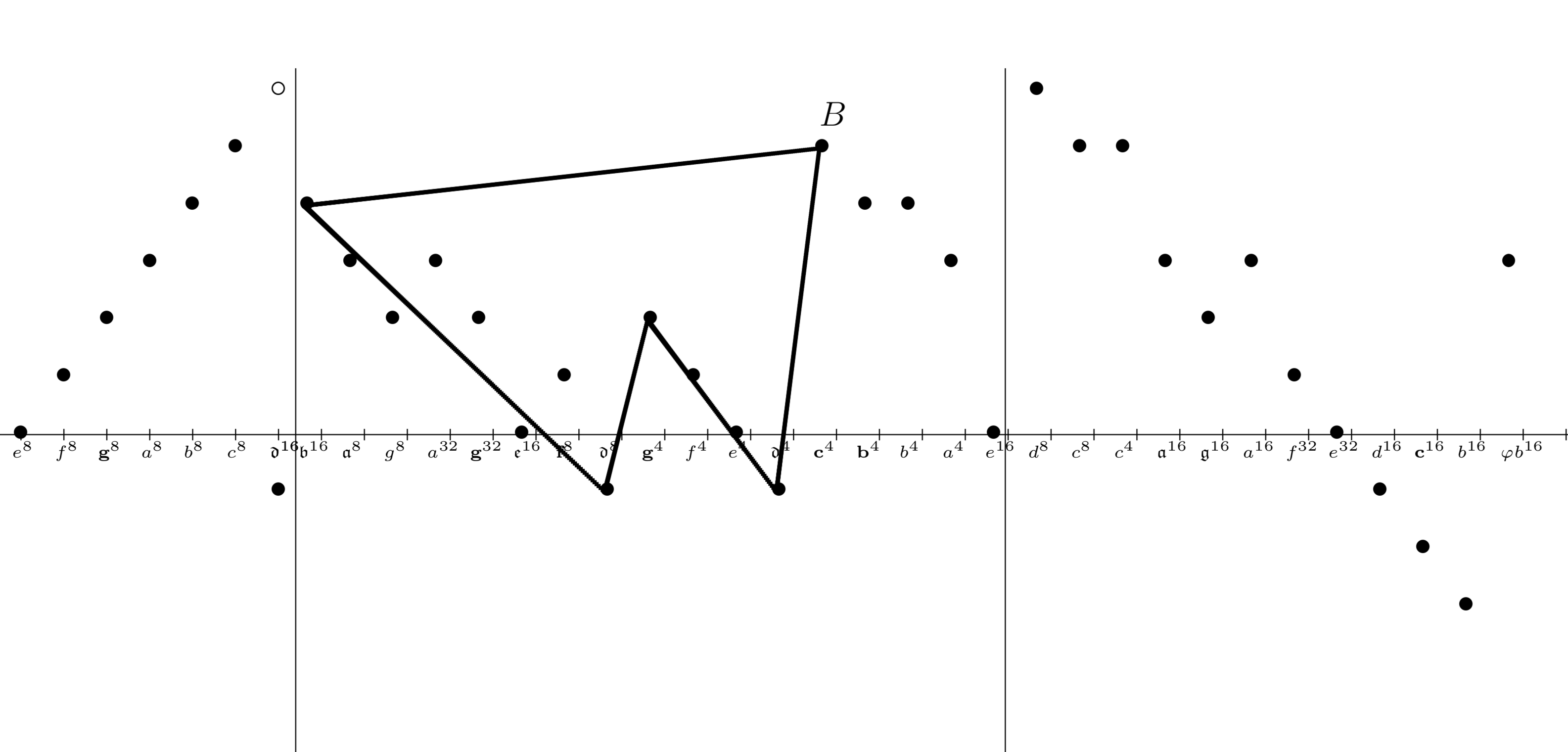}\\
	\includegraphics[scale=0.35]{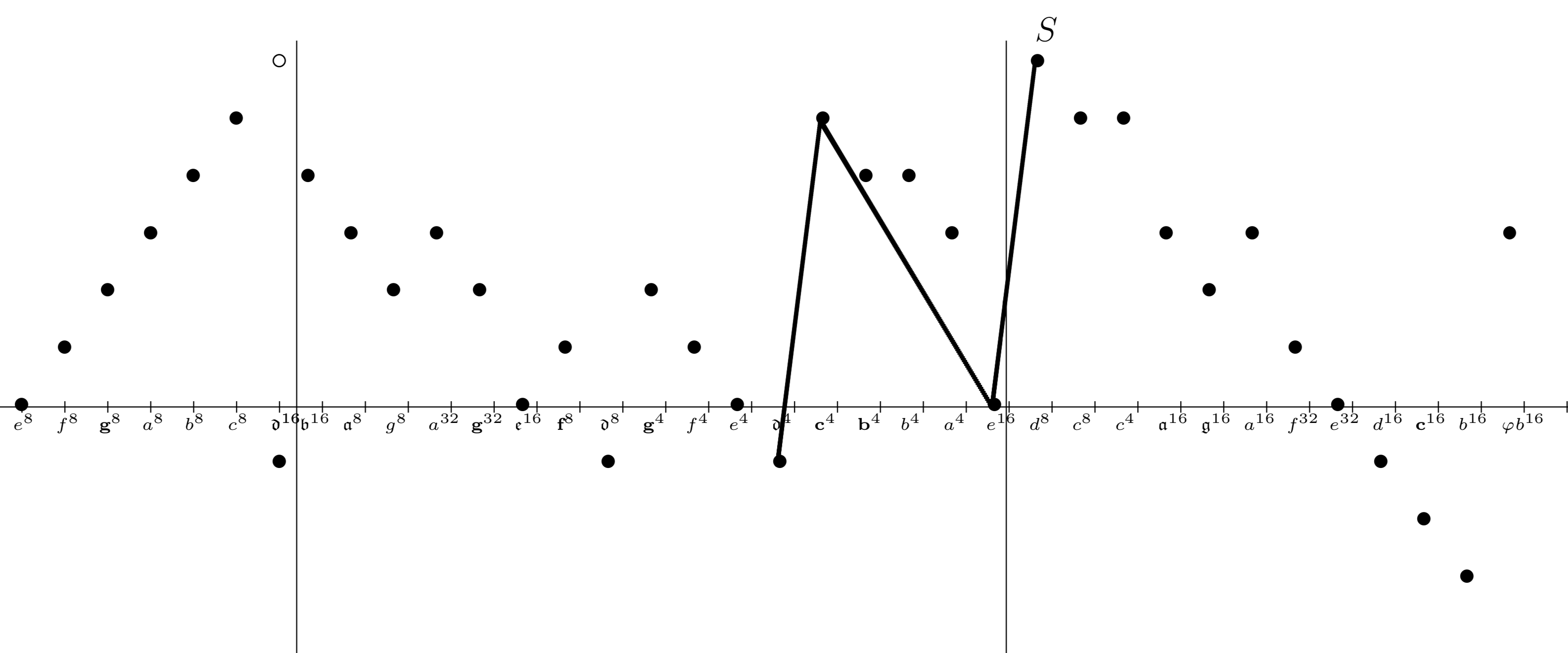}\\
	\includegraphics[scale=0.35]{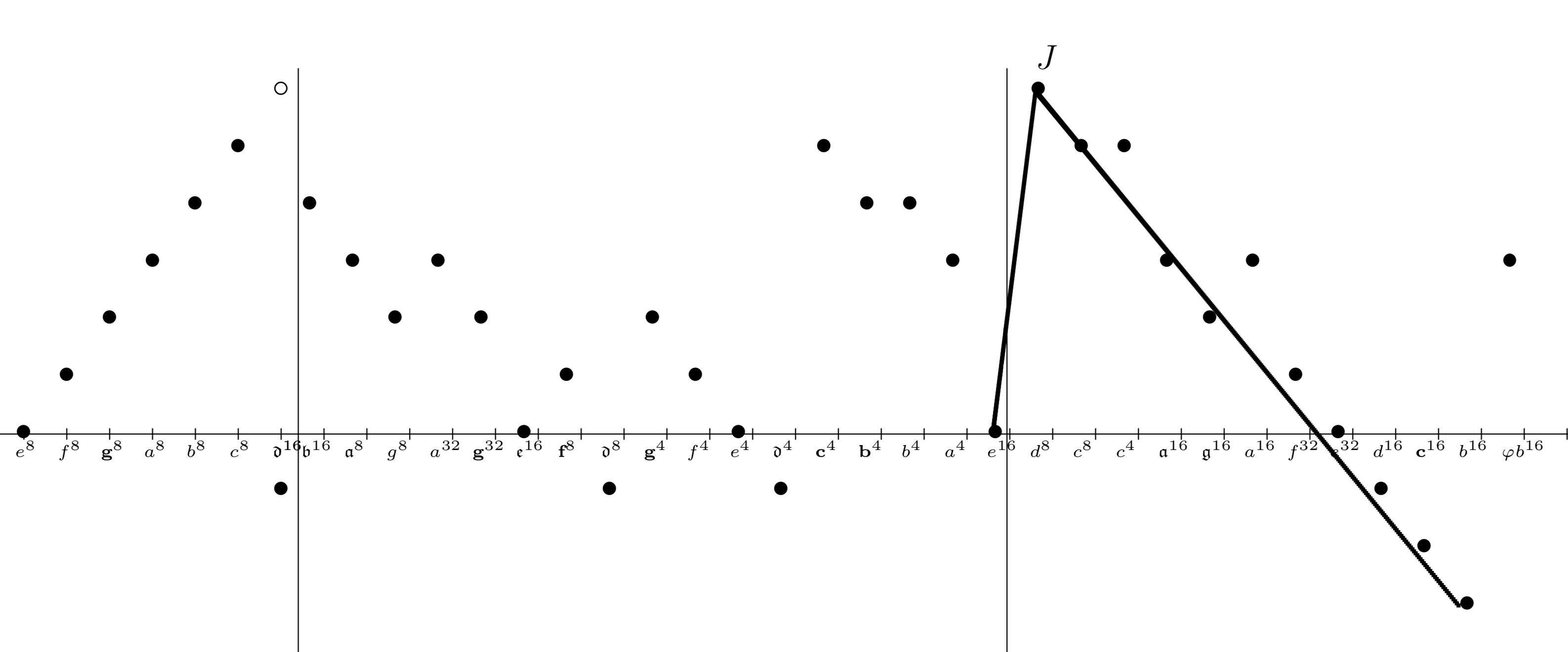}\\
	\includegraphics[scale=0.35]{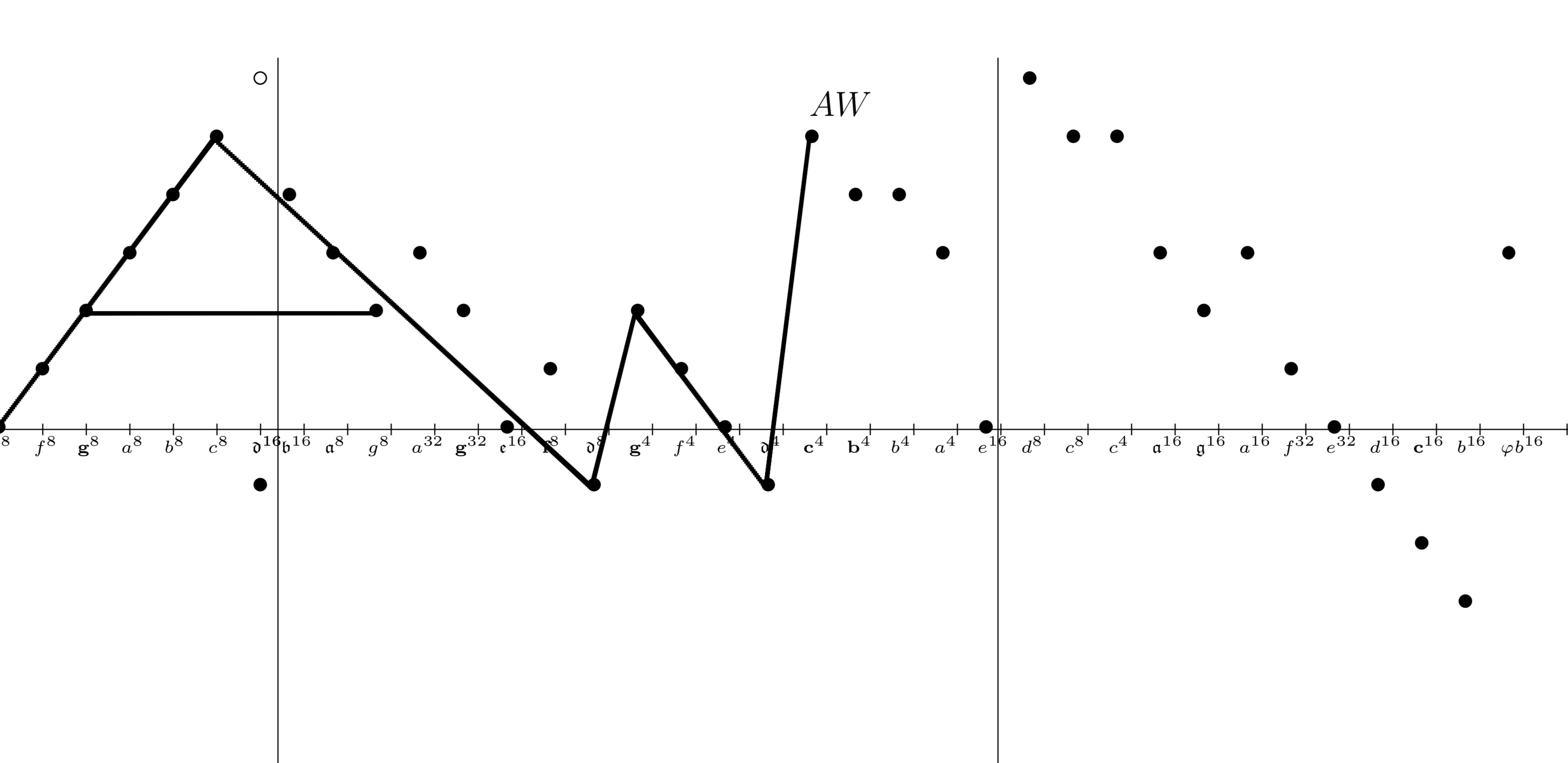}\\
	\caption{Letters J, S, B, A and W also arise from points defined by canon \^a4 Quaerendo Invenietis (see Algorithm \ref{Alg}). }
	\label{QIS1A}
\end{figure}

\begin{figure}[H]
		\centering
		\includegraphics[scale=0.35]{QI-0}\\
			\includegraphics[scale=0.35]{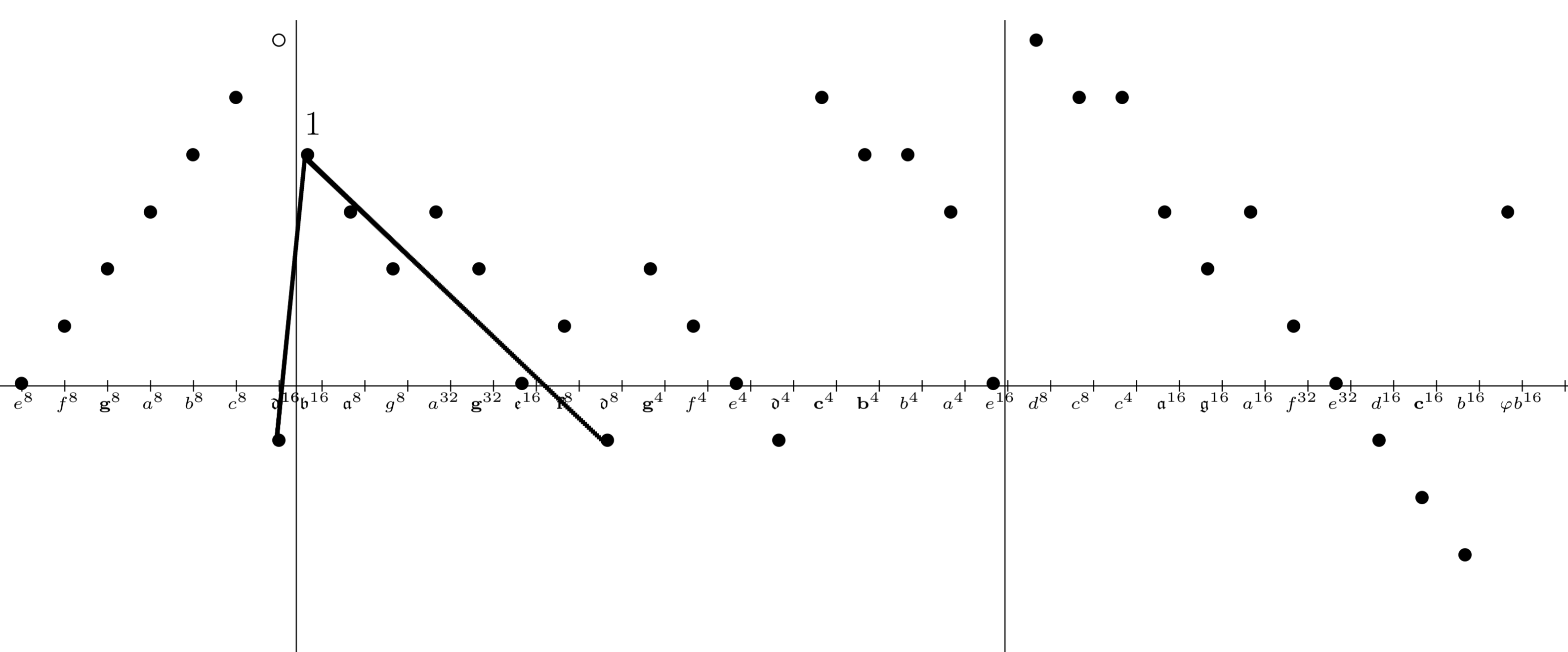}\\
	\includegraphics[scale=0.35]{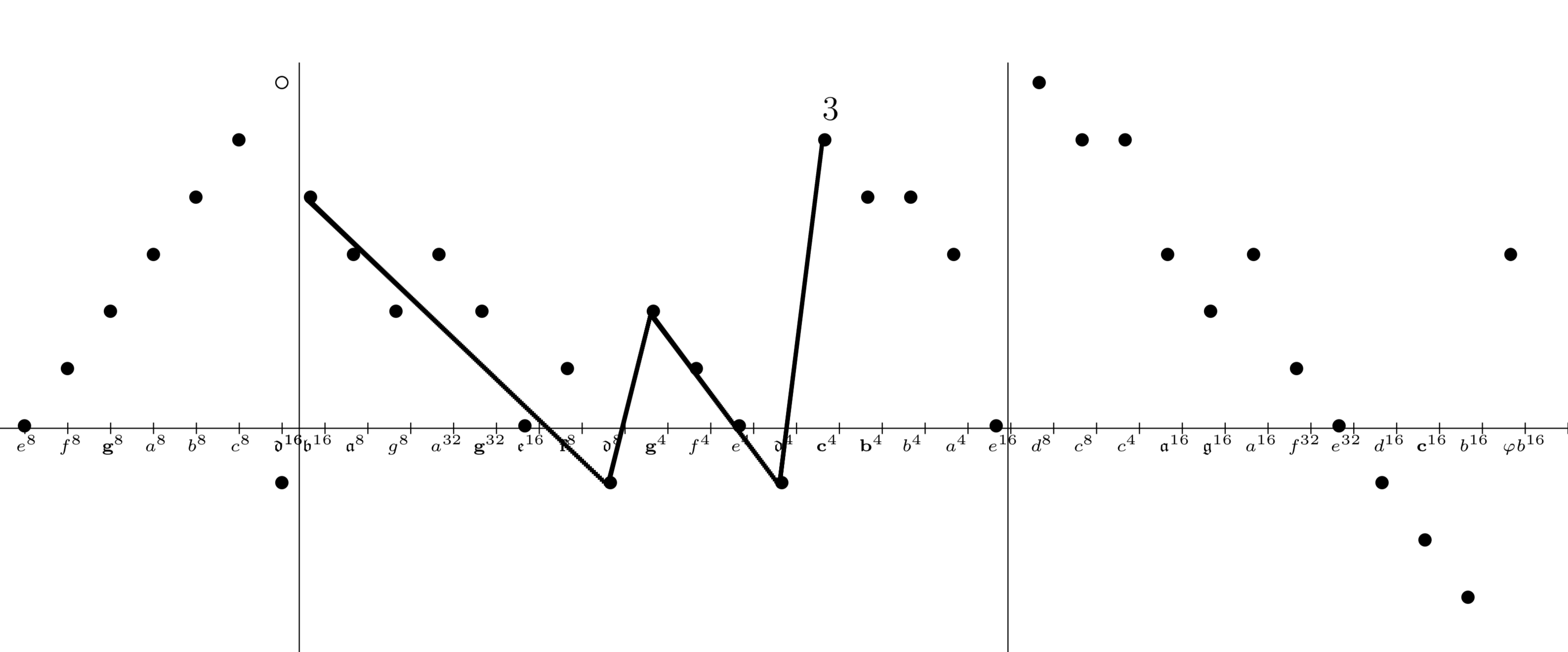}\\
	\includegraphics[scale=0.35]{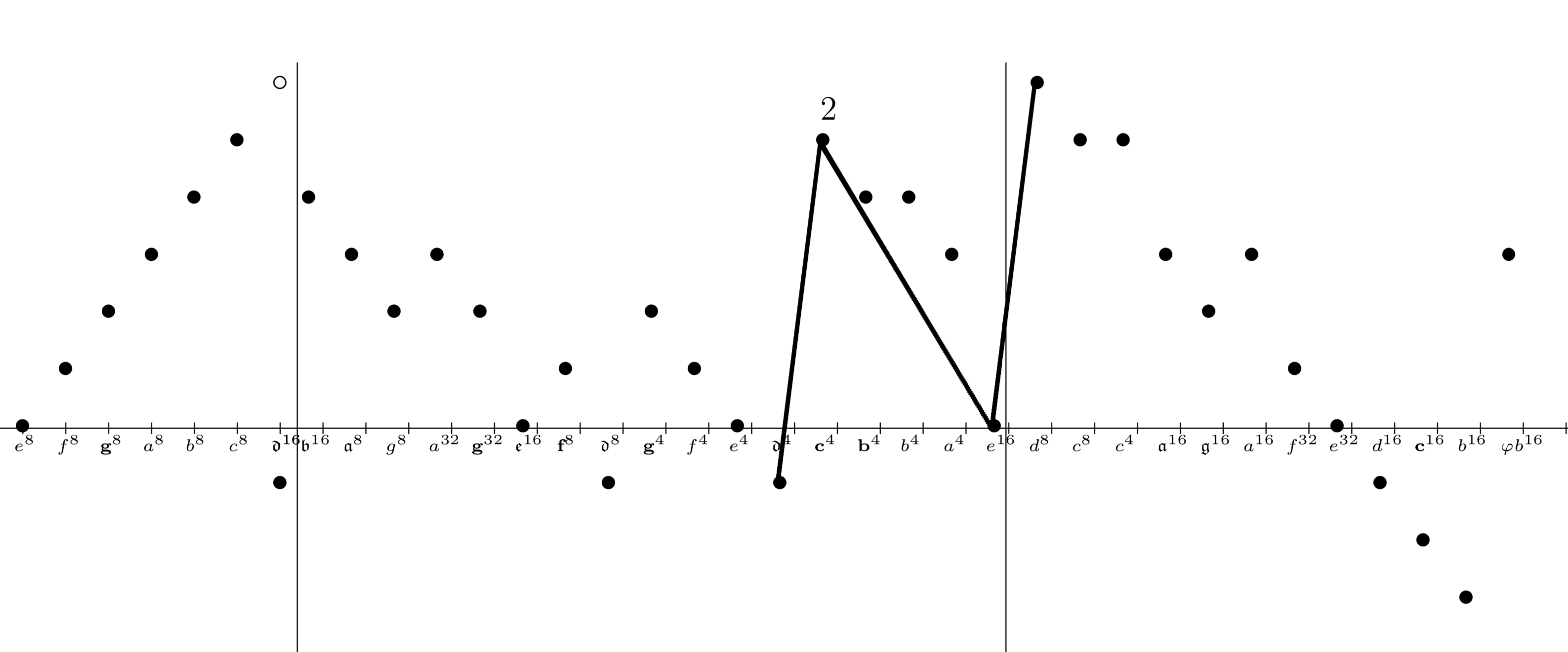}\\
	\includegraphics[scale=0.35]{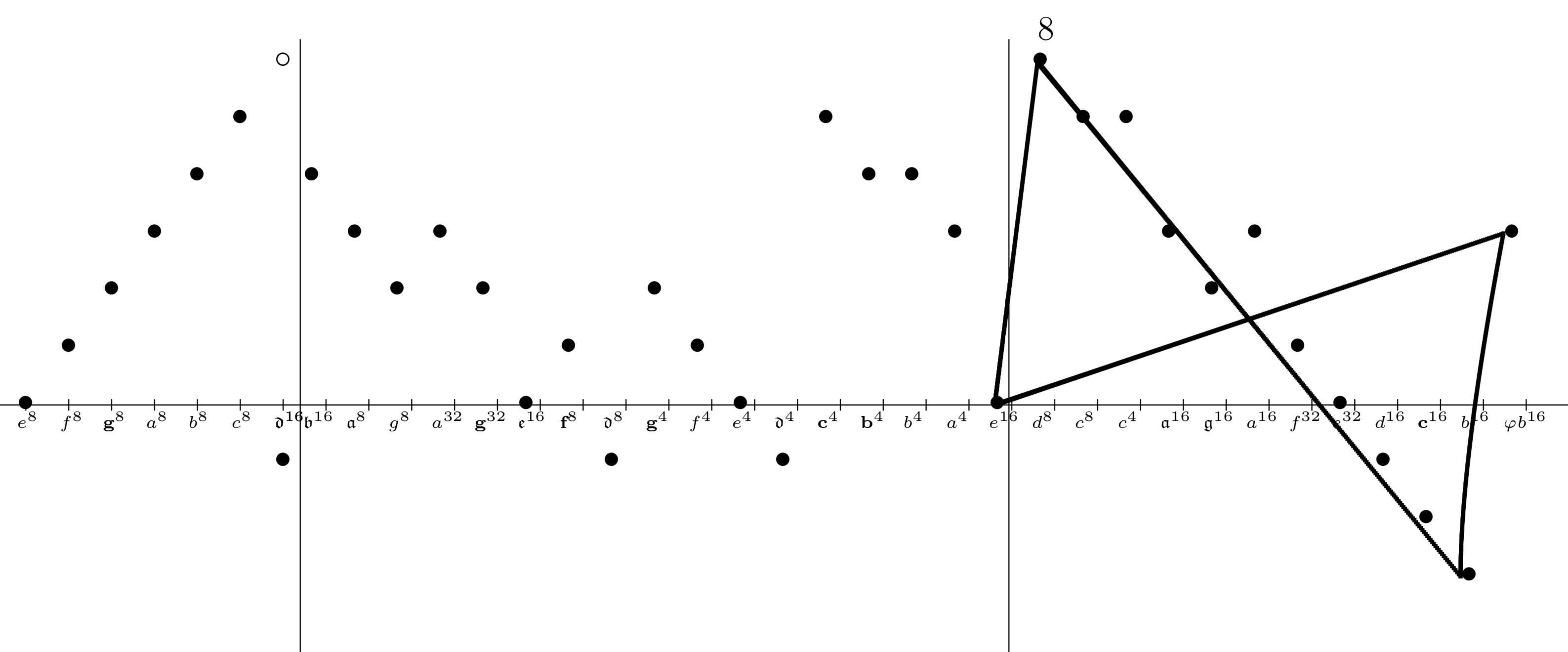}\\
	\caption{Numbers, 1,3,2, and 8 can be built up via canon \^a4 Quaerendo Invenietis. }
	\label{QIS1B}
\end{figure}

  The following points defined by canon \^a4 Quaerendo Invenietis are represented in Figure \ref{QIS1A}. If such points are connected then it is possible to drawn the letters J, S, and B and the letters A and W (see Figure \ref{chrismon(1)}). The scheme allows built numbers 1,2,3, and 8 in Figure \ref{QIS1B}.
  
  \begin{equation}
  \begin{split}
 ((\eighthnote,e),\varepsilon^{0}),&  ((\eighthnote,f),\varepsilon^{1}),((\flat\eighthnote,g),\varepsilon^{2}),((\eighthnote,a),\varepsilon^{3}),((\eighthnote,b),\varepsilon^{4}),((\flat\eighthnote,c),\varepsilon^{5}),((\sharp\Vier,d),\varepsilon^{-1}),\ViPa\\
((\sharp\Vier,b),\varepsilon^{4}),& ((\sharp\eighthnote,a),\varepsilon^{3}),((\natural\eighthnote,g),\varepsilon^{2}), ((\Halb, a), \varepsilon^{3}), ((\flat\Halb, g), \varepsilon^{2}),((\sharp\eighthnote,e),\varepsilon^{0}),((\flat\eighthnote,f),\varepsilon^{1}),\\
((\sharp\eighthnote,d),\varepsilon^{-1}),&((\flat\text{\Sech},g),\varepsilon^{2}),((\text{\Sech},f),\varepsilon^{1}),((\text{\Sech},e),\varepsilon^{0}),((\sharp\text{\Sech},d),\varepsilon^{-1}),((\flat\text{\Sech},c),\varepsilon^{5}),((\text{\Sech},b),\varepsilon^{4}),\\
((\text{\Sech},a),\varepsilon^{3}),& ((\Vier,e),\varepsilon^{0}),\AcPa, ((\eighthnote,d), \varepsilon^{6}),((\natural\eighthnote,c), \varepsilon^{5}),(\natural\text{\Sech},c), \varepsilon^{5}),((\sharp\Vier,a),\varepsilon^{3}), ((\natural\Vier,g), \varepsilon^{2}),\\
((\Vier,a),\varepsilon^{3}),&((\Halb, f), \varepsilon^{1}),((\Halb, e), \varepsilon^{0}),((\Vier,d), \varepsilon^{-1}),((\flat\Vier,c), \varepsilon^{-2}),((\Vier,b), \varepsilon^{-3}),((\underset{\cdot}{\Vier},b),\varepsilon^{11}),\\((\Vier,f), \varepsilon^{1}).&
 \end{split}
 \end{equation}

\vspace{-6pt} 

\section{Concluding Remarks and Future Work}\label{conclusions}

Brauer configuration algebras are a theoretical framework to investigate block cyphers. In particular, Western musical writing can be analyzed via Brauer configuration algebras. Such analysis allows for determining relationships between notes and measures in a musical piece by assuming that these are ciphertexts of classical cryptosystems whose plaintexts are defined by appropriated Brauer configurations. Due to these procedures, a Brauer analysis can be applied to some of the canons proposed by Bach in his Musical Offering. Musical notes in these canons have a structure or form based on some of the commonly used Bach's symbols. \par\bigskip
\begin{centering}
\textit{Future Work}\par\bigskip
\end{centering}

This work applies Brauer's analysis to some of Bach's canons. Analyzing the remaining Bach canons in his Musical Offering is another task to be developed in the future. Proving, in particular, that the structure or form of such canons is also based on Bach's symbols.

\end{document}